\title{Ocneanu Cells and Boltzmann Weights \\ for the $SU(3)$ $\mathcal{ADE}$ Graphs}
\author{
        David E. Evans and Mathew Pugh \\ \\
        School of Mathematics, \\
        Cardiff University, \\
        Senghennydd Road, \\
        Cardiff, CF24 4AG, \\
        Wales, U.K.
}
\date{\today}

\documentclass[12pt]{article}
\usepackage{amssymb}
\usepackage{graphicx}

\newtheorem{Def}{Definition}[section]

\newtheorem{Lemma}[Def]{Lemma}

\newtheorem{Thm}[Def]{Theorem}

\begin{document}
\maketitle

\begin{abstract}
We determine the cells, whose existence has been announced by Ocneanu,  on all the candidate nimrep graphs except $\mathcal{E}_4^{(12)}$ proposed by di Francesco and Zuber for the $SU(3)$ modular invariants classified by Gannon. This enables the Boltzmann weights to be computed for the corresponding integrable statistical mechanical models and provide the framework for studying corresponding braided subfactors to realise all the $SU(3)$ modular invariants as well as a
framework for a new $SU(3)$ planar algebra theory.
\end{abstract}

\maketitle

\section{Introduction}

In the last twenty years, a very fruitful circle of ideas has developed linking the theory of subfactors with modular invariants in conformal field theory. Subfactors have been studied through their paragroups, planar algebras and have serious contact with free probability theory. The understanding and classification of modular invariants is significant for conformal field theory and their underlying statistical mechanical models. These areas are linked through the use of braided subfactors and $\alpha$-induction which in particular for $SU(2)$ subfactors and $SU(2)$ modular invariants invokes $ADE$ classifications on both sides.
This paper is the first of our series to study more precisely these connections in the context of $SU(3)$ subfactors and $SU(3)$ modular invariants. The aim is to understand them not only through braided subfactors and $\alpha$-induction but introduce and develop a pertinent planar algebra theory and free probability.

A group acting on a factor can be recovered from the inclusion of its fixed point algebra.
A general subfactor encodes a more sophisticated symmetry or a way of handling non group like symmetries including but going beyond quantum groups \cite{evans/kawahigashi:1998}.
The classification of subfactors was initiated by Jones \cite{jones:1983} who found that the minimal symmetry to understand the inclusion is through the Temperley-Lieb algebra.
This arises from the representation theory of $SU(2)$ or dually certain representations of Hecke algebras.
All $SU(2)$ modular invariant partition functions were classified by Cappelli, Itzykson and Zuber \cite{cappelli/itzykson/zuber:1987i, cappelli/itzykson/zuber:1987ii} using $ADE$ Coxeter-Dynkin diagrams and their realization by braided subfactors is reviewed and referenced in \cite{evans:2003}. There are a number of invariants (encoding the symmetry) one can assign to a subfactor, and under certain circumstances they are complete at least for hyperfinite subfactors.
Popa \cite{popa:1995} axiomatized the inclusions of relative commutants in the Jones tower, and Jones \cite{jones:planar} showed that this was equivalent to his planar algebra description. Here one is naturally forced to work with nonamenable factors through free probabilistic constructions e.g. \cite{guionnet/jones/shlyakhtenko:2007}. In another vein, Banica and Bisch \cite{banica/bisch:2007} understood the principal graphs, which encode only the multiplicities in the inclusions of the relative commutants, and more generally nimrep graphs in terms of spectral measures, and so provide another way of understanding the subfactor invariants.

In our series of papers we will look at this in the context of $SU(3)$, through the subfactor theory and their modular invariants, beginning here and continuing in \cite{evans/pugh:2009ii, evans/pugh:2009iii, evans/pugh:2009iv, evans/pugh:2009v, evans/pugh:2009vi}. The $SU(3)$ modular invariants were classified by Gannon \cite{gannon:1994}. Ocneanu \cite{ocneanu:2002} announced that all these modular invariants were realised by subfactors, and most of these are understood in the literature and will be reviewed in the sequel \cite{evans/pugh:2009ii}. A braided subfactor automatically gives a modular invariant through $\alpha$-induction.
This $\alpha$-induction yields a representation of the Verlinde algebra or a nimrep - which yields multiplicity graphs associated to the modular invariants (or at least associated to the inclusion, as a modular invariant may be represented by wildly differing inclusions and so may possess inequivalent but isospectral nimreps, as is the case for $\mathcal{E}^{(12)}$). In the case of the $SU(3)$ modular  invariants, candidates of these graphs were proposed by di Francesco and Zuber \cite{di_francesco/zuber:1990} by looking for graphs whose spectrum reproduced the diagonal part of the modular invariant, aided to some degree by first listing the graphs and spectra of fusion graphs of the finite subgroups of $SU(3)$. In the $SU(2)$ situation there is a precise relation between the $ADE$ Coxeter-Dynkin graphs and finite subgroups of $SU(2)$ as part of the McKay correspondence. However, for $SU(3)$, the relation between nimrep graphs and finite subgroups of $SU(3)$ is imprecise and not a perfect match.
For $SU(2)$, an affine Dynkin diagram describing the McKay graph of a finite subgroup gives rise to a Dynkin diagram describing a nimrep or the diagonal part of a modular invariant by removing the vertex corresponding to the identity representation. Di Francesco and Zuber found graphs whose spectrum described the diagonal part of a modular invariant by taking the list of McKay graphs of finite subgroups of $SU(3)$ and removing vertices. Not every modular invariant could be described in this way, and not every finite subgroup yielded a nimrep for a modular invariant. In higher rank $SU(N)$, the number of finite subgroups will increase but the number of exceptional modular invariants should decrease, so this procedure is even less likely to be accurate. Evans and Gannon have suggested an alternative way of associating finite subgroups to modular invariants, by considering the largest finite stabiliser groups \cite{evans/gannon:2008}.

A modular invariant which is realised by a subfactor will yield a graph. To construct these subfactors we will need some input graphs which will actually coincide with the output nimrep graphs - $SU(3)$ $\mathcal{ADE}$ graphs.
The aim of this series of papers is to study the $SU(3)$ $\mathcal{ADE}$ graphs, which appear in the classification of modular invariant partition functions from numerous viewpoints including the determination of their Boltzmann weights in this paper, representations of $SU(3)$-Temperley-Lieb or Hecke algebra \cite{evans/pugh:2009ii}, a new notion of $SU(3)$-planar algebras \cite{evans/pugh:2009iii} and their modules \cite{evans/pugh:2009iv}, and spectral measures \cite{evans/pugh:2009v, evans/pugh:2009vi}.

As pointed out to us by Jean-Bernard Zuber, there is a renewal of interest (by physicists) in these $SU(3)$ and related theories, in connection with topological quantum computing \cite{ardonne/schoutens:2007} and by Joost Slingerland in connection with condensed matter physics \cite{bais/slingerland:2008} where we see that $\alpha$-induction is playing a key role.

We begin however in this paper by computing the numerical values of the Ocneanu cells, and consequently representations of the Hecke algebra, for the $\mathcal{ADE}$ graphs. These cells give numerical weight to Kuperberg's \cite{kuperberg:1996} diagram of trivalent vertices -- corresponding to the fact that the trivial representation is contained in the triple product of the fundamental representation of $SU(3)$ through the determinant. They will yield in a natural way, representations of an $SU(3)$-Temperley-Lieb or Hecke algebra.
(For $SU(2)$ or bipartite graphs, the corresponding weights (associated to the diagrams of cups or caps),
arise in a more straightforward fashion from a Perron-Frobenius eigenvector, giving a natural representation of the Temperley-Lieb algebra or Hecke algebra).
We have been unable thus far to compute the cells for the exceptional graph $\mathcal{E}_4^{(12)}$. This graph is meant to be the nimrep for the modular invariant conjugate to the Moore-Seiberg invariant $\mathcal{E}_{MS}^{(12)}$ \cite{moore/seiberg:1989}. However we will still be able to realise this modular invariant by subfactors in \cite{evans/pugh:2009ii} using \cite{evans/pinto:2003}.
For the orbifold graphs $\mathcal{D}^{(3k)}$, $k = 2,3,\ldots \;$, orbifold conjugate $\mathcal{D}^{(n)\ast}$, $n=6,7,\ldots \;$, and $\mathcal{E}_1^{(12)}$ we compute solutions which satisfy some additional condition, but for the other graphs we compute all the Ocneanu cells, up to equivalence. The existence of these cells has been announced by Ocneanu (e.g. \cite{ocneanu:2000ii, ocneanu:2002}), although the numerical values have remained unpublished. Some of the representations of the Hecke algebra have appeared in the literature and we compare our results.

For the $\mathcal{A}$ graphs, our solution for the Ocneanu cells $W$ gives an identical representation of the Hecke algebra to that of Jimbo et al. \cite{jimbo/miwa/okado:1987} given in (\ref{eqn:Boltzmann-Weights_for_A}). Our cells for the $\mathcal{A}^{(n)\ast}$ graphs give equivalent Boltzmann weights to those given by Behrend and Evans in \cite{behrend/evans:2004}. In \cite{di_francesco/zuber:1990}, di Francesco and Zuber give a representation of the Hecke algebra for the graphs $\mathcal{D}^{(6)\ast}$ and $\mathcal{E}^{(8)}$, whilst in \cite{sochen:1991} a representation of the Hecke algebra is computed for the graphs $\mathcal{E}_1^{(12)}$ and $\mathcal{E}^{(24)}$. Our solutions for the cells $W$ give an identical Hecke representation for $\mathcal{E}^{(8)}$ and an equivalent Hecke representation for $\mathcal{E}_1^{(12)}$. However, for $\mathcal{E}^{(24)}$, our cells give inequivalent Boltzmann weights. In \cite{fendley:1989}, Fendley gives Boltzmann weights for $\mathcal{D}^{(6)}$ which are not equivalent to those we obtain, but which are equivalent if we take one of the weights in \cite{fendley:1989} to be the complex conjugate of what is given.

Subsequently, we will use these weights, their existence and occasionally more precise use of their numerical values. Here we outline some of the flavour of these applications.
We use these cells to define an $SU(3)$ analogue of the Goodman-de la Harpe Jones construction of a subfactor, where we embed the $SU(3)$-Temperley-Lieb or Hecke algebra in an AF path algebra of the $SU(3)$ $\mathcal{ADE}$ graphs. Using this construction, we realize all the $SU(3)$ modular invariants by subfactors \cite{evans/pugh:2009ii}.

We will then \cite{evans/pugh:2009iii, evans/pugh:2009iv} look at the $SU(3)$-Temperley-Lieb algebra and the $SU(3)$-GHJ subfactors from the viewpoint of planar algebras. We give a diagrammatic representation of the $SU(3)$-Temperley-Lieb algebra, and show that it is isomorphic to Wenzl's representation of a Hecke algebra. Generalizing Jones's notion of a planar algebra, we construct an $SU(3)$-planar algebra which will capture the structure contained in the $SU(3)$ $\mathcal{ADE}$ subfactors. We show that the subfactor for an $\mathcal{ADE}$ graph with a flat connection has a description as a flat $SU(3)$-planar algebra. We introduce the notion of modules over an $SU(3)$-planar algebra, and describe certain irreducible Hilbert $SU(3)$-$TL$-modules. A partial decomposition of the $SU(3)$-planar algebras for the $\mathcal{ADE}$ graphs is achieved.
Moreover, in \cite{evans/pugh:2009v, evans/pugh:2009vi} we consider spectral measures for the $ADE$ graphs in terms of probability measures on the circle $\mathbb{T}$. We generalize this to $SU(3)$, and in particular obtain spectral measures for the $SU(3)$ graphs. We also compare various Hilbert series of dimensions associated to $ADE$ models for $SU(2)$, and compute the Hilbert series of certain $q$-deformed Calabi-Yau algebras of dimension 3.

In Section \ref{section2}, we specify the graphs we are interested in, and in Section \ref{section3} recall the notion of cells due to Ocneanu which we will then compute in Sections \ref{Sect:computation_of_weights} - \ref{Sect:computation_of_weights-E(24)}.

\section{$\mathcal{ADE}$ Graphs}\label{section2}

We enumerate the graphs we are interested in. These will eventually provide the nimrep classification graphs for the list of $SU(3)$ modular invariants, but at this point, they will only provide a framework for
some statistical mechanical models with configurations spaces built from these graphs together with some Boltzmann weights which we will need to construct.
However, for the sake of clarity of notation, we start by listing the $SU(3)$ modular invariants. There are four infinite series of $SU(3)$ modular invariants: the diagonal invariants, labelled by $\mathcal{A}$, the orbifold invariants $\mathcal{D}$, the conjugate invariants $\mathcal{A}^{\ast}$, and the orbifold conjugate invariants $\mathcal{D}^{\ast}$. These will provide four infinite  families of graphs, written  as
$\mathcal{A}$, the orbifold graphs $\mathcal{D}$, the conjugate graphs $\mathcal{A}^{\ast}$, and the orbifold conjugate graphs $\mathcal{D}^{\ast}$, shown in Figures \ref{Fig:SU(3)-A(infty)}, \ref{fig:Weights3}, \ref{fig:Weights10}, \ref{fig:A(star)-graphs} and \ref{fig:Weights11}.
There are also exceptional $SU(3)$ modular invariants, i.e. invariants which are not diagonal, orbifold, or their conjugates, and there are only finitely many of these. These are $\mathcal{E}^{(8)}$ and its conjugate, $\mathcal{E}^{(12)}$, $\mathcal{E}_{MS}^{(12)}$ and its conjugate, and $\mathcal{E}^{(24)}$. The exceptional invariants $\mathcal{E}^{(12)}$ and $\mathcal{E}^{(24)}$ are self-conjugate.

The modular invariants arising from $SU(3)_k$ conformal embeddings are:
\begin{itemize}
\item $\mathcal{D}^{(6)}$: $SU(3)_3 \subset SO(8)_1$, also realised as an orbifold $SU(3)_3 / \mathbb{Z}_3$,
\item $\mathcal{E}^{(8)}$: $SU(3)_5 \subset SU(6)_1$, plus its conjugate,
\item $\mathcal{E}^{(12)}$: $SU(3)_9 \subset (\mathrm{E}_6)_1$,
\item $\mathcal{E}_{MS}^{(12)}$: Moore-Seiberg invariant, an automorphism of the orbifold invariant $\mathcal{D}^{(12)} = SU(3)_9 / \mathbb{Z}_3$, plus its conjugate,
\item $\mathcal{E}^{(24)}$: $SU(3)_{21} \subset (\mathrm{E}_7)_1$.
\end{itemize}

These modular invariants will be associated with graphs, as follows. There will be one graph $\mathcal{E}^{(8)}$ for the $\mathcal{E}^{(8)}$ modular invariant and its orbifold graph $\mathcal{E}^{(8)\ast}$ for its conjugate invariant as in Figure \ref{fig:E(8)&E(8)(star)-graphs}.
The modular invariants $\mathcal{E}_{MS}^{(12)}$ and its conjugate will be associated to the graphs $\mathcal{E}_{5}^{(12)}$ and $\mathcal{E}_{4}^{(12)}$ respectively as in Figure \ref{fig:Weights7}.
The exceptional invariant $\mathcal{E}^{(12)}$ is self-conjugate but has associated to it two isospectral graphs $\mathcal{E}_1^{(12)}$  and $\mathcal{E}_2^{(12)}$ as in Figure \ref{fig:labelled_E1(12)&E2(12)graphs}. The invariant $\mathcal{E}^{(24)}$ is also self-conjugate and has associated to it one graph $\mathcal{E}^{(24)}$ as in Figure \ref{fig:Weights8}. The modular invariants themselves play no role in this paper other than to help label these graphs. In the sequel to this paper \cite{evans/pugh:2009ii} we will use the Boltzmann weights obtained here to construct braided subfactors,
which via $\alpha$-induction \cite{bockenhauer/evans:1998, bockenhauer/evans:1999i, bockenhauer/evans:1999ii, bockenhauer/evans:2000, bockenhauer/evans/kawahigashi:1999, bockenhauer/evans/kawahigashi:2000} will realise the corresponding modular invariants. Furthermore, $\alpha$-induction naturally provides a nimrep or representation of the original fusion rules or Verlinde algebra. The corresponding nimreps will then be computed and we will recover the original input graph. The theory of $\alpha$-induction will guarantee that the spectra of these graphs are described by the diagonal part of the corresponding modular invariant. Thus detailed information about the spectra of these graphs will naturally follow from this procedure and does not need to be computed at this stage.
Many of these modular invariants are already realised in the literature and this will be reviewed in the sequel to this paper \cite{evans/pugh:2009ii}.

\section{Ocneanu Cells}\label{section3}

Let $\Gamma$ be $SU(3)$ and $\widehat{\Gamma}$ its irreducible representations. One can associate to $\Gamma$ a McKay graph $\mathcal{G}_{\Gamma}$ whose vertices are labelled by the irreducible representations of $\Gamma$, where for any pair of vertices $i, j \in \widehat{\Gamma}$ the number of edges from $i$ to $j$ are given by the multiplicity of $j$ in the decomposition of $i \otimes \rho$ into irreducible representations, where $\rho$ is the fundamental irreducible representation of $SU(3)$, and which, along with its conjugate representation $\overline{\rho}$, generates $\widehat{\Gamma}$.
The graph $\mathcal{G}_{\Gamma}$ is made of triangles, corresponding to the fact that the fundamental representation $\rho$ satisfies $\rho \otimes \rho \otimes \rho \ni \mathbf{1}$. We define maps $s$, $r$ from the edges of $\mathcal{G}_{\Gamma}$ to its vertices, where for an edge $\gamma$, $s(\gamma)$ denotes the source vertex of $\gamma$ and $r(\gamma)$ its range vertex. For the graph $\mathcal{G}_{\Gamma}$, a triangle $\triangle_{ijk}^{(\alpha \beta \gamma)} = i \stackrel{\alpha}{\longrightarrow} j \stackrel{\beta}{\longrightarrow} k \stackrel{\gamma}{\longrightarrow} i$ is a closed path of length 3 consisting of edges $\alpha$, $\beta$, $\gamma$ of $\mathcal{G}_{\Gamma}$ such that $s(\alpha) = r(\gamma) = i$, $s(\beta) = r(\alpha) = j$ and $s(\gamma) = r(\beta) = k$. For each triangle $\triangle_{ijk}^{(\alpha \beta \gamma)}$, the maps $\alpha$, $\beta$ and $\gamma$ are composed:
$$i \stackrel{\mathrm{id} \otimes \mathrm{det}^{\ast}}{\longrightarrow} i \otimes \rho \otimes \rho \otimes \rho \stackrel{\gamma \otimes \mathrm{id}}{\longrightarrow} k \otimes \rho \otimes \rho \stackrel{\beta \otimes \mathrm{id}}{\longrightarrow} j \otimes \rho \stackrel{\alpha \otimes \mathrm{id}}{\longrightarrow} i,$$
and since $i$ is irreducible, the composition is a scalar. Then for every such triangle on $\mathcal{G}_{\Gamma}$ there is a complex number, called an Ocneanu cell. There is a gauge freedom on the cells, which comes from a unitary change of basis in $\mathrm{Hom}[i \otimes \rho, j]$ for every pair $i$, $j$.

These cells are axiomatized in the context of an arbitrary graph $\mathcal{G}$ whose adjacency matrix has Perron-Frobenius eigenvalue $[3] = [3]_q$, although in practice it will be one of the $\mathcal{ADE}$ graphs. Note however we do not require $\mathcal{G}$ to be three-colourable (e.g. the graphs $\mathcal{A}^{\ast}$ which will be associated to the conjugate modular invariant).
Here the quantum number $[m]_q$ be defined as $[m]_q = (q^m - q^{-m})/(q - q^{-1})$. We will frequently denote the quantum number $[m]_q$ simply by $[m]$, for $m \in \mathbb{N}$. Now $[3]_q = q^2 + 1 + q^{-2}$, so that $q$ is easily determined from the eigenvalue of $\mathcal{G}$. The quantum number $[2] = [2]_q$ is then simply $q + q^{-1}$.
If $\mathcal{G}$ is an $\mathcal{ADE}$ graph, the Coxeter number $n$ of $\mathcal{G}$ is the number in parentheses in the notation for the graph $\mathcal{G}$, e.g. the exceptional graph $\mathcal{E}^{(8)}$ has Coxeter number 8, and $q = e^{\pi i/n}$. With this $q$, the quantum numbers $[m]_q$ satisfy the fusion rules for the irreducible representations of the quantum group $SU(2)_n$, i.e.
\begin{equation}\label{eqn:fusion_rules_for_quantum_numbers}
[a]_q \; [b]_q = \sum_c \; [c]_q,
\end{equation}
where the summation is over all integers $|b-a| \leq c \leq \mathrm{min}(a+b,2n-a-b)$ such that $a+b+c$ is even.

We define a type I frame in an arbitrary $\mathcal{G}$ to be a pair of edges $\alpha$, $\alpha'$ which have the same start and endpoint. A type II frame will be given by four edges $\alpha_i$, $i=1,2,3,4$, such that $s(\alpha_1) = s(\alpha_4)$, $s(\alpha_2) = s(\alpha_3)$, $r(\alpha_1) = r(\alpha_2)$ and $r(\alpha_3) = r(\alpha_4)$.

\begin{Def}[\cite{ocneanu:2002}]\label{cell_system}
Let $\mathcal{G}$ be an arbitrary graph with Perron-Frobenius eigenvalue $[3]$ and Perron-Frobenius eigenvector $(\phi_i)$. A \textbf{cell system} $W$ on $\mathcal{G}$ is a map that associates to each oriented triangle $\triangle_{ijk}^{(\alpha \beta \gamma)}$ in $\mathcal{G}$ a complex number $W \left( \triangle_{ijk}^{(\alpha \beta \gamma)} \right)$ with the following properties:

$(i)$ for any type I frame \includegraphics[width=16mm]{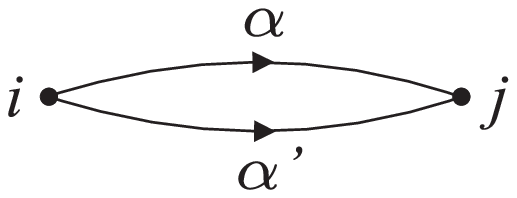} in $\mathcal{G}$ we have \\
\begin{minipage}[b]{11.5cm}
 \begin{minipage}[t]{2cm}
  \mbox{} \\
  \parbox[t]{1cm}{\begin{eqnarray}\label{eqn:typeI_frame}\end{eqnarray}}
 \end{minipage}
 \begin{minipage}[t]{9cm}
  \mbox{} \\
 \includegraphics[width=90mm]{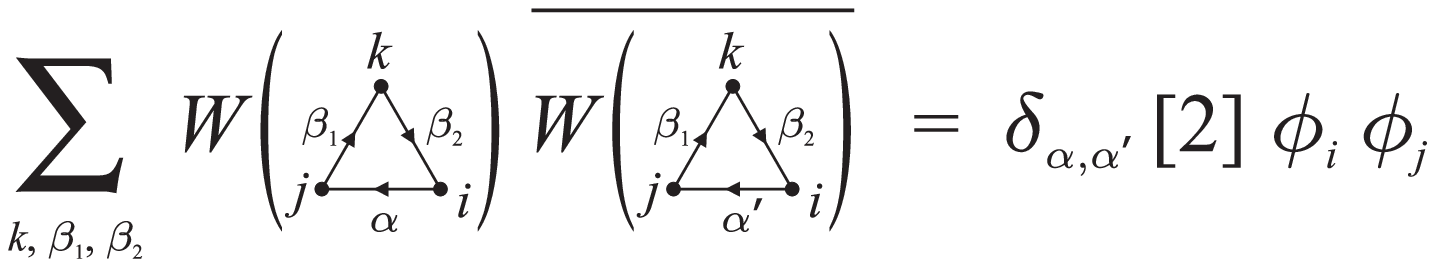}
 \end{minipage}
\end{minipage} \\

$(ii)$ for any type II frame \includegraphics[width=30mm]{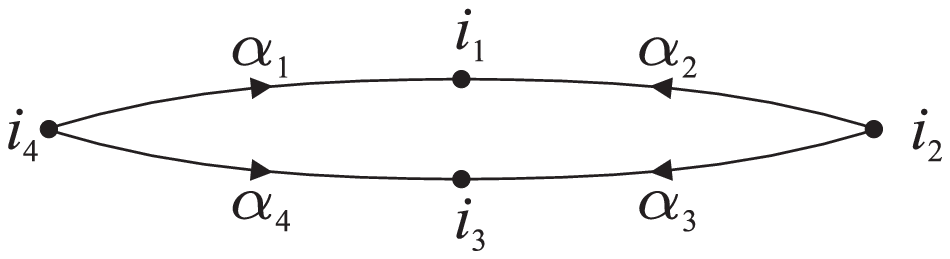} in $\mathcal{G}$ we have \\
\begin{minipage}[b]{11.5cm}
 \begin{minipage}[b]{2cm}
  \mbox{} \\
  \parbox[b]{1cm}{\begin{eqnarray}\label{eqn:typeII_frame}\end{eqnarray}}
 \end{minipage}
 \begin{minipage}[b]{9cm}
  \mbox{} \\
 \includegraphics[width=90mm]{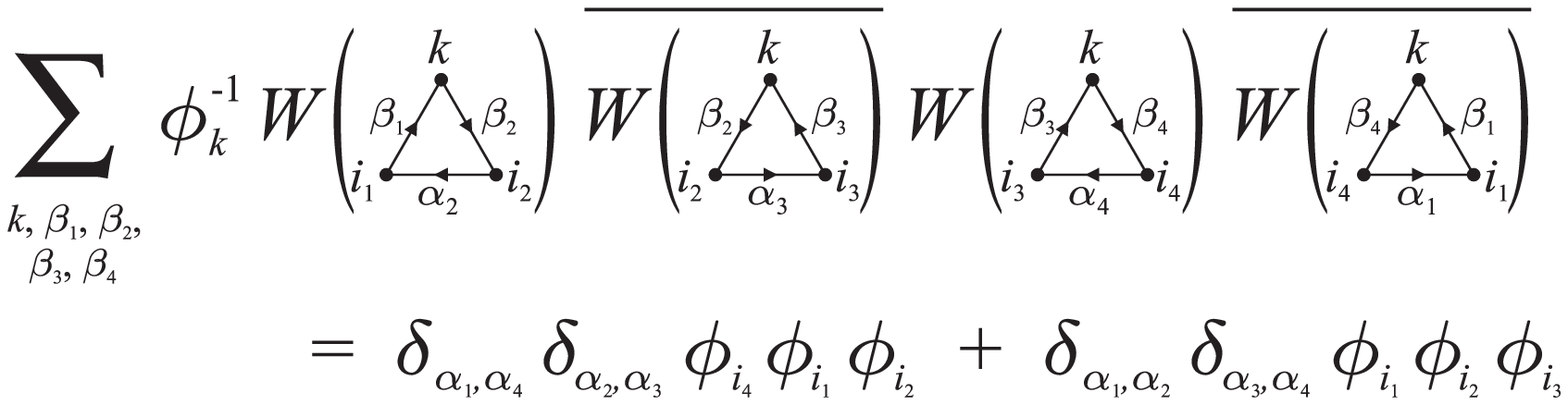}
 \end{minipage}
\end{minipage} \\
\end{Def}

\begin{figure}[bt]
\begin{center}
\includegraphics[width=40mm]{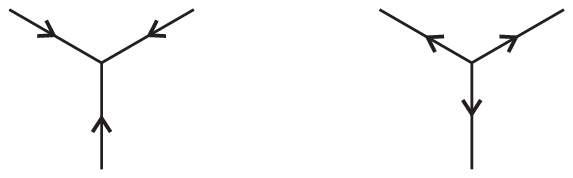}\\
 \caption{$A_2$ webs}\label{fig:A2-webs1}
\end{center}
\end{figure}

In \cite{kuperberg:1996}, Kuperberg defined the notion of a spider- a way of depicting the operations of the representation theory of groups and other group-like objects with certain planar graphs, called webs (hence the term ``spider''). Certain spiders were defined in terms of generators and relations, isomorphic to the representation theories of rank two Lie algebras and the quantum deformations of these representation theories. This formulation generalized a well-known construction for $A_1 = \textrm{su}(2)$ by Kauffman \cite{kauffman:1987}. For the $A_2 = \textrm{su}(3)$ case, the $A_2$ webs are illustrated in Figure \ref{fig:A2-webs1}.

The $A_2$ web space generated by these $A_2$ webs satisfy the Kuperberg relations, which are relations on local parts of the diagrams:
\begin{center}
\begin{minipage}[b]{11.5cm}
 \begin{minipage}[t]{3cm}
  \parbox[t]{2cm}{\begin{eqnarray*}\textrm{K1:}\end{eqnarray*}}
 \end{minipage}
 \begin{minipage}[t]{5.5cm}
  \mbox{} \\
 \includegraphics[width=20mm]{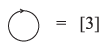}
 \end{minipage}
 \begin{minipage}[t]{2cm}
  \mbox{} \\
  \parbox[t]{1cm}{}
 \end{minipage}
\end{minipage}
\begin{minipage}[b]{11.5cm}
 \begin{minipage}[t]{3cm}
  \parbox[t]{2cm}{\begin{eqnarray*}\textrm{K2:}\end{eqnarray*}}
 \end{minipage}
 \begin{minipage}[t]{5.5cm}
  \mbox{} \\
 \includegraphics[width=23mm]{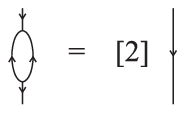}
 \end{minipage}
 \begin{minipage}[t]{2cm}
  \mbox{} \\
  \parbox[t]{1cm}{}
 \end{minipage}
\end{minipage}
\begin{minipage}[b]{11.5cm}
 \begin{minipage}[t]{3cm}
  \parbox[t]{2cm}{\begin{eqnarray*}\textrm{K3:}\end{eqnarray*}}
 \end{minipage}
 \begin{minipage}[t]{5.5cm}
  \mbox{} \\
 \includegraphics[width=55mm]{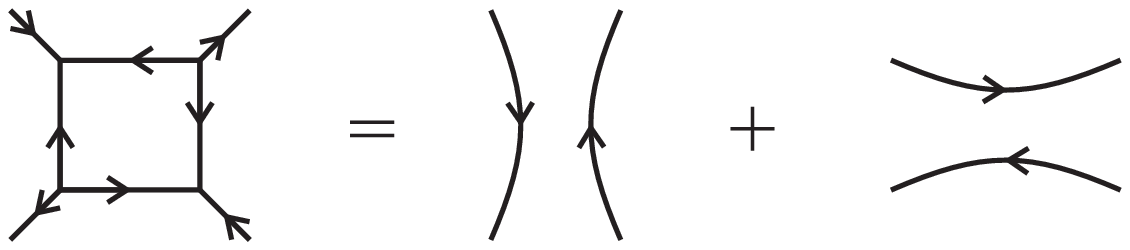}
 \end{minipage}
 \begin{minipage}[t]{2cm}
  \mbox{} \\
  \parbox[t]{1cm}{}
 \end{minipage}
\end{minipage}
\end{center}

\begin{figure}[b]
\begin{center}
\includegraphics[width=100mm]{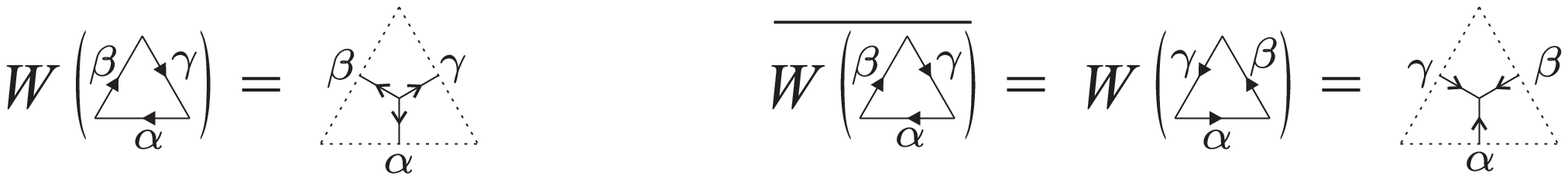}\\
 \caption{Cells associated to trivalent vertices} \label{fig:Oc-Kup}
\end{center}
\end{figure}

The rules (\ref{eqn:typeI_frame}), (\ref{eqn:typeII_frame}) correspond precisely to evaluating the Kuperberg relations K2, K3 respectively, associating a cell $W(\triangle_{\alpha, \beta, \gamma})$ to an incoming trivalent vertex, and $\overline{W(\triangle_{\alpha, \beta, \gamma})}$ to an outgoing trivalent vertex, as in Figure \ref{fig:Oc-Kup}.

We define the connection
$$X^{\rho_1,\rho_2}_{\rho_3,\rho_4} =
\begin{array}{c}
l \; \stackrel{\rho_1}{\longrightarrow} \; i \\
\scriptstyle \rho_3 \textstyle \big\downarrow \qquad \big\downarrow \scriptstyle \rho_2 \\
k \; \stackrel{\textstyle\longrightarrow}{\scriptstyle{\rho_4}} \; j
\end{array}$$
for $\mathcal{G}$ by
\begin{equation} \label{Connection_using_weights_W}
X^{\rho_1,\rho_2}_{\rho_3,\rho_4} = q^{\frac{2}{3}} \delta_{\rho_1,\rho_3} \delta_{\rho_2, \rho_4} - q^{-\frac{1}{3}} \mathcal{U}^{\rho_1,\rho_2}_{\rho_3,\rho_4},
\end{equation}
where $\mathcal{U}^{\rho_1,\rho_2}_{\rho_3,\rho_4}$ is given by the representation of the Hecke algebra, and is defined by
\begin{equation} \label{eqn:HeckeRep}
\mathcal{U}^{\rho_1,\rho_2}_{\rho_3,\rho_4} = \sum_{\lambda} \phi_{s(\rho_1)}^{-1} \phi_{r(\rho_2)}^{-1} W(\triangle_{j,l,k}^{(\lambda, \rho_3, \rho_4)}) \overline{W(\triangle_{j,l,i}^{(\lambda, \rho_1, \rho_2)})}.
\end{equation}
This definition of the connection is really Kuperberg's braiding of \cite{kuperberg:1996}.

\begin{figure}[tb]
\begin{center}
\includegraphics[width=70mm]{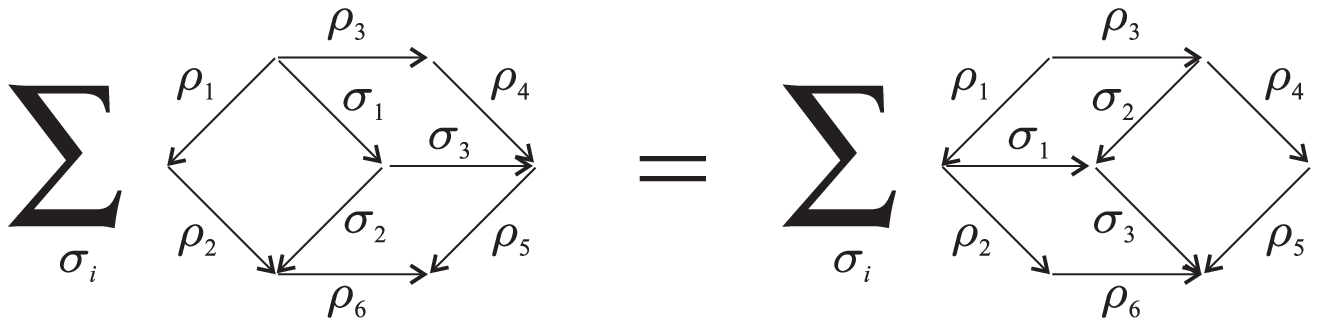}\\
 \caption{The Yang-Baxter equation} \label{fig:YBE2}
\end{center}
\end{figure}

The above connection corresponds to the natural braid generator $g_i$, which is the Boltzmann weight at criticality, and which satisfy
\begin{eqnarray}
g_i g_j & = & g_j g_i, \quad \textrm{if } |j-i|>1, \\
g_i g_{i+1} g_i & = & g_{i+1} g_i g_{i+1}. \label{eqn:braiding_relation}
\end{eqnarray}
It was claimed in \cite{ocneanu:2000ii} that the connection satisfies the unitarity property of connections
\begin{equation} \label{eqn:unitarity_property_of_connections}
\sum_{\rho_3,\rho_4} X^{\rho_1,\rho_2}_{\rho_3,\rho_4} \; \overline{X^{\rho_1',\rho_2'}_{\rho_3,\rho_4}} = \delta_{\rho_1, \rho_1'} \delta_{\rho_2, \rho_2'},
\end{equation}
and the Yang-Baxter equation
\begin{equation} \label{eqn:YBE}
\sum_{\sigma_1, \sigma_2, \sigma_3} X^{\sigma_1,\sigma_2}_{\rho_1,\rho_2} \; X^{\rho_3,\rho_4}_{\sigma_1,\sigma_3} \; X^{\sigma_3,\rho_5}_{\sigma_2,\rho_6} = \sum_{\sigma_1, \sigma_2, \sigma_3} X^{\rho_3,\sigma_2}_{\rho_1,\sigma_1} \; X^{\sigma_1,\sigma_3}_{\rho_2,\rho_6} \; X^{\rho_4,\rho_5}_{\sigma_2,\sigma_3}.
\end{equation}
The Yang-Baxter equation (\ref{eqn:YBE}) is represented graphically in Figure \ref{fig:YBE2}. We give a proof that the connection (\ref{Connection_using_weights_W}) satisfies these two properties.

\begin{Lemma} \label{Lemma:connection_properties}
If the conditions in Definition \ref{cell_system} are satisfied, the connection defined in (\ref{Connection_using_weights_W}) satisfies the unitarity property (\ref{eqn:unitarity_property_of_connections}) and the Yang-Baxter equation (\ref{eqn:YBE}).
\end{Lemma}
\emph{Proof:}
We first show unitarity.
\begin{eqnarray*}
\lefteqn{ \sum_{\rho_3,\rho_4} X^{\rho_1,\rho_2}_{\rho_3,\rho_4} \overline{X^{\rho_1',\rho_2'}_{\rho_3,\rho_4}} } \\
& = & \sum_{\rho_3,\rho_4} \left( q^{\frac{2}{3}} \delta_{\rho_1, \rho_3} \delta_{\rho_2, \rho_4} - q^{-\frac{1}{3}} \sum_{\lambda} \frac{1}{\phi_{s(\rho_1)} \phi_{r(\rho_2)}} W_{\rho_3, \rho_4, \lambda} \overline{W_{\rho_1, \rho_2, \lambda}} \right) \\
& & \qquad \quad \times \left( q^{-\frac{2}{3}} \delta_{\rho_1', \rho_3} \delta_{\rho_2', \rho_4} - q^{\frac{1}{3}} \sum_{\lambda} \frac{1}{\phi_{s(\rho_1)} \phi_{r(\rho_2)}} W_{\rho_1', \rho_2', \lambda} \overline{W_{\rho_3, \rho_4, \lambda}} \right) \\
& = & \delta_{\rho_1, \rho_1'} \delta_{\rho_3, \rho_3'} + \sum_{\stackrel{\rho_3, \rho_4}{\scriptscriptstyle{\lambda, \lambda'}}} \frac{1}{\phi_{s(\rho_1)}^2 \phi_{r(\rho_2)}^2} W_{\rho_3, \rho_4, \lambda} \overline{W_{\rho_1, \rho_2, \lambda}} W_{\rho_1', \rho_2', \lambda} \overline{W_{\rho_3, \rho_4, \lambda}} \\
& & - \sum_{\rho_3, \rho_4, \lambda} \frac{1}{\phi_{s(\rho_1)} \phi_{r(\rho_2)}} \left( q \delta_{\rho_1, \rho_3} \delta_{\rho_2, \rho_4} W_{\rho_1', \rho_2', \lambda} \overline{W_{\rho_3, \rho_4, \lambda}} \right. \\
& & \qquad \qquad \qquad \qquad \qquad \left. + q^{-1} \delta_{\rho_1', \rho_3} \delta_{\rho_2', \rho_4} W_{\rho_3, \rho_4, \lambda} \overline{W_{\rho_1, \rho_2, \lambda}} \right) \\
& = & \delta_{\rho_1, \rho_1'} \delta_{\rho_3, \rho_3'} + \sum_{\lambda, \lambda'} \frac{1}{\phi_{s(\rho_1)}^2 \phi_{r(\rho_2)}^2} \overline{W_{\rho_1, \rho_2, \lambda}} W_{\rho_1', \rho_2', \lambda} [2] \phi_{s(\rho_3)} \phi_{r(\rho_4)} \delta_{\lambda, \lambda'} \\
& & - (q + q^{-1}) \sum_{\lambda} \frac{1}{\phi_{s(\rho_1)} \phi_{r(\rho_2)}} W_{\rho_1', \rho_2', \lambda} \overline{W_{\rho_1, \rho_2, \lambda}} \\
& = & \delta_{\rho_1, \rho_1'} \delta_{\rho_3, \rho_3'},
\end{eqnarray*}
since $q + q^{-1} = [2]$, where we have used Ocneanu's type I equation (\ref{eqn:typeII_frame}) in the penultimate equality.

We now show that the connection satisfies the Yang-Baxter equation. For the left hand side of (\ref{eqn:YBE}) we have
\begin{eqnarray*}
\lefteqn{ \sum_{\sigma_1, \sigma_2, \sigma_3} X^{\sigma_1,\sigma_2}_{\rho_1,\rho_2} X^{\rho_3,\rho_4}_{\sigma_1,\sigma_3} X^{\sigma_3,\rho_5}_{\sigma_2,\rho_6} } \\
& = & \sum_{\sigma_1, \sigma_2, \sigma_3} \left( q^{\frac{2}{3}} \delta_{\rho_1, \sigma_1} \delta_{\rho_2, \sigma_2} - q^{-\frac{1}{3}} \mathcal{U}^{\sigma_1,\sigma_2}_{\rho_1,\rho_1} \right) \left( q^{-\frac{2}{3}} \delta_{\sigma_1, \rho_3} \delta_{\sigma_3, \rho_4} - q^{\frac{1}{3}} \mathcal{U}^{\rho_3,\rho_4}_{\sigma_1, \sigma_3} \right) \\
& & \qquad \qquad \times \left( q^{-\frac{2}{3}} \delta_{\sigma_2, \sigma_3} \delta_{\rho_6, \rho_5} - q^{\frac{1}{3}} \mathcal{U}^{\sigma_3,\rho_5}_{\sigma_2,\rho_6} \right) \\
& = & q^2 \delta_{\rho_1, \rho_3} \delta_{\rho_2, \rho_4} \delta_{\rho_5, \rho_6} - q \delta_{\rho_1, \rho_3} \; \mathcal{U}^{\rho_4,\rho_5}_{\rho_2,\rho_6} - q \delta_{\rho_5, \rho_6} \; \mathcal{U}^{\rho_3,\rho_4}_{\rho_1,\rho_2} - q \delta_{\rho_5, \rho_6} \; \mathcal{U}^{\rho_3,\rho_4}_{\rho_1,\rho_2} \\
& & + \sum_{\sigma_3} \mathcal{U}^{\rho_3,\rho_4}_{\rho_1,\sigma_3} \; \mathcal{U}^{\sigma_3,\rho_5}_{\rho_2,\rho_6} + \sum_{\sigma_2} \mathcal{U}^{\rho_3,\sigma_2}_{\rho_1,\rho_2}  \; \mathcal{U}^{\rho_4,\rho_5}_{\sigma_2,\rho_6} + \delta_{\rho_5, \rho_6} \sum_{\sigma_1, \sigma_2} \mathcal{U}^{\sigma_1,\sigma_2}_{\rho_1,\rho_2}  \; \mathcal{U}^{\rho_3,\rho_4}_{\sigma_1,\sigma_2} \\
& & - q^{-1} \sum_{\sigma_i} \mathcal{U}^{\sigma_1,\sigma_2}_{\rho_1,\rho_2} \; \mathcal{U}^{\rho_3,\rho_4}_{\sigma_1,\sigma_3}  \; \mathcal{U}^{\sigma_3,\rho_5}_{\sigma_2,\rho_6}  \\
& = & q^2 \delta_{\rho_1, \rho_3} \delta_{\rho_2, \rho_4} \delta_{\rho_5, \rho_6} - q \delta_{\rho_1, \rho_3} \; \mathcal{U}^{\rho_4,\rho_5}_{\rho_2,\rho_6}  - 2q \delta_{\rho_5, \rho_6} \; \mathcal{U}^{\rho_3,\rho_4}_{\rho_1,\rho_2} \\
& & + \sum_{\sigma_3} \mathcal{U}^{\rho_3,\rho_4}_{\rho_1,\sigma_3}  \; \mathcal{U}^{\sigma_3,\rho_5}_{\rho_2,\rho_6}  + \sum_{\sigma_2} \mathcal{U}^{\rho_3,\sigma_2}_{\rho_1,\rho_2}  \; \mathcal{U}^{\rho_4,\rho_5}_{\sigma_2,\rho_6}  \\
& & + \delta_{\rho_5, \rho_6} \sum_{\stackrel{{\sigma_1, \sigma_2}}{\scriptscriptstyle{\lambda, \lambda'}}} \frac{1}{\phi_{s(\rho_1)} \phi_{r(\rho_2)} \phi_{s(\rho_3)} \phi_{r(\rho_4)}} W_{\rho_1, \rho_2, \lambda} \overline{W_{\sigma_1, \sigma_2, \lambda}} W_{\sigma_1, \sigma_2, \lambda'} \overline{W_{\rho_3, \rho_4, \lambda'}} \\
& & - q^{-1} \sum_{\stackrel{{\sigma_i, \lambda}}{\scriptscriptstyle{\lambda', \lambda''}}} \frac{1}{\phi_{s(\rho_1)}^2 \phi_{r(\rho_2)} \phi_{r(\rho_4)} \phi_{s(\sigma_2)} \phi_{r(\rho_6)}} W_{\rho_1, \rho_2, \lambda} \overline{W_{\rho_3, \rho_4, \lambda'}} \\
& & \qquad \qquad \qquad \qquad \qquad \qquad \times \overline{W_{\sigma_1, \sigma_2, \lambda}} W_{\sigma_1, \sigma_3, \lambda'}  W_{\sigma_2, \rho_6, \lambda''} \overline{W_{\sigma_1, \rho_5, \lambda''}} \\
& = & q^2 \delta_{\rho_1, \rho_3} \delta_{\rho_2, \rho_4} \delta_{\rho_5, \rho_6} - q \delta_{\rho_1, \rho_3} \; \mathcal{U}^{\rho_4,\rho_5}_{\rho_2,\rho_6}  - 2q \delta_{\rho_5, \rho_6} \; \mathcal{U}^{\rho_3,\rho_4}_{\rho_1,\rho_2} \\
& & + \sum_{\sigma_3} \mathcal{U}^{\rho_3,\rho_4}_{\rho_1,\sigma_3}  \; \mathcal{U}^{\sigma_3,\rho_5}_{\rho_2,\rho_6}  + \sum_{\sigma_2} \mathcal{U}^{\rho_3,\sigma_2}_{\rho_1,\rho_2}  \; \mathcal{U}^{\rho_4,\rho_5}_{\sigma_2,\rho_6}  \\
& & + \delta_{\rho_5, \rho_6} \sum_{\lambda, \lambda'} \frac{1}{\phi_{s(\rho_1)} \phi_{r(\rho_2)} \phi_{s(\rho_3)} \phi_{r(\rho_4)}} W_{\rho_1, \rho_2, \lambda} \overline{W_{\rho_3, \rho_4, \lambda'}} [2] \phi_{r(\rho_2)} \phi_{s(\rho_1)} \delta_{\lambda, \lambda'} \\
& & - q^{-1} \sum_{\lambda, \lambda'} \frac{1}{\phi_{s(\rho_1)}^2 \phi_{r(\rho_2)} \phi_{r(\rho_4)} \phi_{r(\rho_6)}} W_{\rho_1, \rho_2, \lambda} \overline{W_{\rho_3, \rho_4, \lambda'}} \\
& & \qquad \qquad \times \left( \delta_{\lambda, \rho_6} \delta_{\lambda', \rho_5} \phi_{r(\rho_2)} \phi_{r(\rho_6)} \phi_{r(\rho_4)} + \delta_{\lambda, \lambda'} \delta_{\rho_5, \rho_6} \phi_{s(\rho_1)} \phi_{r(\rho_2)} \phi_{r(\rho_6)} \right) \\
& = & q^2 \delta_{\rho_1, \rho_3} \delta_{\rho_2, \rho_4} \delta_{\rho_5, \rho_6} - q \delta_{\rho_1, \rho_3} \; \mathcal{U}^{\rho_4,\rho_5}_{\rho_2,\rho_6}  - q \delta_{\rho_5, \rho_6} \; \mathcal{U}^{\rho_3,\rho_4}_{\rho_1,\rho_2} + \sum_{\sigma_3} \mathcal{U}^{\rho_3,\rho_4}_{\rho_1,\sigma_3} \; \mathcal{U}^{\sigma_3,\rho_5}_{\rho_2,\rho_6} \\
& & + \sum_{\sigma_2} \mathcal{U}^{\rho_3,\sigma_2}_{\rho_1,\rho_2} \; \mathcal{U}^{\rho_4,\rho_5}_{\sigma_2,\rho_6} - q^{-1} \frac{1}{\phi_{s(\rho_1)}} W_{\rho_1, \rho_2, \rho_6} \overline{W_{\rho_3, \rho_4, \rho_5}}.
\end{eqnarray*}
Computing the right hand side of (\ref{eqn:YBE}) in the same way, we arrive at the same expression.
\hfill
$\Box$

\section{Computation of the cells $W$ for $\mathcal{ADE}$ graphs} \label{Sect:computation_of_weights}

In the remaining sections we will compute cells systems $W$ for each $\mathcal{ADE}$ graph $\mathcal{G}$, with the exception of the graph $\mathcal{E}_4^{(12)}$.

Let $\triangle_{i,j,k}^{(\alpha, \beta, \gamma)}$ be the triangle $i \stackrel{\alpha}{\longrightarrow} j \stackrel{\beta}{\longrightarrow} k \stackrel{\gamma}{\longrightarrow} i$ in $\mathcal{G}$. For most of the $\mathcal{ADE}$ graph, using the equations (\ref{eqn:typeI_frame}) and (\ref{eqn:typeII_frame}) only, we can compute the cells up to choice of phase $W(\triangle_{i,j,k}^{(\alpha, \beta, \gamma)}) = \lambda_{i,j,k}^{\alpha, \beta, \gamma} |W(\triangle_{i,j,k}^{(\alpha, \beta, \gamma)})|$ for some $\lambda_{i,j,k}^{\alpha, \beta, \gamma} \in \mathbb{T}$, and also obtain some restrictions on the values which the phases $\lambda_{i,j,k}^{\alpha, \beta, \gamma}$ may take. However, for the graph $\mathcal{D}^{(n)\ast}$, $n=5,6,\ldots \;$, we impose a $\mathbb{Z}_3$ symmetry on our solutions, whilst for the graphs $\mathcal{D}^{(3k)}$, $k=2,3,\ldots \;$, and $\mathcal{E}_1^{(12)}$ we seek an orbifold solution obtained using the identification of the graphs $\mathcal{D}^{(3k)}$, $\mathcal{E}_1^{(12)}$ as $\mathbb{Z}_3$ orbifolds of $\mathcal{A}^{(3k)}$, $\mathcal{E}_2^{(12)}$ respectively.
There is still much freedom in the actual choice of phases, so that the cell system is not unique. We therefore define an equivalence relation between two cell systems:

\begin{Def}
Two families of cells $W_1$, $W_2$ which give a cell system for $\mathcal{G}$ are equivalent if, for each pair of adjacent vertices $i$, $j$ of $\mathcal{G}$, we can find a family of unitary matrices $(u(\sigma_1, \sigma_2))_{\sigma_1, \sigma_2}$, where $\sigma_1$, $\sigma_2$ are any pair of edges from $i$ to $j$, such that
\begin{equation}\label{eqn:def-equvialence_of_W1,W2}
W_1(\triangle_{i,j,k}^{(\sigma, \rho, \gamma)}) = \sum_{\sigma', \rho', \gamma'} u(\sigma, \sigma') u(\rho, \rho') u(\gamma, \gamma') W_2(\triangle_{i,j,k}^{(\sigma', \rho', \gamma')}),
\end{equation}
where the sum is over all edges $\sigma'$ from $i$ to $j$, $\rho'$ from $j$ to $k$, and $\gamma'$ from $k$ to $i$.
\end{Def}

\begin{Lemma}
Let $W_1$, $W_2$ be two equivalent families of cells, and $X^{(1)}$, $X^{(2)}$ the corresponding connections defined using cells $W_1$, $W_2$ respectively. Then $X^{(1)}$ and $X^{(2)}$ are equivalent in the sense of \cite[p.542]{evans/kawahigashi:1998}, i.e. there exists a set of unitary matrices $(u(\rho,\sigma))_{\rho,\sigma}$ such that
$$X^{(1)\rho_1,\rho_2}_{\rho_3,\rho_4} = \sum_{\sigma_i} u(\rho_3,\sigma_3) u(\rho_4,\sigma_4) \overline{u(\rho_1,\sigma_1)} \overline{u(\rho_2,\sigma_2)} X^{(2)\sigma_1,\sigma_2}_{\sigma_3,\sigma_4}.$$
\end{Lemma}

Let $W_l(\triangle_{i,j,k}^{(\sigma, \rho, \gamma)}) = \lambda_{i,j,k}^{(l)\sigma, \rho, \gamma} |W_l(\triangle_{i,j,k}^{(\sigma, \rho, \gamma)})|$, for $l = 1,2$, be two families of cells which give cell systems. If $|W_1(\triangle_{i,j,k}^{(\sigma, \rho, \gamma)})| = |W_2(\triangle_{i,j,k}^{(\sigma, \rho, \gamma)})|$, so that $W_1$ and $W_2$ differ only up to phase choice, then the equation (\ref{eqn:def-equvialence_of_W1,W2}) becomes \begin{equation}\label{eqn:def-equvialence_of_W1,W2(lambdas)}
\lambda_{i,j,k}^{(1)\sigma, \rho, \gamma} = \sum_{\sigma', \rho', \gamma'} u(\sigma, \sigma') u(\rho, \rho') u(\gamma, \gamma') \lambda_{i,j,k}^{(2)\sigma, \rho, \gamma}.
\end{equation}

For graphs with no multiple edges we write $\triangle_{i,j,k}$ for the triangle $\triangle_{i,j,k}^{(\alpha, \beta, \gamma)}$. For such graphs, two solutions $W_1$ and $W_2$ differ only up to phase choice, and (\ref{eqn:def-equvialence_of_W1,W2(lambdas)}) becomes
\begin{equation}\label{eqn:equvialence_of_W1,W2(lambdas)-no_multiple_edges}
\lambda_{i,j,k}^{(1)} = u_{\sigma} u_{\rho} u_{\gamma} \lambda_{i,j,k}^{(2)},
\end{equation}
where $u_{\sigma}, u_{\rho}, u_{\gamma} \in \mathbb{T}$ and $\sigma$ is the edge from $i$ to $j$, $\rho$ the edge from $j$ to $k$ and $\gamma$ the edge from $k$ to $i$.

We will write $U^{(x,y)}$ for the matrix indexed by the vertices of $\mathcal{G}$, with entries given by  $\mathcal{U}^{\rho_1, \rho_2}_{\rho_3, \rho_4}$ for all edges $\rho_i$, $i=1,2,3,4$ on $\mathcal{G}$ such that $s(\rho_1) = s(\rho_3) = x$, $r(\rho_2) = r(\rho_4) = y$, i.e. $[U^{(s(\rho_1),r(\rho_2))}]_{r(\rho_1),r(\rho_3)} = \mathcal{U}^{\rho_1, \rho_2}_{\rho_3, \rho_4}$.

We first present some relations that the quantum numbers $[a]_q$ satisfy, which are easily checked:

\begin{Lemma}
\begin{itemize}
\item[(i)] If $q=\exp (i \pi / n)$ then $[a]_q = [n-a]_q$, for any $a = 1,2,\ldots, n-1$,
\item[(ii)] For any $q$, $[a]_q - [a-2]_q = [2a-2]_q/[a-1]_q$, for any $a \in \mathbb{N}$,
\item[(iii)] For any $q$, $[a]_q^2 - [a-1]_q [a+1]_q = 1$ and $[a]_q [a+b]_q - [a-1]_q [a+b+1]_q = [b+1]_q$, for any $a \in \mathbb{N}$.
\end{itemize}
\end{Lemma}

\section{$\mathcal{A}$ graphs}

\begin{figure}[tb]
\begin{center}
\includegraphics[width=60mm]{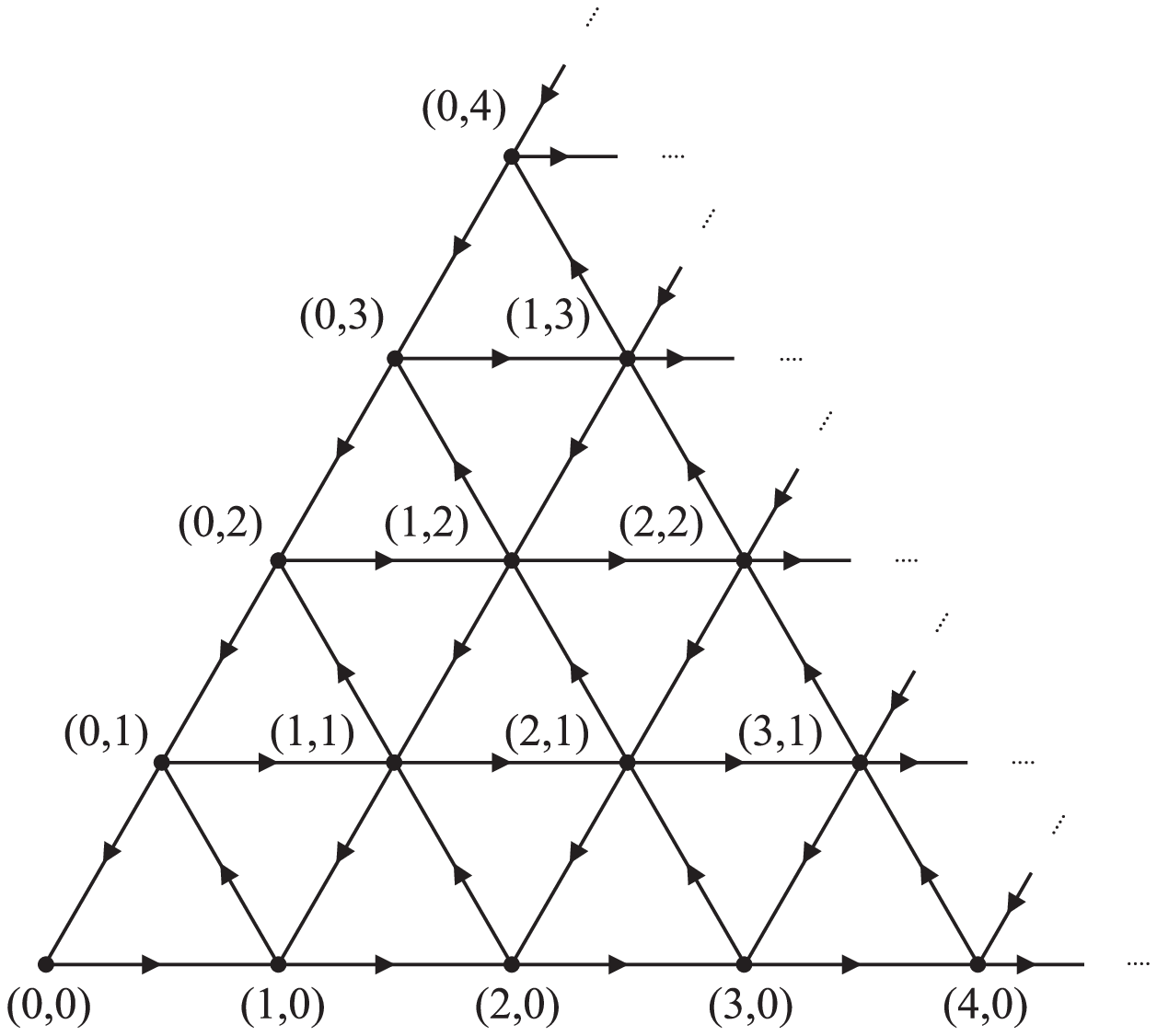}\\
 \caption{The infinite graph $\mathcal{A}^{(\infty)}$} \label{Fig:SU(3)-A(infty)}
\end{center}
\end{figure}

The infinite graph $\mathcal{A}^{(\infty)}$ is illustrated in Figure \ref{Fig:SU(3)-A(infty)}, whilst for finite $n$, the graphs $\mathcal{A}^{(n)}$ are the subgraphs of $\mathcal{A}^{(\infty)}$, given by all the vertices $(\lambda_1, \lambda_2)$ such that $\lambda_1 + \lambda_2 \leq n-3$, and all the edges in $\mathcal{A}^{(\infty)}$ which connect these vertices. The apex vertex $(0,0)$ is the distinguished vertex. For the triangle $\triangle_{(i_1,j_1)(i_2,j_2)(i_3,j_3)} = (i_1,j_1) \rightarrow (i_2,j_2) \rightarrow (i_3,j_3) \rightarrow (i_1,j_1)$ in $\mathcal{A}^{(n)}$ we will use the notation $W_{\triangle(i,j)}$ for the cell $W(\triangle_{(i,j)(i+1,j)(i,j+1)})$ and $W_{\nabla(i,j)}$ for the cell $W(\triangle_{(i+1,j)(i,j+1)(i+1,j+1)})$.

\begin{figure}[tb]
\begin{center}
\includegraphics[width=35mm]{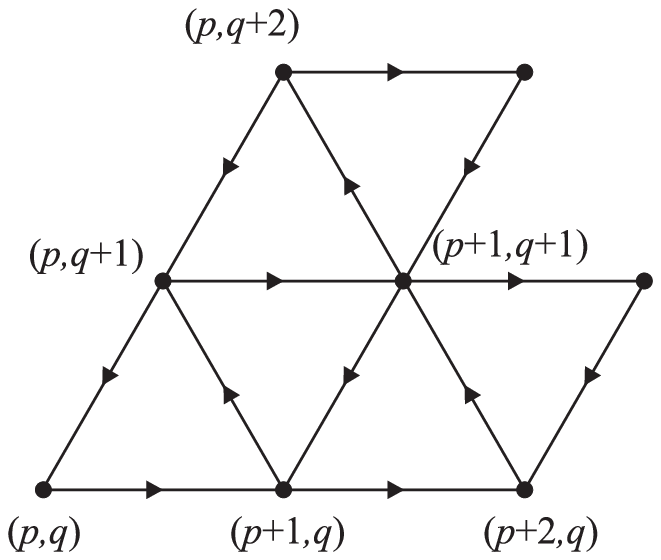}\\
 \caption{Triangles in $\mathcal{A}^{(n)}$}\label{fig:Weights1}
\end{center}
\end{figure}

\begin{Thm} \label{Thm:Aweights}
There is up to equivalence a unique set of cells for $\mathcal{A}^{(n)}$, $n < \infty$, given by
\begin{eqnarray}
\qquad \;\; W_{\triangle(k,m)} & = & \sqrt{[k+1][k+2][m+1][m+2][k+m+1][k+m+2]}/[2], \label{eqn:W_for_A-1} \\
W_{\nabla(k,m)} & = & \sqrt{[k+1][k+2][m+1][m+2][k+m+2][k+m+3]}/[2], \label{eqn:W_for_A-2}
\end{eqnarray}
for all $k,m \geq 0$. For the graph $\mathcal{A}^{(\infty)}$ with Perron-Frobenius eigenvalue $\alpha \geq 3$, there is a solution given by replacing $[j]$ by $[j]_q$ where $q = e^x$ for any $x \in \mathbb{R}$ such that $\alpha = [3]_q$.
\end{Thm}
\emph{Proof:}
Let $n < \infty$. We first prove the equalities
\begin{eqnarray}
\qquad |W_{\triangle(k,m)}| & = & \sqrt{[k+1][k+2][m+1][m+2][k+m+1][k+m+2]}/[2], \label{eqn:|W|_for_A-1}\\
|W_{\nabla(k,m)}| & = & \sqrt{[k+1][k+2][m+1][m+2][k+m+2][k+m+3]}/[2], \label{eqn:|W|_for_A-2}
\end{eqnarray}
by induction on $k,m$. The Perron-Frobenius eigenvector for $\mathcal{A}^{(n)}$ is \cite{di_francesco:1992}:
\begin{equation} \label{PF-evector_for_A}
\phi_{\lambda} = \frac{[\lambda_1+1]_q[\lambda_2+1]_q[\lambda_1+\lambda_2+2]_q}{[2]}.
\end{equation}
For the type I frame $\stackrel{(0,0)}{\bullet} \rightarrow \stackrel{(1,0)}{\bullet}$ equation (\ref{eqn:typeI_frame}) gives $|W_{\triangle(0,0)}|^2 = [2][3]$, whilst from the type I frame $\stackrel{(1,0)}{\bullet} \rightarrow \stackrel{(0,1)}{\bullet}$ we obtain $|W_{\triangle(0,0)}|^2 + |W_{\nabla(0,0)}|^2 = [2][3]^2$, giving $|W_{\nabla(0,0)}|^2 = [3][4]$. We assume (\ref{eqn:|W|_for_A-1}) and (\ref{eqn:|W|_for_A-2}) are true for $(k,m)=(p,q)$. We first show (\ref{eqn:|W|_for_A-1}) is true for $(k,m)=(p+1,q)$ and $(k,m)=(p,q+1)$ (see Figure \ref{fig:Weights1}). From the type I frame $\stackrel{(p+1,q+1)}{\bullet} \rightarrow \stackrel{(p+1,q)}{\bullet}$ we get
$$|W_{\triangle(p+1,q)}|^2 + |W_{\nabla(p,q)}|^2 = [p+2]^2[q+1][q+2][p+q+2][p+q+3]/[2],$$
and substituting in for $|W_{\nabla(p,q)}|^2$ we obtain
\begin{eqnarray*}
\lefteqn{ |W_{\triangle(p+1,q)}|^2 } \\
& = & [p+2][q+1][q+2][p+q+2][p+q+3]([2][p+2]-[p+1])/[2]^2 \\
& = & [p+2][p+2][q+1][q+2][p+q+2][p+q+3]/[2]^2.
\end{eqnarray*}
Similarly, from the type I frame $\stackrel{(p,q+1)}{\bullet} \rightarrow \stackrel{(p+1,q+1)}{\bullet}$ we get
$$|W_{\triangle(p,q+1)}|^2 = [p+1][p+2][q+2][q+3][p+q+2][p+q+3]/[2]^2,$$
as required.

For $k,m \geq 0$, the equality in (\ref{eqn:|W|_for_A-2}) follows from (\ref{eqn:|W|_for_A-1}) by considering the type I frame $\stackrel{(k+1,m)}{\bullet} \rightarrow \stackrel{(k,m+1)}{\bullet}$. We get
$$|W_{\triangle(k,m)}|^2 + |W_{\nabla(k,m)}|^2 = [k+1][k+2][m+1][m+2][k+m+2]^2/[2],$$
and substituting in for $|W_{\triangle(k,m)}|^2$ we obtain
\begin{eqnarray*}
|W_{\nabla(k,m)}|^2 & = & [k+1][k+2][m+1][m+2][k+m+2] \\
& & \quad \times ([2][k+m+2]-[k+m+1])/[2]^2 \\
& = & [k+1][k+2][m+1][m+2][k+m+2][k+m+3]/[2]^2.
\end{eqnarray*}
Hence (\ref{eqn:|W|_for_A-1}) and (\ref{eqn:|W|_for_A-2}) are true for all $k,m \geq 0$.

\begin{figure}[tb]
\begin{center}
\includegraphics[width=60mm]{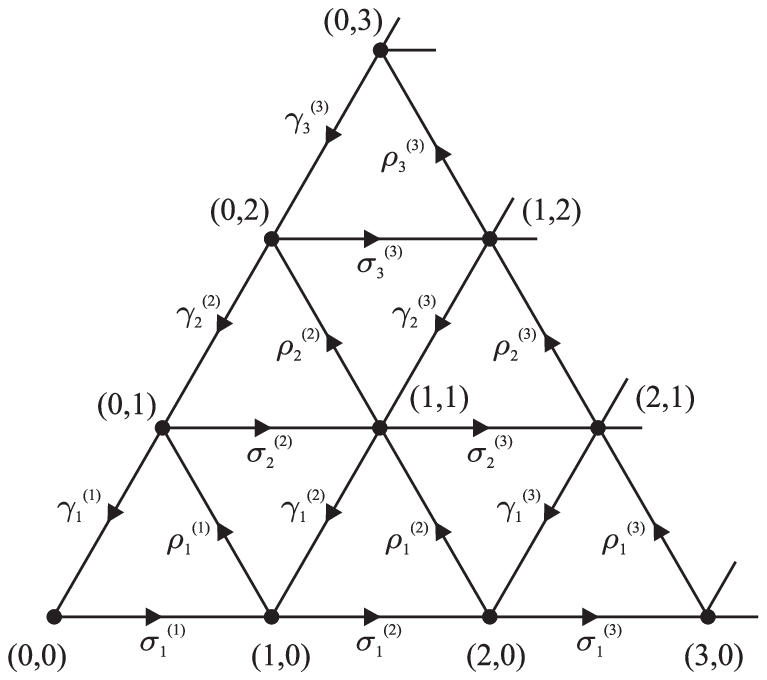}\\
 \caption{Labels for the vertices and edges of $\mathcal{A}^{(n)}$}\label{fig:Weights2}
\end{center}
\end{figure}

There is no restriction on the choice of phase for $\mathcal{A}^{(n)}$, so any choice is a solution. We now turn to the uniqueness of these cells. Let $W^{\sharp}$ be another family of cells, with $W^{\sharp}_{\triangle(k,m)} = \lambda_{(k,m)} |W_{\triangle(k,m)}|$ and $W^{\sharp}_{\nabla(k,m)} = \lambda_{(k,m)}' |W_{\nabla(k,m)}|$ (any other solution must be of this form since there are no double edges on $\mathcal{A}^{(n)}$). We label the edges of $\mathcal{A}^{(n)}$ by $\sigma_i^{(j)}$, $\rho_i^{(j)}$, $\gamma_i^{(j)}$, $j=1,\ldots, n-3$, $i=1,\ldots, j$, as shown in Figure \ref{fig:Weights2}.

Let us start with the triangle $\triangle_{(0,0)(1,0)(0,1)}$. By (\ref{eqn:equvialence_of_W1,W2(lambdas)-no_multiple_edges}) we require $1 = u_{\sigma_1^{(1)}} u_{\rho_1^{(1)}} u_{\gamma_1^{(1)}} \lambda_{(0,0)}$. Choose $u_{\sigma_1^{(1)}} = u_{\gamma_1^{(1)}} = 1$ and $u_{\rho_1^{(1)}} = \overline{\lambda_{(0,0)}}$.

Next consider the triangle $\triangle_{(1,0)(0,1)(1,1)}$. We have $1 = u_{\sigma_2^{(2)}} u_{\gamma_1^{(2)}} \overline{\lambda_{(0,0)}} \lambda_{(0,0)}'$, so choose $u_{\sigma_2^{(2)}} = 1$ and $u_{\gamma_1^{(2)}} = \lambda_{(0,0)} \overline{\lambda_{(0,0)}'}$. Similarly, setting $u_{\sigma_1^{(2)}} = u_{\gamma_2^{(2)}} = 1$, $u_{\rho_1^{(2)}} = \lambda_{(0,0)}' \overline{\lambda_{(0,0)} \lambda_{(1,0)}}$ and $u_{\rho_2^{(2)}} = \overline{\lambda_{(0,1)}}$ then (\ref{eqn:equvialence_of_W1,W2(lambdas)-no_multiple_edges}) is satisfied for the triangles $\triangle_{(1,0)(2,0)(1,1)}$ and $\triangle_{(0,1)(1,1)(0,2)}$.

Continuing in this way we set $u_{\gamma_k^{(k)}} = 1$, $u_{\gamma_i^{(k)}} = \overline{u_{\rho_i^{(k-1)}} \lambda_{(k-i-1,i-1)}'}$, for $i=1, \ldots, k-1$, and $u_{\sigma_i^{(k)}} = 1$, $u_{\rho_i^{(k)}} = u_{\rho_i^{(k-1)}} \lambda_{(k-i-1,i-1)}' \overline{\lambda_{(k-i,i-1)}}$, for $i=1, \ldots, k$, for each $k \leq n-3$. Hence, any choice of $\lambda$ and $\lambda'$ will give an equivalent solution to (\ref{eqn:W_for_A-1}), (\ref{eqn:W_for_A-2}).

For $\mathcal{A}^{(\infty)}$, we have Perron-Frobenius eigenvectors $\phi = (\phi_{\lambda_1, \lambda_2})$ given by
$$\phi_{(\lambda_1, \lambda_2)} = \frac{[\lambda_1 +1]_q [\lambda_2 +1]_q [\lambda_1 + \lambda_2 +2]_q}{[2]_q}.$$
Then the rest of the proof follows as for finite $n$.
\hfill
$\Box$

Using these cells $W$ we obtain the following representation of the Hecke algebra for $\mathcal{A}^{(n)}$. We have written the label for the rows (and columns) in front of each matrix.
\begin{eqnarray}
\quad U^{((\lambda_1, \lambda_2),(\lambda_1, \lambda_2 + 1))} & = & \begin{array}{c} \scriptstyle(\lambda_1 + 1, \lambda_2) \\ \scriptstyle(\lambda_1 - 1, \lambda_2 + 1) \end{array} \left( {\begin{array}{cc}
                 \frac{[\lambda_1 + 2]}{[\lambda_1 + 1]} & \frac{\sqrt{[\lambda_1][\lambda_1 + 2]}}{[\lambda_1 + 1]} \\
                 \frac{\sqrt{[\lambda_1][\lambda_1 + 2]}}{[\lambda_1 + 1]} & \frac{[\lambda_1]}{[\lambda_1 + 1]}
               \end{array} } \right), \\
\quad U^{((\lambda_1, \lambda_2),(\lambda_1 - 1, \lambda_2))} & = & \begin{array}{c} \scriptstyle(\lambda_1 - 1, \lambda_2 + 1) \\ \scriptstyle(\lambda_1, \lambda_2 - 1) \end{array} \left( {\begin{array}{cc}
                 \frac{[\lambda_2 + 2]}{[\lambda_2 + 1]} & \frac{\sqrt{[\lambda_2][\lambda_2 + 2]}}{[\lambda_2 + 1]} \\
                 \frac{\sqrt{[\lambda_2][\lambda_2 + 2]}}{[\lambda_2 + 1]} & \frac{[\lambda_2]}{[\lambda_2 + 1]}
               \end{array} } \right),
\end{eqnarray}
\begin{eqnarray}
\lefteqn{ U^{((\lambda_1, \lambda_2),(\lambda_1 + 1, \lambda_2 - 1))} } \\
& = & \begin{array}{c} \scriptstyle(\lambda_1 + 1, \lambda_2) \\ \scriptstyle(\lambda_1, \lambda_2 - 1) \end{array} \left( {\begin{array}{cc}
                 \frac{[\lambda_1 + \lambda_2 + 3]}{[\lambda_1 + \lambda_2 + 2]} & \frac{\sqrt{[\lambda_1 + \lambda_2 + 1][\lambda_1 + \lambda_2 + 3]}}{[\lambda_1 + \lambda_2 + 2]} \\
                 \frac{\sqrt{[\lambda_1 + \lambda_2 + 1][\lambda_1 + \lambda_2 + 3]}}{[\lambda_1 + \lambda_2 + 2]} & \frac{[\lambda_1 + \lambda_2 + 1]}{[\lambda_1 + \lambda_2 + 2]}
               \end{array} } \right). \nonumber
\end{eqnarray}

Let $e_1, e_2, e_3$ be vectors in the direction of the edges from vertex $(\lambda_1,\lambda_2)$ to the vertices $(\lambda_1+1,\lambda_2)$, $(\lambda_1-1,\lambda_2+1)$, $(\lambda_1,\lambda_2-1)$ respectively, and define an inner-product by $e_j \cdot e_k = \delta_{j,k} - 1/N$.
Wenzl \cite{wenzl:1988} constructed representations of the Hecke algebra, which are given in \cite{di_francesco/zuber:1990} as:
\begin{equation} \label{eqn:Boltzmann-Weights_for_A}
\begin{array}{ccc}
\lambda & \longrightarrow & \lambda + e_k \\
\big\downarrow & & \big\downarrow \\
\lambda + e_j & \longrightarrow & \lambda + e_j + e_k
\end{array} = (1- \delta_{jl}) \frac{\sqrt{s_{jl} (\lambda' + e_j) s_{jl} (\lambda' + e_k) }}{s_{jl} (\lambda')},
\end{equation}
where $\lambda = (\lambda_1,\lambda_2)$ is a vertex on $\mathcal{A}^{(n)}$, $\lambda' = (\lambda_1+1,\lambda_2+1)$, and $s_{jl} (\lambda) = \sin ((\pi / n) (e_j - e_l) \cdot \lambda)$. Note that this weight is 0 when $j=l$.

\begin{Lemma}
The weights in the representation of the Hecke algebra given above for $\mathcal{A}^{(n)}$ are identical to those in (\ref{eqn:Boltzmann-Weights_for_A}).
\end{Lemma}
\emph{Proof:}
For $j=l$ the result is immediate since there is no triangle $\lambda \rightarrow \lambda + e_j \rightarrow \lambda + 2e_j \rightarrow \lambda$ on $\mathcal{A}^{(n)}$, and hence the weight in our representation of the Hecke algebra will be zero also.
For an arbitrary vertex $\lambda = (\lambda_1,\lambda_2)$ of $\mathcal{A}^{(n)}$, $s_{jl}(\lambda') = \sin((\pi / n) (e_j - e_l) \cdot ((\lambda_1+1) e_1 - (\lambda_2+1) e_3))$. We will show the result for $j=1$, $l=2$ (the other cases follow similarly).
We have $s_{12}(\lambda') = \sin((\lambda_1+1) \pi/n)$ and $s_{12}(\lambda' + e_j) = s_{12}(\lambda' + e_1) = \sin((\lambda_1+2)\pi/n)$. We also have $s_{12}(\lambda' + e_2) = \sin(\lambda_1\pi/n)$.
Then for $k=1$, (\ref{eqn:Boltzmann-Weights_for_A}) becomes
$$\frac{\sqrt{\sin^2((\lambda_1+2)\pi/n)}}{\sin((\lambda_1+1) \pi/n)} = \frac{[\lambda_1+2]}{[\lambda_1+1]},$$
whilst for $k=2$, (\ref{eqn:Boltzmann-Weights_for_A}) becomes
$$\frac{\sqrt{\sin((\lambda_1+2)\pi/n)\sin(\lambda_1\pi/n)}}{\sin((\lambda_1+1) \pi/n)} = \frac{\sqrt{[\lambda_1][\lambda_1+2]}}{[\lambda_1+1]},$$
as required.
\hfill
$\Box$

\section{$\mathcal{D}$ graphs}

\begin{figure}[tb]
\begin{center}
\includegraphics[width=115mm]{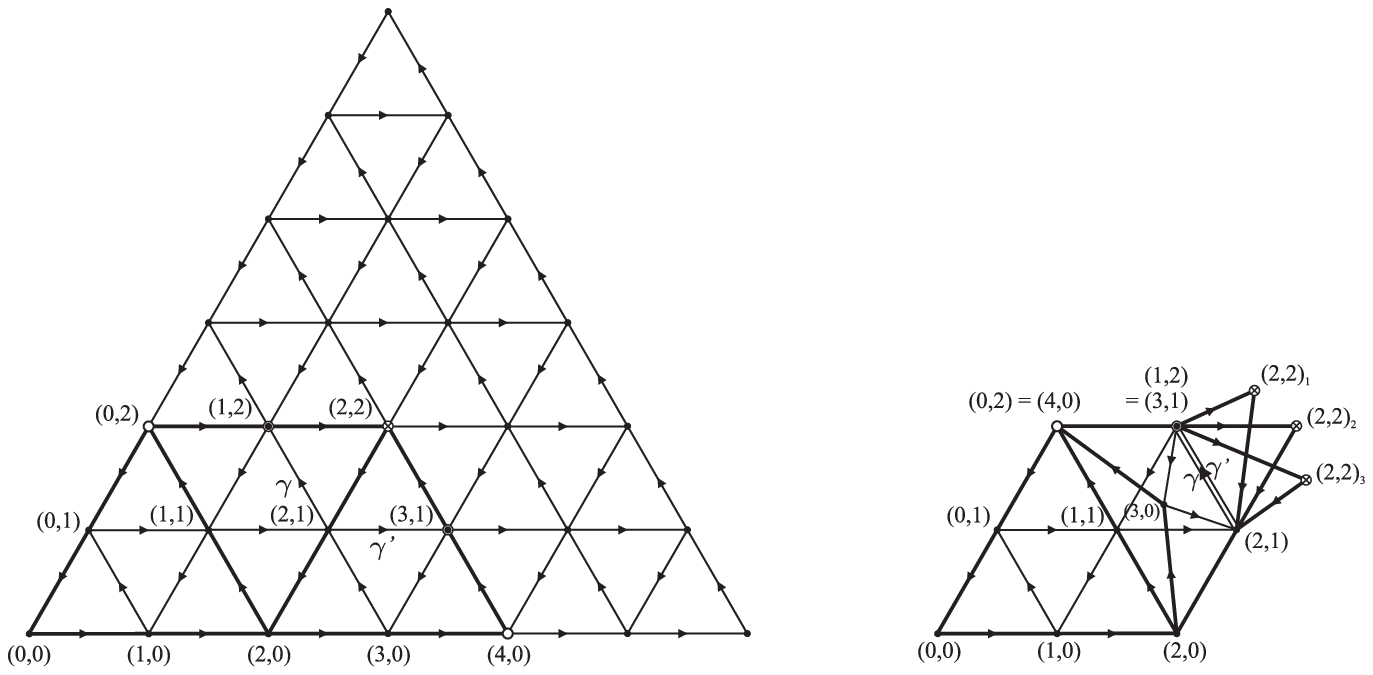}\\
 \caption{$\mathcal{A}^{(9)}$ and its $\mathbb{Z}_3$ orbifold $\mathcal{D}^{(9)}$}\label{fig:Weights3}
\end{center}
\end{figure}

\begin{figure}[tb]
\begin{center}
\includegraphics[width=40mm]{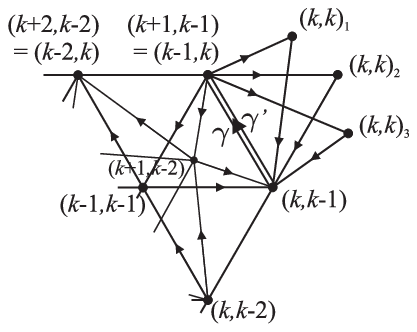}\\
 \caption{Labels for the graph $\mathcal{D}^{(3k+3)}$}\label{fig:Weights4}
\end{center}
\end{figure}

The Perron-Frobenius weights for the vertices of $\mathcal{A}^{(n)}$ are invariant under the $\mathbb{Z}_3$ symmetry given by rotation by $2 \pi / 3$. The graph $\mathcal{D}^{(n)}$ is obtained from the graph $\mathcal{A}^{(n)}$ by taking its $\mathbb{Z}_3$ orbifold, as illustrated in Figure \ref{fig:Weights3} for $n=9$ \cite{evans/kawahigashi:1994}. The Perron-Frobenius weights for the vertices of $\mathcal{D}^{(n)}$ are equal to the corresponding weights in $\mathcal{A}^{(n)}$, except that for $n=3k+3$, for integer $k \geq 1$, the vertices $(k,k)_1$, $(k,k)_2$ and $(k,k)_3$ (see Figure \ref{fig:Weights4}) which come from the fixed point $(k,k)$ of $\mathcal{A}^{(3k+3)}$ under the rotation whose Perron-Frobenius weights are a third of the weight for the vertex $(k,k)$ of $\mathcal{A}^{(3k+3)}$. The absolute values $|W^{\mathcal{A}}|$ of the cells for $\mathcal{A}^{(n)}$ are also invariant under the rotation.

Let $n \geq 5$, $n \not \equiv 0 \textrm{ mod } 3$. We will find one solution (up to a choice of phase) for the cells of $\mathcal{D}^{(n)}$ by identifying the absolute values $|W^{(\mathcal{A})}|$ for the cells in $\mathcal{A}^{(n)}$ with the absolute values $|W^{(\mathcal{D})}|$ for the corresponding cells in $\mathcal{D}^{(n)}$ when taking the orbifold. Each type I frame in $\mathcal{D}^{(n)}$ has a corresponding type I frame in $\mathcal{A}^{(n)}$, and similarly for the type II frames. Since the Perron-Frobenius weights are the same for $\mathcal{A}^{(n)}$ and $\mathcal{D}^{(n)}$, these $|W^{\mathcal{D}}|$ will certainly satisfy (\ref{eqn:typeI_frame}) and (\ref{eqn:typeII_frame}) since the $|W^{\mathcal{A}}|$ do. As in the case of $\mathcal{A}^{(n)}$, there are no restrictions on the choice of phase. Then we have the following theorem:

\begin{Thm}
Every orbifold solution for the cells of $\mathcal{D}^{(n)}$, $n \not \equiv 0 \textrm{ mod } 3$, is equivalent to the solution for which the cells in $\mathcal{D}^{(n)}$ are equal to the corresponding cells in $\mathcal{A}^{(n)}$ given in (\ref{eqn:W_for_A-1}), (\ref{eqn:W_for_A-2}).
\end{Thm}
\emph{Proof:}
The unitaries $u_{i,j} \in \mathbb{T}$, for $i,j$ vertices on $\mathcal{D}^{(n)}$, may be chosen systematically as in the proof of Theorem \ref{Thm:Aweights}, beginning with $u_{(k,k),(k,k)} = \overline{\lambda_{(k,k),(k,k),(k,k)}}^{1/3}$ if $n = 3k+4$ or $u_{(k+1,k),(k+1,k)} = \overline{\lambda_{(k+1,k),(k+1,k),(k+1,k)}}^{1/3}$ if $n = 3k+5$, and proceeding triangle by triangle.
\hfill
$\Box$

Now let $n = 3k+3$ for some integer $k \geq 1$. For $q=e^{i \pi/(3k+3)}$, we have $[(3k+3)/2 + i]_q = [(3k+3)/2 - i]_q$ where $i \in \mathbb{Z}$ for $k$ even and $i \in \mathbb{Z}+\frac{1}{2}$ for $k$ odd. In particular we will use $[2k+2 +j] = [k+1-j]$ for $j \in \mathbb{Z}$. The Perron-Frobenius weights $\phi_{(k,k)_i} = \phi_{(k,k)} /3 = [k+1]^2[2k+2]/(3[2]) = [k+1]^3/(3[2])$, $i=1,2,3$. We again find an orbifold solution for the cells for $\mathcal{D}^{(3k+3)}$, except for those which involve the vertices $(k,k)_i$, $i=1,2,3$, which correspond to the fixed point $(k,k)$ on the graph $\mathcal{A}^{(3k+3)}$. Let $\gamma$, $\gamma'$ be the two edges in the double edge of $\mathcal{D}^{(3k+3)}$, where $\gamma$ is the edge from $(k,k-1)$ to $(k-1,k)$ and $\gamma'$ is the edge from $(k,k-1)$ to $(k+1,k-1)$ in $\mathcal{A}^{(3k+3)}$ (see Figure \ref{fig:Weights3}). We will use the notation $W^{(\xi)}_{v,(k,k-1),(k-1,k)}$ to denote the cell for the triangle $\triangle_{v,(k,k-1),(k-1,k)}$ where the edge $\xi \in \{ \gamma, \gamma' \}$ is used, for $v = (k-1,k-1)$, $(k+1,k-2)$ or $(k,k)_i$, $i=1,2,3$. Then in particular we have the following:
\begin{eqnarray*}
|W^{(\gamma)}_{(k-1,k-1),(k,k-1),(k-1,k)}|^2 & = & \frac{[k]^2[k+1]^2[2k][2k+1]}{[2]^2} \\
& = & \frac{[k]^2[k+1]^2[k+2][k+3]}{[2]^2}, \\
|W^{(\gamma')}_{(k+1,k-2),(k,k-1),(k-1,k)}|^2 & = & \frac{[k-1][k][k+1][k+2][2k+1][2k+2]}{[2]^2} \\
& = & \frac{[k-1][k][k+1]^2[k+2]^2}{[2]^2}.
\end{eqnarray*}
Since $\gamma'$ is not an edge used to form the triangle $\triangle_{(k-1,k-1),(k,k-1),(k-1,k)}$ in $\mathcal{A}^{(3k+3)}$, we obtain $W^{(\gamma')}_{(k-1,k-1),(k,k-1),(k-1,k)} = 0$. Similarly we obtain $W^{(\gamma)}_{(k+1,k-2),(k,k-1),(k-1,k)} = 0$. The cells involving the vertices $(k,k)_i$ coming from the triplicated vertex $(k,k)$ in $\mathcal{A}^{(3k+3)}$ will then be a third of the corresponding cells for $\mathcal{A}^{(3k+3)}$, since the type I frames $\stackrel{(k-1,k)}{\bullet} \rightarrow \stackrel{(k,k)_i}{\bullet}$ give $|W^{(\gamma)}_{(k-1,k),(k,k)_i,(k,k-1)}|^2 + |W^{(\gamma')}_{(k-1,k),(k,k)_i,(k,k-1)}|^2 = [k][k+1]^4[k+2]/(3[2])$ for $i=1,2,3$. So
$$|W^{(\gamma)}_{(k-1,k),(k,k)_i,(k,k-1)}|^2 = \frac{1}{3} |W_{(k-1,k),(k,k),(k,k-1)}|^2 = \frac{1}{3} \frac{[k]^2[k+1]^3[k+2]}{[2]^2},$$
$$|W^{(\gamma')}_{(k-1,k),(k,k)_i,(k,k-1)}|^2 = \frac{1}{3} |W_{(k+1,k-1),(k,k),(k,k-1)}|^2 = \frac{1}{3} \frac{[k][k+1]^3[k+2]^2}{[2]^2}.$$

The phase $\lambda$ of the cell $W$ is the number $\lambda \in \mathbb{T}$ such that $W = \lambda |W|$. Let $\lambda_i, \lambda_i' \in \mathbb{T}$, be the choice of phase for the cells $W^{(\gamma)}_{(k-1,k),(k,k)_i,(k,k-1)}$, $W^{(\gamma')}_{(k-1,k),(k,k)_i,(k,k-1)}$ respectively.
Similarly, let $\lambda^{(\xi)}_{(k-1,k-1),(k,k-1),(k-1,k)}$ be the phase for $W^{(\xi)}_{(k-1,k-1),(k,k-1),(k-1,k)}$,
where $\xi \in \{ \gamma, \gamma' \}$, and $W_{v_1,v_2,v_3} = \lambda_{v_1,v_2,v_3} |W_{v_1,v_2,v_3}|$ for all other triangles $\triangle_{v_1,v_2,v_3}$ of $\mathcal{D}^{(3k+3)}$. The type II frame $\stackrel{(k,k-1)}{\bullet} \rightrightarrows \stackrel{(k-1,k)}{\bullet}$ gives the following restriction on the phases $\lambda_i$, $\lambda_i'$:
\begin{equation} \label{eqn:restriction_on_lambdas(IX)}
\lambda_1 \overline{\lambda_1'} + \lambda_2 \overline{\lambda_2'} + \lambda_3 \overline{\lambda_3'} = 0.
\end{equation}
From the type II frame $\stackrel{(k,k)_i}{\bullet} \rightarrow \stackrel{(k,k-1)}{\bullet} \leftarrow \stackrel{(k,k)_j}{\bullet}$ we obtain $\mathrm{Re}(\lambda_i \lambda_j' \overline{\lambda_i' \lambda_j}) = - 1/2$ for $i \neq j$, giving $\lambda_i \overline{\lambda_i'} = (-1/2 + \varepsilon_{ij} \sqrt{3} i /2) \lambda_j \overline{\lambda_j'}$, $\varepsilon_{ij} \in \{ \pm 1 \}$. Note that $\varepsilon_{ji} = - \varepsilon_{ij}$, and substituting for $\lambda_i \overline{\lambda_i'}$ with $j=i+1$ into (\ref{eqn:restriction_on_lambdas(IX)}) we find $\varepsilon_{12} = \varepsilon_{23} = \varepsilon_{31}$. Then we have
\begin{equation} \label{eqn:restriction_on_lambdas(X)}
\lambda_i \overline{\lambda_i'} = (-\frac{1}{2} + \varepsilon \frac{\sqrt{3} i}{2}) \lambda_{i+1} \overline{\lambda_{i+1}'},
\end{equation}
for $\varepsilon \in \{ \pm 1 \}$, $i = 1,2,3$ (mod 3).
Then there are two solutions for the cells of $\mathcal{D}^{(3k+3)}$, $W$ and its complex conjugate $\overline{W}$. The solution $\overline{W}$ is the solution to the graph where we switch vertices $(k,k)_2 \leftrightarrow (k,k)_3$.

\begin{Thm}
Every orbifold solution for the cells of $\mathcal{D}^{(3k+3)}$ is given, up to equivalence, by the inequivalent solutions $W$ or its complex conjugate $\overline{W}$, where $W$ is given by
$$W^{(\gamma)}_{(k-1,k),(k,k)_i,(k,k-1)} = \epsilon_i \frac{[k]\sqrt{[k+1]^3[k+2]}}{\sqrt{3} \, [2]},$$
$$W^{(\gamma')}_{(k-1,k),(k,k)_i,(k,k-1)} = \overline{\epsilon_i} \frac{[k+2]\sqrt{[k][k+1]^3}}{\sqrt{3} \, [2]},$$
$$W^{(\gamma)}_{(k-1,k-1),(k,k-1),(k-1,k)} = \frac{[k][k+1]\sqrt{[k+2][k+3]]}}{[2]},$$
$$W^{(\gamma')}_{(k+1,k-2),(k,k-1),(k-1,k)} = \frac{[k+1][k+2]\sqrt{[k-1][k]}}{[2]},$$
$$W^{(\gamma')}_{(k-1,k-1),(k,k-1),(k-1,k)} = W^{(\gamma)}_{(k+1,k-2),(k,k-1),(k-1,k)} = 0,$$
where $\epsilon_1 = 1$, $\epsilon_2 = e^{2 \pi i/3} = \overline{\epsilon_3}$, and all other cells are equal to the corresponding cells in $\mathcal{A}^{(3k+3)}$ given in (\ref{eqn:W_for_A-1}), (\ref{eqn:W_for_A-2}).
\end{Thm}
\emph{Proof:}
Let $W^{\sharp}$ be any orbifold solution for the cells of $\mathcal{D}^{(3k+3)}$. Then $W^{\sharp}$ is given, for $i=1,2,3$, by
\begin{eqnarray*}
W^{\sharp(\gamma)}_{(k-1,k),(k,k)_i,(k,k-1)} & = & \lambda^{\sharp}_i |W^{(\gamma)}_{(k-1,k),(k,k)_i,(k,k-1)}|, \\ W^{\sharp(\gamma')}_{(k-1,k),(k,k)_i,(k,k-1)} & = & \lambda^{\sharp}_i{}' |W^{(\gamma')}_{(k-1,k),(k,k)_i,(k,k-1)}|, \end{eqnarray*}
\begin{eqnarray*}
\lefteqn{ W^{\sharp(\xi)}_{(k-1,k-1),(k,k-1),(k-1,k)} } \\
& = & \lambda^{\sharp(\xi)}_{(k-1,k-1),(k,k-1),(k-1,k)} |W^{\sharp(\xi)}_{(k-1,k-1),(k,k-1),(k-1,k)}|,
\end{eqnarray*}
where $\xi \in \{ \gamma, \gamma' \}$, and $W^{\sharp}_{v_1,v_2,v_3} = \lambda^{\sharp}_{v_1,v_2,v_3} |W_{v_1,v_2,v_3}|$ for all other triangles $\triangle_{v_1,v_2,v_3}$ of $\mathcal{D}^{(3k+3)}$, and where the choice of $\lambda^{\sharp}_i$, $\lambda^{\sharp}_i{}'$ satisfy condition (\ref{eqn:restriction_on_lambdas(X)}) with $\varepsilon = 1$.
We need to find a family of unitaries $\{ u_{\rho} \}$ for edges $\rho \neq \gamma'$ of $\mathcal{D}^{(3k+3)}$, where $u_{\gamma} = (u_{\gamma}(\xi, \xi'))$, $\xi, \xi' \in \{ \gamma, \gamma' \}$, is a $2 \times 2$ unitary matrix, and $u_{\rho} \in \mathbb{T}$ for all other $\rho$. These unitaries must satisfy (\ref{eqn:def-equvialence_of_W1,W2(lambdas)}) and (\ref{eqn:equvialence_of_W1,W2(lambdas)-no_multiple_edges}), i.e.
$\epsilon_l = u_{\mu_l} u_{\mu_l'} (u_{\gamma}(\gamma, \gamma) \lambda^{\sharp}_l + u_{\gamma}(\gamma, \gamma') \lambda^{\sharp}_l{}')$ and $\overline{\epsilon_l} = u_{\mu_l} u_{\mu_l'} (u_{\gamma}(\gamma', \gamma) \lambda^{\sharp}_l + u_{\gamma}(\gamma', \gamma') \lambda^{\sharp}_l{}')$, for $l=1,2,3$, and
\begin{eqnarray*}
1 & = & u_{\sigma_1} u_{\sigma_2} \sum_{\xi'} u(\xi, \xi') \lambda^{\sharp(\xi')}_{(k-1,k-1),(k,k-1),(k-1,k)}, \\
1 & = & u_{\sigma_1'} u_{\sigma_2'} \sum_{\xi'} u(\xi, \xi') \lambda^{\sharp(\xi')}_{(k+1,k-2),(k,k-1),(k-1,k)}.
\end{eqnarray*}
For all other triangles $\triangle_{p_1,p_2,p_3}^{(\rho_1, \rho_2, \rho_3)}$ of $\mathcal{D}^{(3k+3)}$ we require $1 = u_{\rho_1} u_{\rho_2} u_{\rho_3} \lambda^{\sharp}_{p_1,p_2,p_3}$.

For $u_{\gamma}$ we choose $u_{\gamma}(\gamma, \gamma) = 1$, $u_{\gamma}(\gamma, \gamma') = u_{\gamma}(\gamma', \gamma) = 0$ and $u_{\gamma}(\gamma' \gamma') = \lambda^{\sharp}_1 \overline{\lambda^{\sharp}_1}$. We set $u_{\mu_l'} = 1$ and $u_{\mu_l} = \epsilon_l \overline{\lambda^{\sharp}_l}$, for $l=1,2,3$, and $u_{\sigma_1} = u_{\sigma_1'} = 1$, $u_{\sigma_2} = \overline{\lambda^{\sharp(\gamma)}_{(k-1,k-1),(k,k-1),(k-1,k)}}$ and $u_{\sigma_2'} = \overline{\lambda^{\sharp(\gamma')}_{(k+1,k-2),(k,k-1),(k-1,k)}}$.

\begin{figure}[tb]
\begin{center}
\includegraphics[width=60mm]{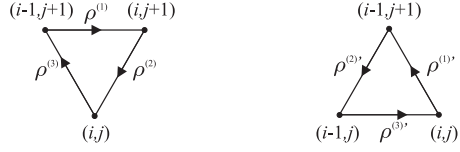}\\
 \caption{Triangles $\triangle_{(i,j),(i-1,j+1),(i,j+1)}^{(\rho^{(1)}, \rho^{(2)}, \rho^{(3)})}$ and $\triangle_{(i-1,j),(i,j),(i-1,j+1)}^{(\rho^{(1)}{}', \rho^{(2)}{}', \rho^{(3)}{}')}$}\label{fig:Weights5}
\end{center}
\end{figure}

For the remaining triangles we proceed as follows. Let $m=2k-2$. For each triangle $\triangle_{(i,j),(i-1,j+1),(i,j+1)}^{(\rho^{(1)}, \rho^{(2)}, \rho^{(3)})}$ as in Figure \ref{fig:Weights5} (and similarly for triangles $\triangle_{(i,j),(i-1,j+1),(i,j+1)}$) such that $i+j=m$, if either $u_{\rho^{(1)}}$ or $u_{\rho^{(2)}}$ hasn't yet been assigned a value we set it to be 1, and set $u_{\rho^{(3)}} = \overline{u_{\rho^{(1)}} u_{\rho^{(2)}} \lambda^{\sharp}_{(i,j),(i-1,j+1),(i,j+1)}}$. Next, for each triangle $\triangle_{(i-1,j),(i,j),(i-1,j+1)}^{(\rho^{(1)}{}', \rho^{(2)}{}', \rho^{(3)}{}')}$ as in Figure \ref{fig:Weights5} (and similarly for triangles $\triangle_{(i+1,j-1),(i,j),(i+1,j)}$) such that $i+j=m$, if either $u_{\rho^{(1)}{}'}$ or $u_{\rho^{(2)}{}'}$ hasn't yet been assigned a value we set it to be 1, and set $u_{\rho^{(3)}{}'} = \overline{u_{\rho^{(1)}{}'} u_{\rho^{(2)}{}'} \lambda^{\sharp}_{(i-1,j),(i,j),(i-1,j+1)}}$. We then set $m = 2k-3$ and repeat the above steps. Continuing in this way, for $m=2k-4, \ldots, 3$, we find the required unitaries $\{ u_{\rho} \}$. The proof for the uniqueness of the complex conjugate solution can be shown similarly.

For the solutions $W$ and $\overline{W}$ to be equivalent, we require unitaries as above such that
\begin{eqnarray*}
\epsilon_l & = & u_{\mu_l} u_{\mu_l'} (u_{\gamma}(\gamma, \gamma) \overline{\epsilon}_l + \frac{\sqrt{[k+2]}}{\sqrt{[k]}} u_{\gamma}(\gamma, \gamma') \epsilon_l), \\
\overline{\epsilon}_l & = & u_{\mu_l} u_{\mu_l'} (\frac{\sqrt{[k]}}{\sqrt{[k+2]}} u_{\gamma}(\gamma', \gamma) \overline{\epsilon}_l + u_{\gamma}(\gamma', \gamma') \epsilon_l),
\end{eqnarray*}
for $l=1,2,3$. This forces $u_{\gamma}(\gamma, \gamma) = u_{\gamma}(\gamma', \gamma') = 0$, $u_{\gamma}(\gamma, \gamma') = \sqrt{[k]}/\sqrt{[k+2]}$ and $u_{\gamma}(\gamma', \gamma) = \sqrt{[k+2]}/\sqrt{[k]}$. But then $u_{\gamma}$ is not a unitary.
\hfill
$\Box$

Using the cells $W$ we obtain the following representation of the Hecke algebra for $\mathcal{D}^{(3k+3)}$, we use the notation $v^{(\gamma)}$ if the path uses the edge $\gamma$, where $v$ is a vertex of $\mathcal{D}^{(3k+3)}$.:
\begin{eqnarray*}
\lefteqn{ U^{((k-1,k-1),(k-1,k))} \;\; = \;\; \begin{array}{c} \scriptstyle(k,k-1)^{(\gamma)} \\ \scriptstyle(k,k-1)^{(\gamma')} \\ \scriptstyle(k-2,k) \end{array} \left( {\begin{array}{ccc}
                 \frac{[k+1]}{[k]} & 0 & \frac{\sqrt{[k-1][k+1]}}{[k]} \\
                 0 & 0 & 0 \\
                 \frac{\sqrt{[k-1][k+1]}}{[k]} & 0 & \frac{[k-1]}{[k]}
               \end{array} } \right), } \\
& & = \;\; U^{((k,k-1),(k-1,k-1))} \quad \textrm{ with rows labelled by } (k-1,k)^{(\gamma)}, (k-1,k)^{(\gamma')},
\end{eqnarray*}
\begin{eqnarray*}
\lefteqn{ U^{((k+1,k-2),(k-1,k))} \;\; = \;\; \begin{array}{c} \scriptstyle(k,k-1)^{(\gamma)} \\ \scriptstyle(k,k-1)^{(\gamma')} \\ \scriptstyle(k-2,k) \end{array} \left( {\begin{array}{ccc}
                 0 & 0 & 0 \\
                 0 & \frac{[k+1]}{[k+2]} & \frac{\sqrt{[k+1][k+3]}}{[k+2]} \\
                 0 & \frac{\sqrt{[k+1][k+3]}}{[k+2]} & \frac{[k+3]}{[k+2]}
               \end{array} } \right), } \\
& & = \;\; U^{((k,k-1),(k+1,k-2))} \quad \textrm{ with rows labelled by } (k-1,k)^{(\gamma)}, (k-1,k)^{(\gamma')}, \\
& & \qquad (k,k-2), \\
\lefteqn{ U^{((k,k-1),(k,k)_i)} \;\; = \;\; \begin{array}{c} \scriptstyle(k-1,k)^{(\gamma)} \\ \scriptstyle(k-1,k)^{(\gamma')} \end{array} \left( {\begin{array}{cc}
                 \frac{[k]}{[k+1]} & \overline{\epsilon}_i \frac{\sqrt{[k][k+2]}}{[k+1]} \\
                 \epsilon_i \frac{\sqrt{[k][k+2]}}{[k+1]} & \frac{[k+2]}{[k+1]}
               \end{array} } \right), } \\
& & = \;\; U^{((k,k)_i,(k-1,k))} \quad \textrm{ with rows labelled by } (k,k-1)^{(\gamma)}, (k,k-1)^{(\gamma')},
\end{eqnarray*}
\begin{eqnarray*}
\lefteqn {U^{((k-1,k),(k,k-1))} } \\
& = & \begin{array}{c} \scriptstyle(k,k)_1 \\ \scriptstyle(k,k)_2 \\ \scriptstyle(k,k)_3 \\ \scriptstyle(k-1,k-1) \\ \scriptstyle(k+1,k-2) \end{array} \left( {\begin{array}{ccccc}
                 [2][k+1]a & \overline{\epsilon}a & \epsilon a & b & c \\
                 \epsilon a & [2][k+1]a & \overline{\epsilon}a & \epsilon_2 b & \overline{\epsilon}_2 c \\
                 \overline{\epsilon}a & \epsilon a & [2][k+1]a & \overline{\epsilon}_2 b & \epsilon_2 c \\
                 b & \overline{\epsilon}_2 b & \epsilon_2 b & \frac{[k+3]}{[k+2]} & 0 \\
                 c & \epsilon_2 c & \overline{\epsilon}_2 c & 0 & \frac{[k-1]}{[k]}
               \end{array} } \right),
\end{eqnarray*}
where $\epsilon = \epsilon_2[k] + \overline{\epsilon}_2[k+2]$ and
$$a = \frac{[k+1]}{3[k][k+2]}, \qquad b = \frac{\sqrt{[k+1][k+3]}}{\sqrt{3} \; [k+2]}, \qquad c = \frac{\sqrt{[k-1][k+1]}}{\sqrt{3} \; [k]}.$$
Another representation of the Hecke algebra is given by taking the complex conjugates of the weights in the representation above.

\begin{figure}[htb]
\begin{center}
\includegraphics[width=20mm]{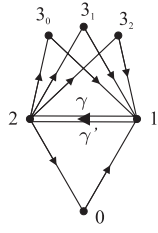}\\
 \caption{Labelling the graph $\mathcal{D}^{(6)}$}\label{fig:Weights10}
\end{center}
\end{figure}

In \cite{fendley:1989}, Fendley gives Boltzmann weights for $\mathcal{D}^{(6)}$, which at criticality and with the parameter $u=1$, give a representation of the Hecke algebra. However these Boltzmann weights are not equivalent to the representation of the Hecke algebra using the cells $W$ or $\overline{W}$. To see this, we use a similar labelling for the graph $\mathcal{D}^{(6)}$ as in \cite{fendley:1989}- see Figure \ref{fig:Weights10}.

Consider the weight $[\overline{U}^{(3_r,2)}]_{\gamma,\gamma'}$, where we label the rows and columns by $\gamma$, $\gamma'$ to denote which edge from $1$ to $2$ is used for the path of length 2 from $3_r$ to $2$, $r=0,1,2$, and the weight $\overline{U}$ is the complex conjugate of that given above, i.e. it is the weight given by the solution $\overline{W}$ for the cells of $\mathcal{D}^{(6)}$. Then for equivalence we require a unitary $u_{3_r,1} \in \mathbb{T}$ and a $2 \times 2$ unitary matrix $u_{\gamma}$ such that
\begin{eqnarray}
\epsilon_2^r \frac{\sqrt{[3]}}{[2]} & = & |u_{3_r,1}|^2 \left( u_{\gamma}(\gamma,\gamma) \overline{u_{\gamma}(\gamma',\gamma)} \frac{1}{[2]} + u_{\gamma}(\gamma,\gamma) \overline{u_{\gamma}(\gamma',\gamma')} \epsilon_2^r \frac{\sqrt{[3]}}{[2]} \right. \nonumber \\
& & \qquad \qquad \left. + u_{\gamma}(\gamma,\gamma') \overline{u_{\gamma}(\gamma',\gamma)} \overline{\epsilon}_2^r \frac{\sqrt{[3]}}{[2]} + u_{\gamma}(\gamma,\gamma') \overline{u_{\gamma}(\gamma',\gamma')} \frac{[3]}{[2]} \right).  \label{eqn:D(6)inequivalence}
\end{eqnarray}
Since $u_{\gamma}$ is independent of $r$, for (\ref{eqn:D(6)inequivalence}) to be satisfied for each $r=0,1,2$, we require $u_{\gamma}(\gamma,\gamma) \overline{u_{\gamma}(\gamma',\gamma')} = 1$ and the other terms to be zero, which gives $u_{\gamma}(\gamma,\gamma') = u_{\gamma}(\gamma',\gamma) = 0$ and $u_{\gamma}(\gamma',\gamma') = (u_{\gamma}(\gamma,\gamma))^{-1}$. But now if we consider the weight $[\overline{U}^{(1,3_r)}]_{\gamma,\gamma'}$, with $u_{2,3_r} \in \mathbb{T}$, we have
\begin{eqnarray*}
\overline{\epsilon}_2^r \frac{\sqrt{[3]}}{[2]} & = & |u_{2,3_r}|^2 \left( u_{\gamma}(\gamma,\gamma) \overline{u_{\gamma}(\gamma',\gamma)} \frac{1}{[2]} + u_{\gamma}(\gamma,\gamma) \overline{u_{\gamma}(\gamma',\gamma')} \overline{\epsilon}_2^r \frac{\sqrt{[3]}}{[2]} \right. \\
& & \qquad \qquad \left. + u_{\gamma}(\gamma,\gamma') \overline{u_{\gamma}(\gamma',\gamma)} \epsilon_2^r \frac{\sqrt{[3]}}{[2]} + u_{\gamma}(\gamma,\gamma') \overline{u_{\gamma}(\gamma',\gamma')} \frac{[3]}{[2]} \right),
\end{eqnarray*}
but $[\overline{U}^{(1,3_r)}]_{\gamma,\gamma'} = \epsilon_2^r \frac{\sqrt{[3]}}{[2]}$, for $r=0,1,2$. We obtain a similar contradiction when considering the weights $U$ defined using the solution $W$ for the cells.

Suppose however that the Boltzmann weight denoted by $\widetilde{W}^{(\widetilde{1},\widetilde{3_r})}_{\widetilde{2},\widetilde{2}}$ in \cite{fendley:1989} is the complex conjugate of that given. Then the Boltzmann weights at criticality of Fendley \cite{fendley:1989} are equivalent to the representation of the Hecke algebra given by the solution $\overline{W}$ for the cells of $\mathcal{D}^{(6)}$. We choose a family of unitaries $u_{0,1} = u_{2,0} = u_{2,3_r} = 1$, $u_{3_r,1} = \overline{\epsilon}_2^r$, $r=0,1,2$, and choose $u_{\gamma}$ to be the $2 \times 2$ identity matrix.

\section{$\mathcal{A}^{\ast}$ graphs}

\begin{figure}[tb]
\begin{center}
\includegraphics[width=100mm]{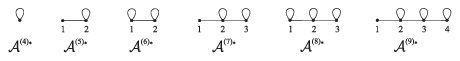}\\
 \caption{$\mathcal{A}^{(n)\ast}$ for $n = 4,5,6,7,8,9$}\label{fig:A(star)-graphs}
\end{center}
\end{figure}

The infinite series of graphs $\mathcal{A}^{(n)\ast}$ are illustrated in Figure \ref{fig:A(star)-graphs}. The graphs $\mathcal{A}^{(2n+1)\ast}$ and $\mathcal{A}^{(2n)\ast}$ are slightly different.

First we consider the graphs $\mathcal{A}^{(2n+1)\ast}$. The Perron-Frobenius weights on the vertices are given by $\phi_{i} = [2i-1]$, $i=1,\ldots,n$.

\begin{Thm}
There is up to equivalence a unique set of cells for $\mathcal{A}^{(2n+1)\ast}$, $n < \infty$, given by
\begin{eqnarray*}
W_{i-1,i,i} & = & \frac{\sqrt{[i][2i-3][2i-1]}}{\sqrt{[i-1]}}, \qquad \qquad i=2,\ldots,n, \\
W_{i,i,i+1} & = & \frac{\sqrt{[i-1][2i-1][2i+1]}}{\sqrt{[i]}}, \qquad \, i=2,\ldots,n-1, \\
W_{i,i,i} & = & (-1)^{i+1} \frac{[2i-1]}{\sqrt{[i-1][i]}}, \qquad \qquad i=2,\ldots,n.
\end{eqnarray*}
\end{Thm}
\emph{Proof:}
Using (\ref{eqn:typeI_frame}), (\ref{eqn:typeII_frame}) we obtain
\begin{eqnarray}
|W_{i-1,i,i}|^2 & = & \frac{[i][2i-3][2i-1]}{[i-1]}, \qquad \qquad i=2,\ldots,n, \label{eqn:|W|_for_A(2n+1star)-1} \\
|W_{i,i,i+1}|^2 & = & \frac{[i-1][2i-1][2i+1]}{[i]}, \qquad \, i=2,\ldots,n-1, \label{eqn:|W|_for_A(2n+1star)-2} \\
|W_{i,i,i}|^2 & = & \frac{[2i-1]^2}{[i-1][i]}, \qquad \qquad \qquad i=2,\ldots,n. \label{eqn:|W|_for_A(2n+1star)-3}
\end{eqnarray}

Let $W_{i,j,k} = \lambda_{i,j,k} |W_{i,j,k}|$ for $\lambda_{i,j,k} \in \mathbb{T}$. From type II frames we have the restriction
\begin{equation} \label{eqn:restriction_on_lambdas-A(2n+1star)}
\lambda_{i,i,i+1}^3 \lambda_{i+1,i+1,i+1} = - \lambda_{i,i+1,i+1}^3 \lambda_{i,i,i},
\end{equation}
for $i = 2, \ldots, n-1$.
Let $W^{\sharp}_{i,j,k} = \lambda^{\sharp}_{i,j,k} |W_{i,j,k}|$ be any other solution to the cells, where the $\lambda^{\sharp}$ satisfy (\ref{eqn:restriction_on_lambdas-A(2n+1star)}). We need to find a family of unitaries $\{ u_{i,j} \}$, where $u_{i,j}$ is the unitary for the edge from vertex $i$ to vertex $j$ on $\mathcal{A}^{(2n+1)\ast}$, which satisfy (\ref{eqn:equvialence_of_W1,W2(lambdas)-no_multiple_edges}), i.e. $-1 = u_{2l,2l}^3 \lambda^{\sharp}_{2l,2l,2l}$ for $l=1,\ldots,\lfloor n/2 \rfloor$, and $1 = u_{i} u_{j} u_{k} \lambda^{\sharp}_{i,j,k}$ for all other triangles $\triangle_{i,j,k}$.
We choose $u_{1,2} = 1$, $u_{2,1} = - (\lambda^{\sharp}_{2,2,2})^{1/3} \overline{\lambda^{\sharp}_{1,2,2}}$, $u_{2,2} = - (\overline{\lambda^{\sharp}_{2,2,2}})^{1/3}$,
and for $i=2,\ldots,n-1$, $u_{i,i+1} = 1$
$u_{i+1,i} = - (\overline{\lambda^{\sharp}_{2,2,2}})^{1/3} \lambda^{\sharp}_{2,3,3} \lambda^{\sharp}_{3,4,4} \cdots \lambda^{\sharp}_{i-1,i,i} \overline{\lambda^{\sharp}_{2,2,3}} \overline{\lambda^{\sharp}_{3,3,4}} \cdots \overline{\lambda^{\sharp}_{i,i,i+1}}$, and
$u_{i+1,i+1} = - (\overline{\lambda^{\sharp}_{2,2,2}})^{1/3} \lambda^{\sharp}_{2,2,3} \lambda^{\sharp}_{3,3,4} \cdots \lambda^{\sharp}_{i,i,i+1} \overline{\lambda^{\sharp}_{2,3,3}} \overline{\lambda^{\sharp}_{3,4,4}} \cdots \overline{\lambda^{\sharp}_{i,i+1,i+1}}$.
\hfill
$\Box$

For $\mathcal{A}^{(2n+1)\ast}$, the above cells $W$ give the following representation of the Hecke algebra:
\begin{eqnarray*}
U^{(i,i+1)} & = & \begin{array}{c} \scriptstyle i \\ \scriptstyle i+1 \end{array} \left( {\begin{array}{cc}
                 \frac{[i-1]}{[i]} & \frac{\sqrt{[i-1][i+1]}}{[i]} \\
                 \frac{\sqrt{[i-1][i+1]}}{[i]} & \frac{[i+1]}{[i]}
               \end{array} } \right), \\
U^{(i,i-1)} & = & \begin{array}{c} \scriptstyle i-1 \\ \scriptstyle i \end{array} \left( {\begin{array}{cc}
                 \frac{[i-2]}{[i-1]} & \frac{\sqrt{[i-2][i]}}{[i-1]} \\
                 \frac{\sqrt{[i-2][i]}}{[i-1]} & \frac{[i]}{[i-1]}
               \end{array} } \right), \\
U^{(i,i)} & = & \begin{array}{c} \scriptstyle i-1 \\ \scriptstyle i \\ \scriptstyle i+1 \end{array} \left( {\begin{array}{ccc}
                 \frac{[i][2i-3]}{[i-1][2i-1]} & \frac{(-1)^{i+1}\sqrt{[2i-3]}}{[i-1]\sqrt{[2i-1]}} & \frac{\sqrt{[2i-3][2i+1]}}{[2i-1]} \\
                 \frac{(-1)^{i+1}\sqrt{[2i-3]}}{[i-1]\sqrt{[2i-1]}} & \frac{1}{[i-1][i]} & \frac{(-1)^{i+1}\sqrt{[2i+1]}}{[i]\sqrt{[2i-1]}} \\
                 \frac{\sqrt{[2i-3][2i+1]}}{[2i-1]} & \frac{(-1)^{i+1}\sqrt{[2i+1]}}{[i]\sqrt{[2i-1]}} & \frac{[i-1][2i+1]}{[i][2i-1]}
               \end{array} } \right).
\end{eqnarray*}

In \cite{behrend/evans:2004}, Behrend and Evans give Boltzmann weights
$$W \left( \left. \begin{array}{cc} a & d \\ b & c \end{array} \right| u \right),$$
which at criticality, with $u=1$, give a representation of the Hecke algebra. (Note, these Boltzmann weights are not to be confused with the Ocneanu cells $W$.)

\begin{Lemma}
The weights in the representation of the Hecke algebra given above for $\mathcal{A}^{(2n+1)\ast}$ are equivalent to the Boltzmann weights at criticality given by Behrend-Evans in \cite{behrend/evans:2004}.
\end{Lemma}
\emph{Proof:}
To make our notation the same as that of \cite{behrend/evans:2004} one replaces $i$ with $(a+1)/2$. Then it is easily checked that the absolute values of our weights given above are equal to those for the Boltzmann weights in \cite{behrend/evans:2004}, setting $q=0$, in all but a few cases. We will show that the absolute values in these other cases are also equal. For $[U^{(i,i)}]_{i+1,i+1}$, the Boltzmann weight in \cite{behrend/evans:2004} is
$$\frac{[a+2] - [a+2]/[a]}{[a+1]} = \frac{[a+2]}{[a][a+1]}([a]-[1]) = \frac{[a+2]}{[a][a+1]}\frac{[\frac{1}{2}(a-1)][a+1]}{[\frac{1}{2}(a+1)]},$$
which is equal to our weight, and similarly for $[U^{(i,i)}]_{i-1,i-1}$. For $[U^{(i,i)}]_{i,i}$ we have to do the most work. From \cite{behrend/evans:2004} its value is
\begin{equation} \label{eqn:Boltzmann_weight_for_A2n+1(star)}
\frac{1}{[3]} \left( [2] - \frac{[a+2][\frac{1}{2}(a-5)]}{[a][\frac{1}{2}(a+1)]} - \frac{[a-2][\frac{1}{2}(a+5)]}{[a][\frac{1}{2}(a-1)]} \right).
\end{equation}
Writing this expression over a common denominator, and using (\ref{eqn:fusion_rules_for_quantum_numbers}), we can write the numerator as
\begin{eqnarray*}
\lefteqn{ [2][a]([2]+[4]+\cdots+[a-1]) - [a+2]([3]+[5]+\cdots+[a-4]) } \\
& & - [a-2]([3]+[5]+\cdots+[a+2]) \\
& = & [a]([1]+[3]+[3]+[5]+\cdots+[a-2]+[a]) \\
& & - ([a+2]+[a-2])([3]+[5]+\cdots+[a-4]) \\
& & - [a-2]([a-2]+[a]+[a+2]) \\
& = & [a] + (2[a]-[a+2]-[a-2])([3]+[5]+\cdots+[a-4]+[a-2]) \\
& & + [a]^2 - [a-2][a] \\
& = & [a] + ([a]-[a+2])([3]+[5]+\cdots+[a-2]) \\
& & + ([a]-[a-2])([3]+[5]+\cdots+[a-2]+[a]) \\
& = & [a] + [(a-3)/2][(a+1)/2]([a]-[a+2]) \\
& & + [(a-1)/2][(a+3)/2]([a]-[a-2]).
\end{eqnarray*}
Now
\begin{eqnarray*}
\lefteqn{ [(a-3)/2][(a+1)/2]([a]-[a+2]) } \\
& \qquad = & [(a-3)/2]([(a+1)/2] + [(a+5)/2] + \cdots + [(3a-1)/2] \\
& & \qquad \qquad \quad - [(a+5)/2] - [(a+9)/2] - \cdots - [(3a+3)/2]) \\
& \qquad = & [(a-3)/2]([(a+1)/2] - [(3a+3)/2]) \\
& \qquad = & [3]+[5]+\cdots+[a-2] -[a+4]-[a+6]-\cdots-[2a-1],
\end{eqnarray*}
and
\begin{eqnarray*}
\lefteqn{ [(a-1)/2][(a+3)/2]([a]-[a-2]) } \\
& \qquad = & [(a-1)/2]([(a+1)/2] + [(a+3)/2] + \cdots + [(3a+1)/2] \\
& & \qquad \qquad \quad - [(a-5)/2] - [(a-1)/2] - \cdots - [(3a-3)/2]) \\
& \qquad = & [(a-1)/2]([(3a+1)/2] - [(a-5)/2]) \\
& \qquad = & [a+2]+[a+4]+\cdots+[2a-1] -[3]-[5]-\cdots-[a-4].
\end{eqnarray*}
Then we find that the numerator is given by $[a] + [a-2] + [a+2] = [3][a]$, and (\ref{eqn:Boltzmann_weight_for_A2n+1(star)}) becomes
$$\frac{[3][a]}{[3][a][\frac{1}{2}(a-1)][\frac{1}{2}(a+1)]} = \frac{1}{[\frac{1}{2}(a-1)][\frac{1}{2}(a+1)]}$$
as required.
To show equivalence, we need unitaries $u_{i,j} \in \mathbb{T}$, for vertices $i$, $j$ of $\mathcal{A}^{(n)\ast}$ such that
$$1 = u_{i,i}u_{i+1,i+1}, \qquad 1 = u_{i,i}u_{i-1,i-1}, \qquad -1 = u_{i,i-1} u_{i-1,i} \overline{u_{i,i+1} u_{i+1,i}},$$
$$(-1)^i = u_{i,i}^2 \overline{u_{i,i+1} u_{i+1,i}}, \qquad  (-1)^{i+1} = u_{i,i}^2 \overline{u_{i,i-1} u_{i-1,i}}.$$
Then we set $u_{i,i} = 1$ for all $i$, and for $m = 0, \ldots, (n-2)/2$, $u_{2m+1,2m} = u_{2m,2m+1} = u_{2m+2,2m+1} = 1$ and $u_{2m+1,2m+2} = -1$.
\hfill
$\Box$

For the graphs $\mathcal{A}^{(4n)\ast}$ (illustrated in Figure \ref{fig:A(star)-graphs}) the Perron-Frobenius weights on the vertices are given by $\phi_{i} = [2i]/[2]$, $i=1,\ldots,2n-1$. There are now two solutions $W^+$, $W^-$ for the cells for $\mathcal{A}^{(4n)\ast}$, which are not equivalent since $|W^+| \neq |W^-|$ and the graph $\mathcal{A}^{(4n)\ast}$ does not contain any multiple edges.

\begin{Thm}
The cells for $\mathcal{A}^{(4n)\ast}$, $n < \infty$, are given, up to equivalence, by the inequivalent solutions $W^+$, $W^-$:\\
$$W^{\pm}_{i,i,i+1} = \frac{\sqrt{[2i][2i+2]}}{[2]\sqrt{[2i+1]}}\sqrt{[2i] \mp [1]}, \qquad \qquad i=1,\ldots,2n-2,$$
$$W^{\pm}_{i,i+1,i+1} = \frac{\sqrt{[2i][2i+2]}}{[2]\sqrt{[2i+1]}}\sqrt{[2i+2] \pm [1]}, \qquad \qquad i=1,\ldots,2n-2,$$
$$W^{\pm}_{i,i,i} = \left\{
\begin{array}{l}
(-1)^{i+1} \frac{\displaystyle \sqrt{[2i]}}{\displaystyle [2]\sqrt{[2i-1][2i+1]}}\sqrt{[2][2i] \pm [4i]}, \qquad i = 1, \ldots, n-1, \\
(-1)^{n+1} \frac{\displaystyle [2n]}{\displaystyle \sqrt{[2][2n-1][2n+1]}}, \qquad \qquad \qquad \quad i = n, \\
(-1)^{i+1} \frac{\displaystyle \sqrt{[2i]}}{\displaystyle [2]\sqrt{[2i-1][2i+1]}}\sqrt{[2][2i] \mp [8n-4i]}, \\ \hfill i = n+1 , \ldots, 2n-1.
\end{array} \right. $$
\end{Thm}
\emph{Proof:}
The proof follows in a similar way to the $\mathcal{A}^{(2n+1)\ast}$ case.
\hfill
$\Box$

For the graphs $\mathcal{A}^{(4n+2)\ast}$ (illustrated in Figure \ref{fig:A(star)-graphs}) the Perron-Frobenius weights on the vertices are again given by $\phi_{i} = [2i]/[2]$, $i=1,\ldots,2n$. There are again two inequivalent solutions $W^+$, $W^-$ for the cells of $\mathcal{A}^{(4n+2)\ast}$.

\begin{Thm}
The cells for $\mathcal{A}^{(4n+2)\ast}$, $n < \infty$, are given, up to equivalence, by the inequivalent solutions $W^+$, $W^-$:\\
$$W^{\pm}_{i,i,i+1} = \frac{\sqrt{[2i][2i+2]}}{[2]\sqrt{[2i+1]}}\sqrt{[2i] \mp [1]}, \qquad \qquad i=1,\ldots,2n-1,$$
$$W^{\pm}_{i,i+1,i+1} = \frac{\sqrt{[2i][2i+2]}}{[2]\sqrt{[2i+1]}}\sqrt{[2i+2] \pm [1]}, \qquad \qquad i=1,\ldots,2n-1,$$
$$W^{\pm}_{i,i,i} = \left\{
\begin{array}{l}
(-1)^{i+1} \frac{\displaystyle \sqrt{[2i]}}{\displaystyle [2]\sqrt{[2i-1][2i+1]}}\sqrt{[2][2i] \pm [4i]}, \qquad \qquad i = 1, \ldots, n, \\
(-1)^{i+1} \frac{\displaystyle \sqrt{[2i]}}{\displaystyle [2]\sqrt{[2i-1][2i+1]}}\sqrt{[2][2i] \mp [8n+4-4i]}, \\
\hfill i = n+1 , \ldots, 2n.
\end{array} \right. $$
\end{Thm}
\emph{Proof:}
The proof again follows in a similar way to the $\mathcal{A}^{(2n+1)\ast}$ case.
\hfill
$\Box$

For $\mathcal{A}^{(2n)\ast}$, the cells $W^+$ above give the following representation of the Hecke algebra:
\begin{eqnarray*}
U^{(i,i+1)} & = & \begin{array}{c} \scriptstyle i \\ \scriptstyle i+1 \end{array} \left( {\begin{array}{cc}
                 \frac{[2i]-[1]}{[2i+1]} & \frac{\sqrt{([2i]-[1])([2i+2]+[1])}}{[2i+1]} \\
                 \frac{\sqrt{([2i]-[1])([2i+2]+[1])}}{[2i+1]} & \frac{[2i+2]+[1]}{[2i+1]}
               \end{array} } \right), \\
U^{(i,i-1)} & = & \begin{array}{c} \scriptstyle i-1 \\ \scriptstyle i \end{array} \left( {\begin{array}{cc}
                 \frac{[2i-2]-[1]}{[2i-1]} & \frac{\sqrt{([2i-2]-[1])([2i]+[1])}}{[2i-1]} \\
                 \frac{\sqrt{([2i-2]-[1])([2i]+[1])}}{[2i-1]} & \frac{[2i]+[1]}{[2i-1]}
               \end{array} } \right), \\
\end{eqnarray*}
\begin{eqnarray*}
\lefteqn{ U^{(i,i)} } \\
& = & \begin{array}{c} \scriptstyle i-1 \\ \scriptstyle i \\ \scriptstyle i+1 \end{array} \left( {\begin{array}{ccc}
                 \frac{[2i-2]([2i]+[1])}{[2i][2i+1]} & (-1)^{i+1}\sqrt{x \, a_+} & \frac{\sqrt{[2i-2][2i-1][2i+2]}}{[2i]\sqrt{[2i+1]}} \\
                 (-1)^{i+1} \sqrt{x \, a_+} & x & (-1)^{i+1}\sqrt{x \, a_-} \\
                 \frac{\sqrt{[2i-2][2i-1][2i+2]}}{[2i]\sqrt{[2i+1]}} & (-1)^{i+1}\sqrt{x \, a_-} & \frac{[2i+2]([2i]-[1])}{[2i][2i+1]}
               \end{array} } \right),
\end{eqnarray*}
where, $a_{\pm} = [2i \mp 2]([2i] \pm [1])/[2i][2i+1]$, and for $m>0$, if $n=2m$,
$$x= \left\{ {\begin{array}{cl}
                 \frac{[2][2i]+[4i]}{[2i-1][2i][2i+1]} & \textrm{ for } i=1,\ldots,m-1, \\
                 \frac{[2]}{[2m-1]^2} & \textrm{ for } i=m, \\
                 \frac{[2][2i]-[4n-4i]}{[2i-1][2i][2i+1]} & \textrm{ for } i=m+1,\ldots,2m-1,
               \end{array} } \right. ,$$
and if $n=2m+1$,
$$x= \left\{ {\begin{array}{cl}
                 \frac{[2][2i]+[4i]}{[2i-1][2i][2i+1]} & \textrm{ for } i=1,\ldots,m, \\
                 \frac{[2][2i]-[4n-4i]}{[2i-1][2i][2i+1]} & \textrm{ for } i=m+1,\ldots,2m,
               \end{array} } \right. .$$

\begin{Lemma}
The weights in the representation of the Hecke algebra given above for $\mathcal{A}^{(2n)\ast}$ are equivalent to the Boltzmann weights at criticality given by Behrend-Evans in \cite{behrend/evans:2004}.
\end{Lemma}
\emph{Proof:}
To make our notation the same as that of \cite{behrend/evans:2004} one replaces $i$ with $a/2$. To see that the absolute values of our weights are equal to those of the Boltzmann weights in \cite{behrend/evans:2004} one needs the following relations on the quantum numbers:
$$[2i]+[1] = \frac{[2i+1]_{q'}[4i+2]_{q'}}{[2i-1]_{q'}}, \qquad [2i]-[1] = \frac{[2i-1]_{q'}[4i+2]_{q'}}{[2i+1]_{q'}},$$
where $q' = \sqrt{q}$ ($q=e^{i \pi/n}$). Again, a bit more work is required for $[U^{(i,i)}]_{i,i}$. For equivalence we make the same choice of $(u_{i,j})_{i,j}$ as for $\mathcal{A}^{(2n+1)\ast}$.
\hfill
$\Box$

\section{$\mathcal{D}^{\ast}$ graphs}

\begin{figure}[tb]
\begin{center}
\includegraphics[width=115mm]{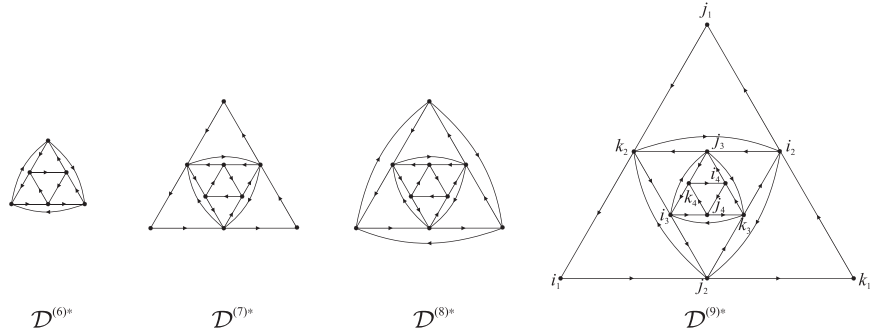}\\
 \caption{$\mathcal{D}^{(n)\ast}$ for $n = 6,7,8,9$} \label{fig:Weights11}
\end{center}
\end{figure}

The graphs $\mathcal{D}^{(n)\ast}$ are illustrated in Figure \ref{fig:Weights11}. We label its vertices by $i_l$, $j_l$ and $k_l$, $l=1, \ldots, \lfloor (n-1)/2 \rfloor$, which we have illustrated in Figure \ref{fig:Weights11} for $n=9$.

We consider first the graphs $\mathcal{D}^{(2n+1)\ast}$. The Perron-Frobenius weights are $\phi_{i_l} = \phi_{j_l} = \phi_{k_l} = [2l-1]$, $l=1, \ldots, n$. Since the graph has a $\mathbb{Z}_3$ symmetry, we will seek $\mathbb{Z}_3$-symmetric solutions (up to choice of phase), i.e. $|W_{i_p,j_q,k_r}|^2 = |W_{i_q,j_r,k_p}|^2 = |W_{i_r,j_p,k_q}|^2 =: |W_{p,q,r}|^2$, $p,q,r \in \{ 1, \ldots, n \}$. Using this notation, we have the following equations from type I frames:
\begin{eqnarray}
& & |W_{1,2,2}|^2 = [2][3], \label{eqn:D(2n+1star)_Weights-I(i)}\\
& & |W_{l,l,l+1}|^2 + |W_{l,l+1,l+1}|^2 = [2][2l-1][2l+1], \qquad \qquad l = 2, \ldots, n-1, \label{eqn:D(2n+1star)_Weights-I(ii)}\\
& & |W_{l-1,l,l}|^2 + |W_{l,l,l}|^2 + |W_{l,l,l+1}|^2 = [2][2l-1]^2, \qquad \qquad l = 2, \ldots, n-1, \label{eqn:D(2n+1star)_Weights-I(iii)}\\
& & |W_{n-1,n,n}|^2 + |W_{n,n,n}|^2 = [2]^3, \label{eqn:D(2n+1star)_Weights-I(iv)}
\end{eqnarray}
and from type II frames we have:
\begin{equation}
|W_{l-1,l,l}|^2 |W_{l,l,l+1}|^2 = [2l-3][2l-1]^2[2l+1], \label{eqn:D(2n+1star)_Weights-II(i)}
\end{equation}
for $l = 2, \ldots, n-1$, and
\begin{equation}
|W_{l-1,l,l}|^2(\frac{1}{[2l-3]} |W_{l-1,l-1,l}|^2 + \frac{1}{[2l-1]} |W_{l,l,l}|^2) = [2l-3][2l-1]^2, \label{eqn:D(2n+1star)_Weights-II(ii)}
\end{equation}
for $l = 2, \ldots, n$,
which are exactly those for the type I and type II frames for the graph $\mathcal{A}^{(2n+1)\ast}$. Since the Perron-Frobenius weights and Coxeter number are also the same as for $\mathcal{A}^{(2n+1)\ast}$, the cells $|W_{p,q,r,}|$ follow.

From the type II frame consisting of the vertices $i_l$, $j_l$, $i_{l+1}$ and $j_{l+1}$ we have the following restriction on the choice of phase
\begin{eqnarray} \label{eqn:restriction_on_lambdas-D(2n+1star)}
\lefteqn{ \lambda_{i_l,j_l,k_{l+1}} \lambda_{i_l,j_{l+1},k_l} \lambda_{i_{l+1},j_l,k_l} \lambda_{i_{l+1},j_{l+1},k_{l+1}} } \\
& = & - \lambda_{i_l,j_l,k_l} \lambda_{i_l,j_{l+1},k_{l+1}} \lambda_{i_{l+1},j_l,k_{l+1}} \lambda_{i_{l+1},j_{l+1},k_l}. \nonumber
\end{eqnarray}

\begin{Thm}
Every $\mathbb{Z}_3$-symmetric solution for the cells $W$ of $\mathcal{D}^{(2n+1)\ast}$, $n < \infty$, is equivalent to the solution
$$W_{i_1,j_2,k_2} = W_{i_2,j_1,k_2} = W_{i_2,j_2,k_1} = \sqrt{[2][3]},$$
$$W_{i_l,j_{l+1},k_{l+1}} = W_{i_{l+1},j_l,k_{l+1}} = W_{i_{l+1},j_{l+1},k_l} = \frac{\sqrt{[l+1][2l-1][2l+1]}}{\sqrt{[l]}},$$
$$W_{i_l,j_l,k_{l+1}} = W_{i_l,j_{l+1},k_l} = W_{i_{l+1},j_l,k_l} = \frac{\sqrt{[l-1][2l-1][2l+1]}}{\sqrt{[l]}},$$
$$W_{i_l,j_l,k_l} = (-1)^{l+1} \frac{[2l-1]}{\sqrt{[l-1][l]}}, \qquad W_{i_n,j_n,k_n} = (-1)^{n+1} \frac{[2n-1]}{\sqrt{[n-1][n]}},$$
for $l=2,\ldots,n-1.$
\end{Thm}
\emph{Proof:}
Let $W^{\sharp}$ be any $\mathbb{Z}_3$-symmetric solution for the cells of $\mathcal{D}^{(2n+1)\ast}$, where the choice of phase satisfies the condition (\ref{eqn:restriction_on_lambdas-D(2n+1star)}). Since $\mathcal{D}^{(2n+1)\ast}$ does not contain any multiple edges, we must have $|W^{\sharp}_{ijk}| = |W_{ijk}|$ for every triangle $\triangle_{ijk}$ of $\mathcal{D}^{(2n+1)\ast}$. We need to find a family of unitaries $\{ u_{p,q} \}$, where $u_{p,q}$ is the unitary for the edge from vertex $p$ to vertex $q$ on $\mathcal{D}^{(2n+1)\ast}$, which satisfy (\ref{eqn:equvialence_of_W1,W2(lambdas)-no_multiple_edges}), i.e. $-1 = u_{i_{2l},j_{2l}} u_{j_{2l},k_{2l}} u_{k_{2l},i_{2l}} \lambda^{\sharp}_{i_{2l},j_{2l},k_{2l}}$ for the triangle $\triangle_{i_{2l},j_{2l},k_{2l}}$, $l=1,\ldots, \lfloor n/2 \rfloor$, and $1 = u_{p_1} u_{p_2} u_{p_3} \lambda_{p_1,p_2,p_3}$ for all other triangles on $\mathcal{D}^{(2n+1)\ast}$.
For triangles involving the outermost vertices, we require $1 = u_{i_1,j_2} u_{j_2,k_2} u_{k_2,i_1} \lambda^{\sharp}_{i_1,j_2,k_1}$, $1 = u_{i_2,j_1} u_{j_1,k_2} u_{k_2,i_2} \lambda^{\sharp}_{i_2,j_1,k_2}$, $1 = u_{i_2,j_2} u_{j_2,k_1} u_{k_1,i_2} \lambda^{\sharp}_{i_2,j_2,k_1}$ and $-1 = u_{i_2,j_2} u_{j_2,k_2} u_{k_2,i_2} \lambda^{\sharp}_{i_2,j_2,k_2}$. So we choose $u_{i_1,j_2} = u_{j_1,k_2} = u_{k_1,i_2} = u_{j_2,k_2} = u_{k_2,i_2} = 1$, $u_{i_2,j_1} = \overline{\lambda^{\sharp}_{i_2,j_1,k_2}}$, $u_{k_2,i_1} = \overline{\lambda^{\sharp}_{i_1,j_2,k_2}}$, $u_{i_2,j_2} = - \overline{\lambda^{\sharp}_{i_2,j_2,k_2}}$ and $u_{j_2,k_1} = - \lambda^{\sharp}_{i_2,j_2,k_2} \overline{\lambda^{\sharp}_{i_2,j_2,k_1}}$. Next consider the equations $1 = u_{i_2,j_3} u_{j_3,k_2} u_{k_2,i_2} \lambda^{\sharp}_{i_2,j_3,k_2}$, $1 = u_{i_3,j_2} u_{j_2,k_2} u_{k_2,i_3} \lambda^{\sharp}_{i_3,j_2,k_2}$ and $1 = u_{i_2,j_2} u_{j_2,k_3} u_{k_3,i_2} \lambda^{\sharp}_{i_2,j_2,k_3}$. We make the following choices: $u_{i_2,j_3} = u_{j_2,k_3} = u_{k_2,i_3} = 1$, $u_{i_3,j_2} = \overline{\lambda^{\sharp}_{i_3,j_2,k_2}}$, $u_{j_3,k_2} = \overline{\lambda^{\sharp}_{i_2,j_3,k_2}}$ and $u_{k_3,i_2} = - \lambda^{\sharp}_{i_2,j_2,k_2} \overline{\lambda^{\sharp}_{i_2,j_2,k_3}}$. Next we consider the equations
$$1 = u_{i_2,j_3} u_{j_3,k_3} u_{k_3,i_2} \lambda^{\sharp}_{i_2,j_3,k_3} = - u_{j_3,k_3} \lambda^{\sharp}_{i_2,j_2,k_2} \overline{\lambda^{\sharp}_{i_2,j_2,k_3}} \lambda^{\sharp}_{i_2,j_3,k_3},$$
$$1 = u_{i_3,j_2} u_{j_2,k_3} u_{k_3,i_3} \lambda^{\sharp}_{i_3,j_2,k_3} =  u_{k_3,i_3} \overline{\lambda^{\sharp}_{i_3,j_2,k_2}} \lambda^{\sharp}_{i_3,j_2,k_3},$$
$$1 = u_{i_3,j_3} u_{j_3,k_2} u_{k_2,i_3} \lambda^{\sharp}_{i_3,j_3,k_2} = u_{i_3,j_3} \overline{\lambda^{\sharp}_{i_2,j_3,k_2}} \lambda^{\sharp}_{i_3,j_3,k_2}.$$

We make the choices $u_{i_3,j_3} = \lambda^{\sharp}_{i_2,j_3,k_2} \overline{\lambda^{\sharp}_{i_3,j_3,k_2}}$, $u_{k_3,i_3} = \lambda^{\sharp}_{i_3,j_2,k_2} \overline{\lambda^{\sharp}_{i_3,j_2,k_3}}$ and $u_{j_3,k_3} = - \lambda^{\sharp}_{i_2,j_2,k_3} \overline{\lambda^{\sharp}_{i_2,j_2,k_2}} \overline{\lambda^{\sharp}_{i_2,j_3,k_3}}$. Then
\begin{eqnarray*}
\lefteqn{ u_{i_3,j_3} u_{j_3,k_3} u_{k_3,i_3} \lambda^{\sharp}_{i_3,j_3,k_3} } \\
& = & - \lambda^{\sharp}_{i_2,j_3,k_2} \overline{\lambda^{\sharp}_{i_3,j_3,k_2}} \lambda^{\sharp}_{i_2,j_2,k_3} \overline{\lambda^{\sharp}_{i_2,j_2,k_2}} \overline{\lambda^{\sharp}_{i_2,j_3,k_3}} \lambda^{\sharp}_{i_3,j_2,k_2} \overline{\lambda^{\sharp}_{i_3,j_2,k_3}} = -1,
\end{eqnarray*}
by (\ref{eqn:restriction_on_lambdas-D(2n+1star)}), as required. Continuing in this way we are done.
\hfill
$\Box$

For $\mathcal{D}^{(2n+1)\ast}$, the Hecke representation for the cells $W$ above is given by the Hecke representation for $\mathcal{A}^{(2n+1)\ast}$, where
$[U^{(i_l,k_r)}]_{j_m,j_p} = [U^{(j_l,i_r)}]_{k_m,k_p} = [U^{(k_l,j_r)}]_{i_m,i_p}$
are given by the weights $[U^{(l,r)}]_{m,p}$ for $\mathcal{A}^{(2n+1)\ast}$, for any $l,m,p,r$ allowed by the graph.

We now consider the graphs $\mathcal{D}^{(2n)\ast}$. The Perron-Frobenius weights are $\phi_{i_l} = \phi_{j_l} = \phi_{k_l} = [2l]/[2]$, and we again assume $|W_{i_p,j_q,k_r}|^2 = |W_{i_q,j_r,k_p}|^2 = |W_{i_r,j_p,k_q}|^2 =: |W_{p,q,r}|^2$, where $p,q,r \in \{ 1, \ldots, n-1 \}$. Then as for $\mathcal{D}^{(2n+1)\ast}$, the $\mathbb{Z}_3$-symmetric solution for the cells follows from the solution for $\mathcal{A}^{(2n)\ast}$, and we have the same restriction (\ref{eqn:restriction_on_lambdas-D(2n+1star)}) on the choice of phase. So we have

\begin{Thm}
For $n < \infty$, the $\mathbb{Z}_3$-symmetric solution for the cells of $\mathcal{D}^{(4n)\ast}$ are given by
$$W^{\pm}_{i_l,j_l,k_{l+1}} = W^{\pm}_{i_l,j_{l+1},k_l} = W^{\pm}_{i_{l+1},j_l,k_l} = \frac{\sqrt{[2l][2l+2]}}{[2]\sqrt{[2l+1]}}\sqrt{[2l] \mp [1]},$$
$$\hspace{9cm} l=2,\ldots,2n-2,$$
$$W^{\pm}_{i_l,j_{l+1},k_{l+1}} = W^{\pm}_{i_{l+1},j_l,k_{l+1}} = W^{\pm}_{i_{l+1},j_{l+1},k_l} = \frac{\sqrt{[2l][2l+2]}}{[2]\sqrt{[2l+1]}}\sqrt{[2l+2] \pm [1]},$$
$$\hspace{9cm} l=1,\ldots,2n-2,$$
$$W^{\pm}_{i_l,j_l,k_l} = \left\{
\begin{array}{l}
(-1)^{l+1} \frac{\displaystyle \sqrt{[2l]}}{\displaystyle [2]\sqrt{[2l-1][2l+1]}}\sqrt{[2][2l] \pm [4l]}, \quad \;\; l = 1, \ldots, n-1, \\
(-1)^{n+1} \frac{\displaystyle [2n]}{\displaystyle \sqrt{[2][2n-1][2n+1]}}, \hfill l = n, \\
(-1)^{l+1} \frac{\displaystyle \sqrt{[2l]}}{\displaystyle [2]\sqrt{[2l-1][2l+1]}}\sqrt{[2][2l] \mp [8n-4l]}, \\
\hfill l = n+1 , \ldots, 2n-1,
\end{array} \right. $$
and the $\mathbb{Z}_3$-symmetric solution for the cells of $\mathcal{D}^{(4n+2)\ast}$ are \\
$$W^{\pm}_{i_l,j_l,k_{l+1}} = W^{\pm}_{i_l,j_{l+1},k_l} = W^{\pm}_{i_{l+1},j_l,k_l} = \frac{\sqrt{[2l][2l+2]}}{[2]\sqrt{[2l+1]}}\sqrt{[2l] \mp [1]},$$
$$\hspace{9cm} l=2,\ldots,2n-1,$$
$$W^{\pm}_{i_l,j_{l+1},k_{l+1}} = W^{\pm}_{i_{l+1},j_l,k_{l+1}} = W^{\pm}_{i_{l+1},j_{l+1},k_l} = \frac{\sqrt{[2l][2l+2]}}{[2]\sqrt{[2l+1]}}\sqrt{[2l+2] \pm [1]},$$
$$\hspace{9cm} l=1,\ldots,2n-1,$$
$$W^{\pm}_{i_l,j_l,k_l} = \left\{
\begin{array}{l}
(-1)^{l+1} \frac{\displaystyle \sqrt{[2l]}}{\displaystyle [2]\sqrt{[2l-1][2l+1]}}\sqrt{[2][2l] \pm [4l]}, \qquad \;\; l = 1, \ldots, n, \\
(-1)^{l+1} \frac{\displaystyle \sqrt{[2l]}}{\displaystyle [2]\sqrt{[2l-1][2l+1]}}\sqrt{[2][2l] \mp [8n+4-4l]}, \\
\hfill l = n+1 , \ldots, 2n.
\end{array} \right. $$
\end{Thm}

The uniqueness of these solutions follows in the same way as for $\mathcal{D}^{(2n+1)\ast}$. If $W^+$ is a solution for the cells of $\mathcal{D}^{(2n)\ast}$, then $W^-$ is a solution for the cells of the graph where we switch vertices $i_l \leftrightarrow i_{n-l}$, $j_l \leftrightarrow j_{n-l}$ and $k_l \leftrightarrow k_{n-l}$, for all $l=1,\ldots,n-1$.

For $\mathcal{D}^{(2n)\ast}$, the Hecke representation for the cells $W^+$ above is given by the Hecke representation for $\mathcal{A}^{(2n)\ast}$, where
$[U^{(i_l,k_r)}]_{j_m,j_p} = [U^{(j_l,i_r)}]_{k_m,k_p} = [U^{(k_l,j_r)}]_{i_m,i_p}$
are given by the weights $[U^{(l,r)}]_{m,p}$ for $\mathcal{A}^{(2n)\ast}$, for any $l,m,p,r$ allowed by the graph.

In \cite{di_francesco/zuber:1990}, di Francesco and Zuber gave a representation of the Hecke algebra for the graph $\mathcal{D}^{(6)\ast}$, with the absolute values of the weights there equal to those for our weights given above. The two Hecke representations are not identical as the weights in \cite{di_francesco/zuber:1990} involve the complex variable $i$. However it has not been possible to determine whether or not the two representations are equivalent as there are known to be a number of typographical errors in the representation in \cite{di_francesco/zuber:1990}.

\section{$\mathcal{E}^{(8)}$}

We will label the vertices of the exceptional graph $\mathcal{E}^{(8)}$ in the following way. We will label the six outmost vertices by $i_l$ and the six inmost vertices by $j_l$, $l=1, \ldots, 6$, such that there are edges from $i_l$ to $j_l$ and from $j_l$ to $i_{l+1}$. The Perron-Frobenius weights on the vertices are $\phi_{i_l} = 1$, $\phi_{j_l} = [3]$. With $[a]=[a]_q$, $q=e^{i \pi/8}$, we have $[4]/[2] = \sqrt{2}$.

\begin{figure}[tb]
\begin{center}
\includegraphics[width=70mm]{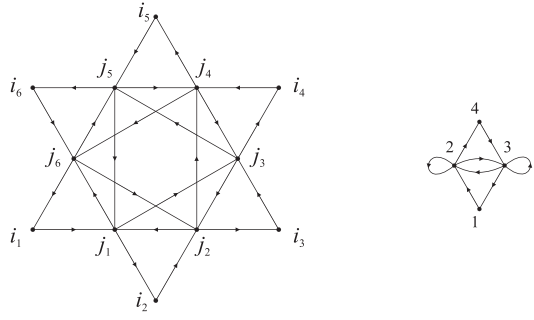}\\
 \caption{$\mathcal{E}^{(8)}$ and its $\mathbb{Z}_3$ orbifold $\mathcal{E}^{(8)\ast}$}\label{fig:E(8)&E(8)(star)-graphs}
\end{center}
\end{figure}

We will again use the notation $W_{i,j,k}$ for $W(\triangle_{i,j,k})$. Then from the type I frames on the graph we have the following equations:
$$|W_{i_l,j_l,j_{l-1}}|^2 = [2] \phi_{i_l} \phi_{j_l} = [2][3],$$
$$|W_{i_l,j_l,j_{l-1}}|^2 + |W_{j_{l+1},j_l,j_{l-1}}|^2 + |W_{j_l,j_{l-1},j_{l-2}}|^2 = [2] \phi_{j_l} \phi_{j_{l-1}} = [2][3]^2.$$
Then $|W_{j_{l+1},j_l,j_{l-1}}|^2 + |W_{j_l,j_{l-1},j_{l-2}}|^2 = [3][4]$. Since there is a $\mathbb{Z}_6$ symmetry of $\mathcal{E}^{(8)}$ we assume $|W_{j_{l+1},j_l,j_{l-1}}|^2 = |W_{j_{k+1},j_k,j_{k-1}}|^2$ for all $k,l$, giving
$$|W_{j_{l+1},j_l,j_{l-1}}|^2 = \frac{1}{2} [3][4] = \frac{[2]^2[3]}{[4]}.$$
The $\mathbb{Z}_6$ symmetry of the cells can be deduced from equation (\ref{eqn:eqn_for_a_typeIIframe_in_E(8)}). Finally, for the type I frames $\stackrel{j_l}{\bullet} \rightarrow \stackrel{j_{l+2}}{\bullet}$ we have $|W_{j_{l+2},j_{l+1},j_l}|^2 + |W_{j_l,j_{l+2},j_{l+4}}|^2 = [2][3]^2$ giving
$$|W_{j_l,j_{l+2},j_{l+4}}|^2 = [2][3]^2 - \frac{[2]^2[3]}{[4]} = \frac{[2]^2[3]^2}{[4]}.$$
Let
\begin{eqnarray*}
W_{i_l,j_l,j_{l-1}} & = & \lambda_{i_l} \sqrt{[2][3]}, \quad l=1, \ldots, 6,\\
W_{j_l,j_{l-1},j_{l-2}} & = & \lambda_{j_l}^{(1)} \frac{[2]\sqrt{[3]}}{\sqrt{[4]}}, \quad l=1, \ldots, 6,\\
W_{j_l,j_{l+2},j_{l+4}} & = & \lambda_{j_l}^{(2)} \frac{[2][3]}{\sqrt{[4]}}, \quad l=1,2.
\end{eqnarray*}

The only type II frames that yield anything new are those for the frame involving the vertices $j_{l-2}$, $j_{l-3}$($= j_{l+3}$), $j_{l+1}$ and $j_l$:
\begin{eqnarray}
0 & = & \phi_{j_{l-1}}^{-1} W_{j_{l-2},j_{l-1},j_l} \overline{W_{j_{l+1},j_l,j_{l-1}}} W_{j_{l-1},j_{l+1},j_{l+3}} \overline{W_{j_{l-1},j_{l-2},j_{l-3}}} \nonumber \\
& & \quad + \phi_{j_{l+2}}^{-1} W_{j_{l-2},j_l,j_{l+2}} \overline{W_{j_{l+2},j_{l+1},j_l}} W_{j_{l+3},j_{l+2},j_{l+1}} \overline{W_{j_{l-2},j_{l-3},j_{l+2}}} \nonumber \\
& = & \frac{[2]^4\sqrt{[3]^3}}{[4]^2} \lambda_{j_l}^{(1)} \lambda_{j_{l+2}}^{(1)} \lambda_{j_{l+4}}^{(1)} \lambda_{j_{l-1}}^{(2)} + \frac{[2]^4\sqrt{[3]^3}}{[4]^2} \lambda_{j_{l-1}}^{(1)} \lambda_{j_{l+1}}^{(1)} \lambda_{j_{l+3}}^{(1)} \lambda_{j_l}^{(2)}, \label{eqn:eqn_for_a_typeIIframe_in_E(8)}
\end{eqnarray}
which for any $l = 1, \ldots, 6$ gives
\begin{equation} \label{eqn:restriction_on_lambdas_for_Weights_E(8)}
\lambda_{j_1}^{(1)} \lambda_{j_3}^{(1)} \lambda_{j_5}^{(1)} \lambda_{j_2}^{(2)} = - \lambda_{j_2}^{(1)} \lambda_{j_4}^{(1)} \lambda_{j_6}^{(1)} \lambda_{j_1}^{(2)}.
\end{equation}
From the type II frame above we see that there must be a $\mathbb{Z}_6$ symmetry on the cells, $|W_{j_{l+1},j_l,j_{l-1}}|^2 = |W_{j_{k+1},j_k,j_{k-1}}|^2$ for all $k,l$, is correct since otherwise the coefficients of the two terms in equation (\ref{eqn:eqn_for_a_typeIIframe_in_E(8)}) would be different, and (\ref{eqn:restriction_on_lambdas_for_Weights_E(8)}) would be
$$\lambda_{j_1}^{(1)} \lambda_{j_3}^{(1)} \lambda_{j_5}^{(1)} \lambda_{j_2}^{(2)} = - c \lambda_{j_2}^{(1)} \lambda_{j_4}^{(1)} \lambda_{j_6}^{(1)} \lambda_{j_1}^{(2)},$$
for some constant $c \in \mathbb{R}$ with $|c| \neq 1$, which is impossible.

\begin{Thm}
There is up to equivalence a unique set of cells for $\mathcal{E}^{(8)}$ given by \\
\begin{eqnarray*}
W_{i_l,j_l,j_{l-1}} = \sqrt{[2][3]}, & & \qquad W_{j_l,j_{l-1},j_{l-2}} = \frac{[2]\sqrt{[3]}}{\sqrt{[4]}}, \quad l=1, \ldots, 6,\\
W_{j_1,j_3,j_5} = \frac{[2][3]}{\sqrt{[4]}}, & & \qquad W_{j_2,j_4,j_6} = - \frac{[2][3]}{\sqrt{[4]}}.
\end{eqnarray*}
\end{Thm}
\emph{Proof:}
Let $W^{\sharp}$ be any solution for the for the cells for $\mathcal{E}^{(8)}$, where the choice of phase satisfies the condition (\ref{eqn:restriction_on_lambdas_for_Weights_E(8)}). We need to find a family of unitaries $\{ u_{p,q} \}$, where $u_{p,q}$ is the unitary for the edge from vertex $p$ to vertex $q$ on $\mathcal{E}^{(8)}$, which satisfy (\ref{eqn:equvialence_of_W1,W2(lambdas)-no_multiple_edges}), i.e. $-1 = u_{j_2,j_4} u_{j_4,j_6} u_{j_6,j_2} \lambda_{j_2}^{(2)}$ for the triangle $\triangle_{j_2, j_4, j_6}$, and $1 = u_{p_1} u_{p_2} u_{p_3} \lambda_{p_1,p_2,p_3}$ for all other triangles, where $\lambda_{p_1,p_2,p_3}$ is the phase associated to triangle  $\triangle_{p_1,p_2,p_3}$.
We make the choices $u_{i_l,j_l} = \overline{u_{j_l,j_{l-1}} \lambda_{i_l}}$,
$u_{j_l,j_{l+1}} = 1$ for $l = 1, \ldots, 6$,
$u_{j_2,j_1} = u_{j_5,j_4} = 1$,
$u_{j_1,j_6} = \overline{\lambda_{j_2}^{(1)}}$,
$u_{j_3,j_2} = \lambda_{j_2}^{(1)} \lambda_{j_6}^{(1)} \lambda_{j_1}^{(2)} \overline{\lambda_{j_1}^{(1)} \lambda_{j_3}^{(1)}}$,
$u_{j_4,j_3} = \overline{\lambda_{j_5}^{(1)}}$,
$u_{j_6,j_5} = \overline{\lambda_{j_6}^{(1)}}$,
$u_{j_3,j_5} = u_{j_4,j_6} = u_{j_6,j_2} = 1$,
$u_{j_1,j_3} = \lambda_{j_2}^{(1)} \overline{\lambda_{j_2}^{(1)}} \overline{\lambda_{j_6}^{(1)}} \overline{\lambda_{j_1}^{(2)}}$,
$u_{j_2,j_4} = \lambda_{j_2}^{(1)} \lambda_{j_3}^{(1)} \lambda_{j_5}^{(1)} \overline{\lambda_{j_2}^{(1)}} \overline{\lambda_{j_4}^{(1)}} \overline{\lambda_{j_6}^{(1)}} \overline{\lambda_{j_1}^{(2)}}$
and $u_{j_5,j_1} = \lambda_{j_2}^{(1)} \lambda_{j_6}^{(1)} \overline{\lambda_{j_1}^{(1)}}$.
\hfill
$\Box$

For $\mathcal{E}^{(8)}$, the above cells $W$ give the following representation of the Hecke algebra:
\begin{eqnarray*}
U^{(i_l,j_{l-1})} & = & U^{(j_l,i_l)} \;\; = \;\; [2], \\
U^{(j_l,j_{l-2})} & = & \begin{array}{c} \scriptstyle j_{l-1} \\ \scriptstyle j_{l+2} \end{array} \left( {\begin{array}{cc}
                 \frac{1}{[2]} & \frac{(-1)^{l+1}\sqrt{[3]}}{[2]} \\
                 \frac{(-1)^{l+1}\sqrt{[3]}}{[2]} & \frac{[3]}{[2]}
               \end{array} } \right), \\
U^{(j_l,j_{l+1})} & = & \begin{array}{c} \scriptstyle j_{l-1} \\ \scriptstyle j_{l+2} \\ \scriptstyle i_{l+1} \end{array} \left( {\begin{array}{ccc}
                 \frac{1}{[2]} & \frac{1}{[2]} & \frac{1}{\sqrt{[3]}} \\
                 \frac{1}{[2]} & \frac{1}{[2]} & \frac{1}{\sqrt{[3]}} \\
                 \frac{1}{\sqrt{[3]}} & \frac{1}{\sqrt{[3]}} & \frac{[2]}{[3]}
               \end{array} } \right),
\end{eqnarray*}
for $l=1,\ldots,6$ (mod 6). This representation is identical to that given by di Francesco-Zuber in \cite{di_francesco/zuber:1990}. (The representation in \cite{di_francesco/zuber:1990} is given for the graph $\mathcal{E}^{(8)\ast}$, and the representation for $\mathcal{E}^{(8)}$ is obtained by an unfolding of the graph $\mathcal{E}^{(8)\ast}$.)

\section{$\mathcal{E}^{(8)\ast}$}

We will label the vertices of the graph $\mathcal{E}^{(8)\ast}$ as in Figure \ref{fig:E(8)&E(8)(star)-graphs}. The Perron-Frobenius weights for $\mathcal{E}^{(8)\ast}$ are $\phi_1 = \phi_4 = 1$, $\phi_2 = \phi_3 = [3]$. As with the graphs $\mathcal{A}^{(n)}$ and $\mathcal{E}^{(8)}$ we easily find $|W_{123}|^2 = [2][3]$ and $|W_{234}|^2 = [2][3]$. Then by the type II frame $\stackrel{1}{\bullet} \rightarrow \stackrel{2}{\bullet} \leftarrow \stackrel{2}{\bullet}$ we have $[3]^{-1} |W_{123}|^2 |W_{223}|^2 = [3]^2$, and so $|W_{223}|^2 = [3]^2/[2]$. Similarly $|W_{233}|^2 = [3]^2/[2]$. From the type I frame $\stackrel{2}{\bullet} \rightarrow \stackrel{2}{\bullet}$ we get $|W_{222}|^2 + |W_{223}|^2 = [2][3]^2$, giving $|W_{222}|^2 = [3]^3/[2]$, and similarly $|W_{333}|^2 = [3]^3/[2]$. Let $W_{ijk} = \lambda_{ijk} |W_{ijk}|$. Then from the type II frame consisting of the vertices 2,2,3,3 we obtain the following restriction on the choice of phase:
\begin{equation} \label{eqn:restriction_on_lambdas-E(8star)}
\lambda_{222} \lambda_{233}^3 = - \lambda_{333} \lambda_{223}^3.
\end{equation}

\begin{Thm}
There is up to equivalence a unique set of cells for $\mathcal{E}^{(8)\ast}$ given by
\begin{eqnarray*}
W_{123} & = & W_{234} \;\; = \;\; \sqrt{[2][3]}, \\
W_{223} & = & W_{233} \;\; = \;\; \frac{[3]}{\sqrt{[2]}}, \\
W_{222} & = & \frac{\sqrt{[3]^3}}{\sqrt{[2]}}, \qquad W_{333} \;\; = \;\; - \frac{\sqrt{[3]^3}}{\sqrt{[2]}}.
\end{eqnarray*}
\end{Thm}
\emph{Proof:}
Let $W^{\sharp}$ be any solution for the cells for $\mathcal{E}^{(8)\ast}$, where the choice of phase satisfies the condition (\ref{eqn:restriction_on_lambdas-E(8star)}). We need to find a family of unitaries $\{ u_{p,q} \}$, where $u_{p,q}$ is the unitary for the edge from vertex $p$ to vertex $q$ on $\mathcal{E}^{(8)\ast}$, which satisfy (\ref{eqn:equvialence_of_W1,W2(lambdas)-no_multiple_edges}), i.e. $-1 = u_{3,3}^3 \lambda_{333}$ for the triangle $\triangle_{3,3,3}$, and $1 = u_{i,j} u_{j,k} u_{k,i} \lambda_{ijk}$ for all other triangles, where $\lambda_{ijk}$ is the phase associated to triangle  $\triangle_{i,j,k}$.
We choose $u_{3,1} = u_{3,2} = u_{4,3} = 1$,
$u_{2,4} = \overline{\lambda_{234}}$,
$u_{3,3} = - \overline{\lambda_{333}}^{\frac{1}{3}}$,
$u_{2,3} = - \lambda_{333}^{\frac{1}{3}} \overline{\lambda_{233}}$,
$u_{1,2} = - \lambda_{233} \overline{\lambda_{123}} \overline{\lambda_{333}}^{\frac{1}{3}}$ and
$u_{2,2} = - \lambda_{233} \overline{\lambda_{223}} \overline{\lambda_{333}}^{\frac{1}{3}}$.
\hfill
$\Box$

For $\mathcal{E}^{(8)\ast}$, the above cells $W$ give the following Hecke representation:
\begin{eqnarray*}
U^{(1,3)} & = & U^{(2,1)} \;\; = \;\; U^{(3,4)} \;\; = \;\; U^{(4,2)} \;\; = \;\; [2], \\
U^{(2,2)} & = & \begin{array}{c} \scriptstyle 3 \\ \scriptstyle 2 \end{array} \left( {\begin{array}{cc}
                 \frac{1}{[2]} & \frac{\sqrt{[3]}}{[2]} \\
                 \frac{\sqrt{[3]}}{[2]} & \frac{[3]}{[2]}
               \end{array} } \right), \\
U^{(3,3)} & = & \begin{array}{c} \scriptstyle 2 \\ \scriptstyle 3 \end{array} \left( {\begin{array}{cc}
                 \frac{1}{[2]} & - \frac{\sqrt{[3]}}{[2]} \\
                 - \frac{\sqrt{[3]}}{[2]} & \frac{[3]}{[2]}
               \end{array} } \right), \\
U^{(2,3)} & = & \begin{array}{c} \scriptstyle 2 \\ \scriptstyle 3 \\ \scriptstyle 4 \end{array} \left( {\begin{array}{ccc}
                 \frac{1}{[2]} & \frac{1}{[2]} & \frac{1}{\sqrt{[3]}} \\
                 \frac{1}{[2]} & \frac{1}{[2]} & \frac{1}{\sqrt{[3]}} \\
                 \frac{1}{\sqrt{[3]}} & \frac{1}{\sqrt{[3]}} & \frac{[2]}{[3]}
               \end{array} } \right). \\
\lefteqn{ = \;\; U^{(3,2)} \quad \textrm{ with rows labelled by } 2,3,1. }
\end{eqnarray*}
This representation is identical to that given by di Francesco-Zuber in \cite{di_francesco/zuber:1990}.

\section{$\mathcal{E}_2^{(12)}$}

\begin{figure}[tb]
\begin{center}
\includegraphics[width=100mm]{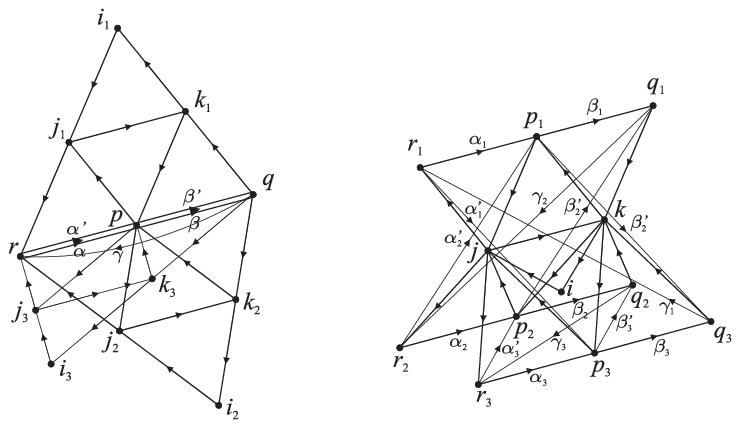}\\
 \caption{$\mathcal{E}_1^{(12)}$ and $\mathcal{E}_2^{(12)}$}\label{fig:labelled_E1(12)&E2(12)graphs}
\end{center}
\end{figure}

We label the vertices and edges of the graph $\mathcal{E}_2^{(12)}$ as in Figure \ref{fig:labelled_E1(12)&E2(12)graphs}. The Perron-Frobenius weights for $\mathcal{E}_2^{(12)}$ are
$$\phi_i = 1, \qquad \phi_j = \phi_k = [3], \qquad \phi_{p_l}=\frac{[2]^3}{[4]}, \qquad \phi_{q_l} = \phi_{r_l} = \frac{[2][3]}{[4]}, \quad l=1,2,3.$$

Let $W_{v_1,v_2,v_3} = \lambda_{v_1,v_2,v_3} |W_{v_1,v_2,v_3}|$ for vertices $v_1$, $v_2$, $v_3$ of $\mathcal{E}_2^{(12)}$. The type II frames consisting of the vertices $p_l$, $k$, $p_{l-1}$ and $r_l$ give a restriction on the phases $\lambda_{v_1,v_2,v_3}$:
\begin{eqnarray*}
0 & = & \phi_{q_{l-1}}^{-1} W_{p_{l-1},q_{l-1},r_l} \overline{W_{p_{l-1},q_{l-1},k}} W_{p_{l},q_{l-1},k} \overline{W_{p_{l},q_{l-1},r_l}} \\
& & \quad + \phi_{j}^{-1} W_{p_{l-1},j,r_l} \overline{W_{p_{l-1},j,k}} W_{p_{l},j,k} \overline{W_{p_{l},j,r_l}} \\
& = & \sqrt{\frac{[2]^9[3]^3}{[4]^5}} \lambda_{p_{l-1},q_{l-1},r_l} \lambda_{p_{l},q_{l-1},k} \overline{\lambda_{p_{l-1},q_{l-1},k}} \overline{\lambda_{p_{l},q_{l-1},r_l}} \\
& & \quad + \sqrt{\frac{[2]^9[3]^3}{[4]^5}} \lambda_{p_{l-1},j,r_l} \lambda_{p_{l},j,k} \overline{\lambda_{p_{l-1},j,k}} \overline{\lambda_{p_{l},j,r_l}},
\end{eqnarray*}
so we have, for $l=1,2,3$,
\begin{equation} \label{eqn:restriction_on_lambdas_for_Weights_E2(12)}
\lambda_{p_{l-1},q_{l-1},r_l} \lambda_{p_{l},q_{l-1},k} \overline{\lambda_{p_{l-1},q_{l-1},k}} \overline{\lambda_{p_{l},q_{l-1},r_l}} = - \lambda_{p_{l-1},j,r_l} \lambda_{p_{l},j,k} \overline{\lambda_{p_{l-1},j,k}} \overline{\lambda_{p_{l},j,r_l}}.
\end{equation}

Then there are two solutions $W^+$, $W^-$ for the cell system for $\mathcal{E}_2^{(12)}$.

\begin{Thm}
Every solution for the cells of $\mathcal{E}_2^{(12)}$ is either equivalent to the solution $W^+$ or the inequivalent conjugate solution $W^-$, given by
$$W_{i,j,k}^{\pm} = \sqrt{[2][3]}, \qquad W_{p_l,j,k}^{\pm} = \frac{[2]\sqrt{[3]}}{\sqrt{[4]}},$$
$$W_{p_l,q_{l-1},r_l}^{\pm} = \frac{\sqrt{[2]}^3}{[4]}\sqrt{[2]^2 \pm \sqrt{[2][4]}}, \quad \;\; W_{p_l,q_{l},r_{l+1}}^{\pm} = - \frac{\sqrt{[2]}^3}{[4]}\sqrt{[2]^2 \mp \sqrt{[2][4]}},$$
$$W_{p_l,q_l,k}^{\pm} = W_{p_l,j,r_{l+1}}^{\pm} = \frac{\sqrt{[2]}^3}{[4]}\sqrt{[2][4] \pm \sqrt{[2][4]}},$$
$$W_{p_l,q_{l-1},k}^{\pm} = W_{p_l,j,r_{l}}^{\pm} = \frac{\sqrt{[2]}^3}{[4]}\sqrt{[2][4] \mp \sqrt{[2][4]}},$$
for $l=1,2,3$.
\end{Thm}
\emph{Proof:}
Let $W^{\sharp}$ be another solution for the cells of $\mathcal{E}_2^{(12)}$, which must be given by $W^{\sharp}_{v_1,v_2,v_3} = \lambda^{\sharp}_{v_1,v_2,v_3} |W^{+}_{v_1,v_2,v_3}|$ where the $\lambda^{\sharp}$'s satisfy the condition (\ref{eqn:restriction_on_lambdas_for_Weights_E2(12)}). We need to find unitaries $u_{v_1,v_2} \in \mathbb{T}$, for $v_1$, $v_2$ vertices of $\mathcal{E}_2^{(12)}$, such that $u_{p_l,q_l} u_{q_l,r_{l+1}} u_{r_{l+1},p_l} \lambda^{\sharp}_{p_l,q_l,r_{l+1}} = -1$, $l=1,2,3$, and $u_{v_1,v_2} u_{v_2,v_3} u_{v_3,v_1} \lambda^{\sharp}_{v_1,v_2,v_3} = 1$ for all other triangles $\triangle_{v_1,v_2,v_3}$ on $\mathcal{E}_2^{(12)}$.
We choose $u_{j,k} = u_{k,i} = u_{j,r_l} = u_{q_l,k} = u_{r_{l+1},p_l} = 1$,
$u_{i,j} = \overline{\lambda^{\sharp}_{i,j,k}}$,
$u_{p_l,j} = \overline{\lambda^{\sharp}_{p_l,j,r_{l+1}}}$,
$u_{k,p_l} = \lambda^{\sharp}_{p_l,j,r_{l+1}} \overline{\lambda^{\sharp}_{p_l,j,k}}$,
$u_{r_l,p_l} = \lambda^{\sharp}_{p_l,j,r_{l+1}} \overline{\lambda^{\sharp}_{p_l,j,r_l}}$,
$u_{p_l,q_l} = \lambda^{\sharp}_{p_l,j,k} \overline{\lambda^{\sharp}_{p_l,q_l,k}} \overline{\lambda^{\sharp}_{p_l,j,r_{l+1}}}$,
$u_{p_l,q_{l-1}} = \lambda^{\sharp}_{p_l,j,k} \overline{\lambda^{\sharp}_{p_l,q_{l-1},k}} \overline{\lambda^{\sharp}_{p_l,j,r_{l+1}}}$,
and finally $u_{q_l,r_{l+1}} = - \lambda^{\sharp}_{p_l,j,r_{l+1}} \lambda^{\sharp}_{p_l,q_l,k} \overline{\lambda^{\sharp}_{p_l,j,k}} \overline{\lambda^{\sharp}_{p_l,j,r_{l+1}}}$, for $l=1,2,3$.

Similarly, for any solution $W^{\sharp \sharp}$ with $|W^{\sharp \sharp}_{v_1,v_2,v_3}| = |W^{-}_{v_1,v_2,v_3}|$.

The solutions $W^+$ and $W^-$ are not equivalent since $|W^+| \neq |W^-|$, and there are no double edges on $\mathcal{E}_2^{(12)}$. We remark that the complex conjugate solutions $\overline{W^{\pm}}$ are equivalent to the solutions $W^{\mp}$: we choose a family of unitaries which satisfy (\ref{eqn:def-equvialence_of_W1,W2}) by $u_{i_l,j_l} = u_{j_l,k_l} = u_{k_l,i_l} = u_{p,j_l} = u_{j_l,r} = u_{q,k_l} = u_{k_l,p} = 1$, $u_{q,r} = -1$, and $2 \times 2$ unitary matrices $u_{\alpha} = u_{\beta} = u$ where $u$ is given by $u(i,j) = 1-\delta_{i,j}$.
\hfill
$\Box$

For $\mathcal{E}_2^{(12)}$, the cells $W^+$ above give the following representation of the Hecke algebra, where $l=1,2,3$ (mod 3):
\begin{eqnarray*}
U^{(i,k)} & = & U^{(j,i)} \;\; = \;\; [2],
\end{eqnarray*}
\begin{eqnarray*}
U^{(k,j)} & = & \begin{array}{c} \scriptstyle i \\ \scriptstyle p_l \end{array} \left( {\begin{array}{cc}
                 \frac{[2]}{[3]} & \frac{\sqrt{[2]^3}}{[3]\sqrt{[4]}} \\
                \frac{\sqrt{[2]^3}}{[3]\sqrt{[4]}} & \frac{[2]^2}{[3][4]}
               \end{array} } \right),
\end{eqnarray*}
\begin{eqnarray*}
U^{(r_l,j)} & = & \begin{array}{c} \scriptstyle p_{l-1} \\ \scriptstyle p_l \end{array} \left( {\begin{array}{cc}
                 \frac{[2]^2 ([2][4]+\sqrt{[2][4]})}{[3]^2[4]} & \frac{\sqrt{[2]^3}}{\sqrt{[3][4]}}  \\
                 \frac{\sqrt{[2]^3}}{\sqrt{[3][4]}}  & \frac{[2]^2 ([2][4]-\sqrt{[2][4]})}{[3]^2[4]}
               \end{array} } \right), \\
\lefteqn{ = \;\; U^{(k,q_l)} \quad \textrm{ with rows labelled by } p_, p_{l+1}, }
\end{eqnarray*}
\begin{eqnarray*}
U^{(q_l,p_l)} & = & \begin{array}{c} \scriptstyle k \\ \scriptstyle r_{l+1} \end{array} \left( {\begin{array}{cc}
                 \frac{[2][4]+\sqrt{[2][4]}}{[2][3]} & \frac{-\sqrt{[2][4]-\sqrt{[2][4]}}}{[2]\sqrt{[3]}} \\
                 \frac{-\sqrt{[2][4]-\sqrt{[2][4]}}}{[2]\sqrt{[3]}} & \frac{[2]^2-\sqrt{[2][4]}}{[2][3]}
               \end{array} } \right), \\
\lefteqn{ = \;\; U^{(p_l,r_{l+1})} \quad \textrm{ with rows labelled by } j, q_l, }
\end{eqnarray*}
\begin{eqnarray*}
U^{(p_l,r_l)} & = & \begin{array}{c} \scriptstyle j \\ \scriptstyle q_{l-1} \end{array} \left( {\begin{array}{cc}
                 \frac{[2][4]-\sqrt{[2][4]}}{[2][3]} & \frac{\sqrt{[2][4]-\sqrt{[2][4]}}}{[2]\sqrt{[3]}} \\
                 \frac{\sqrt{[2][4]-\sqrt{[2][4]}}}{[2]\sqrt{[3]}} & \frac{[2]^2+\sqrt{[2][4]}}{[2][3]}
               \end{array} } \right), \\
\lefteqn{ = \;\; U^{(q_{l-1},p_l)} \quad \textrm{ with rows labelled by } k, r_l, }
\end{eqnarray*}
\begin{eqnarray*}
U^{(r_{l+1},q_l)} & = & \begin{array}{c} \scriptstyle p_l \\ \scriptstyle p_{l+1} \end{array} \left( {\begin{array}{cc}
                 \frac{[2]([2]^2-\sqrt{[2][4]})}{[3]^2} & \frac{-[2]}{\sqrt{[6]}} \\
                 \frac{-[2]}{\sqrt{[6]}} & \frac{[2]([2]^2+\sqrt{[2][4]})}{[3]^2}
               \end{array} } \right),
\end{eqnarray*}
\begin{eqnarray*}
U^{(p_l,k)} & = & \begin{array}{c} \scriptstyle j \\ \scriptstyle q_{l-1} \\ \scriptstyle q_l \end{array} \left( {\begin{array}{ccc}
                 \frac{1}{[2]} & \frac{\sqrt{[2][4] - \sqrt{[2][4]}}}{\sqrt{[2][3][4]}} & \frac{\sqrt{[2][4] + \sqrt{[2][4]}}}{\sqrt{[2][3][4]}} \\
                 \frac{\sqrt{[2][4] - \sqrt{[2][4]}}}{\sqrt{[2][3][4]}} & \frac{[2][4] - \sqrt{[2][4]}}{[3][4]} & \frac{\sqrt{[6]}}{\sqrt{[3][4]}} \\
                 \frac{\sqrt{[2][4] + \sqrt{[2][4]}}}{\sqrt{[2][3][4]}} & \frac{\sqrt{[6]}}{\sqrt{[3][4]}} & \frac{[2][4] + \sqrt{[2][4]}}{[3][4]}
               \end{array} } \right).
\end{eqnarray*}

\section{$\mathcal{E}_1^{(12)}$}

For the graph $\mathcal{E}_1^{(12)}$ (illustrated in Figure \ref{fig:labelled_E1(12)&E2(12)graphs}), we will use the notation $W^{(1)}_{v_1,v_2,v_3}$ for the cell of the triangle $\triangle_{v_1,v_2,v_3}$ where there are no double edges between any of the vertices $v_1$, $v_2$, $v_3$. For triangles that involve the double edges $\alpha, \alpha'$ or $\beta, \beta'$ we will specify which of the double edges is used by the notation $\triangle_{v_1,v_2,v_3}^{(\xi)}$, and $W_{v_1,v_2,v_3(\xi)} := W(\triangle_{v_1,v_2,v_3}^{\xi})$. Since the graph $\mathcal{E}_1^{(12)}$ is a $\mathbb{Z}_3$-orbifold of the graph $\mathcal{E}_2^{(12)}$, we can obtain an orbifold solution for the cells for $\mathcal{E}_1^{(12)}$as follows. We take the $\mathbb{Z}_3$-orbifold of $\mathcal{E}_2^{(12)}$ with the vertices $i$, $j$ and $k$ all fixed points- these are thus triplicated and become the vertices $i_l$, $j_l$ and $k_l$, $l=1,2,3$, on $\mathcal{E}_1^{(12)}$. The vertices $p_1$, $p_2$ and $p_3$ on $\mathcal{E}_2^{(12)}$ are identified and become the vertex $p$ on $\mathcal{E}_1^{(12)}$, and similarly the $q_l$ and $r_l$ become $q$ and $r$. The edges $\alpha_1$, $\alpha_2$ and $\alpha_3$ are identified and become the edge $\alpha$ on $\mathcal{E}_1^{(12)}$, also the edges $\alpha_1'$, $\alpha_2'$ and $\alpha_3'$ are identified and become the edge $\alpha'$. Similarly the edges $\beta_l$, $\beta_l'$ and $\gamma_l$ become the edges $\beta$, $\beta'$ and $\gamma$ respectively on $\mathcal{E}_1^{(12)}$. The Perron-Frobenius weights for the vertices are $\phi_{i_l} = 1$, $\phi_{j_l} = \phi_{k_l} = [3]$, $l=1,2,3$, $\phi_p = [2][4]$ and $\phi_q = \phi_r = [3][4]/[2]$. Note that these are equal to the Perron-Frobenius weights for the corresponding vertices of $\mathcal{E}_2^{(12)}$ up to a scalar factor of $[4]/[2]$.

From the type I frames $\stackrel{i_l}{\bullet} \rightarrow \stackrel{j_l}{\bullet}$, $l=1,2,3$, we have $|W^{(1)}_{i_l,j_l,k_l}|^2 = [2][3]$ (which is equal to $([4]/[2])^2 |W^{(2)}_{i,j,k}|^2 /3$). Then the type I frame $\stackrel{j_l}{\bullet} \rightarrow \stackrel{k_l}{\bullet}$, $l=1,2,3$, gives $|W^{(1)}_{p,j_l,k_l}|^2 = [3][4]$ ($= ([4]/[2])^2 |W^{(2)}_{p_l,j,k}|^2 /3$).
Since the triangle $\triangle_{p,j_l,r}^{(\alpha)}$ in $\mathcal{E}_1^{(12)}$ comes from the triangle $\triangle_{p_l,j,r_l}$ in $\mathcal{E}_2^{(12)}$, then
$$|W^{(1)}_{p,j_l,r(\alpha)}|^2 = \frac{[4]^2}{[2]^2} |W^{(2)}_{p_l,j,r_l}|^2 = [2]([2][4] \mp \sqrt{[2][4]}).$$
The triangle $\triangle_{p,j_l,r}^{(\alpha')}$ in $\mathcal{E}_1^{(12)}$ comes from the triangle $\triangle_{p_l,j,r_{l+1}}$ in $\mathcal{E}_2^{(12)}$, giving
$$|W^{(1)}_{p,j_l,r(\alpha')}|^2 = \frac{[4]^2}{[2]^2} |W^{(2)}_{p_l,j,r_{l+1}}|^2 = [2]([2][4] \pm \sqrt{[2][4]}).$$
Similarly
\begin{eqnarray*}
|W^{(1)}_{p,q,k_l(\beta)}|^2 = \frac{[4]^2}{[2]^2} |W^{(2)}_{p_l,q_l,k}|^2 = [2]([2][4] \pm \sqrt{[2][4]}), \\
|W^{(1)}_{p,q,k_l(\beta')}|^2 = \frac{[4]^2}{[2]^2} |W^{(2)}_{p_l,q_{l-1},k}|^2 = [2]([2][4] \mp \sqrt{[2][4]}).
\end{eqnarray*}
The three triangles $\triangle_{p_l,q_l,r_{l+1}}$, $l=1,2,3$, in $\mathcal{E}_2^{(12)}$ are identified in $\mathcal{E}_1^{(12)}$ and give the triangle $\triangle_{p,q,r}^{(\alpha', \beta)}$, so that $|W^{(1)}_{p,q,r(\alpha', \beta)}|^2 = 3 ([4]/[2])^2 |W^{(2)}_{p_l,q_l,r_{l+1}}|^2 = ([4]/[2])^2([2]^2 \mp \sqrt{[2][4]})$. Similarly $|W^{(1)}_{p,q,r(\alpha, \beta')}|^2 = 3 ([4]/[2])^2 |W^{(2)}_{p_l,q_{l-1},r_l}|^2 = ([4]/[2])^2([2]^2 \pm \sqrt{[2][4]})$.
Considering the type I frame $\stackrel{q}{\bullet} \rightarrow \stackrel{r}{\bullet}$ gives the equation $|W^{(1)}_{p,q,r(\alpha, \beta)}|^2 + |W^{(1)}_{p,q,r(\alpha, \beta')}|^2 + |W^{(1)}_{p,q,r(\alpha', \beta)}|^2 + |W^{(1)}_{p,q,r(\alpha', \beta')}|^2 = [3]^2[4]^2/[2]$. Substituting in for $|W^{(1)}_{p,q,r(\alpha', \beta)}|^2$ and $|W^{(1)}_{p,q,r(\alpha, \beta')}|^2$ we find $|W^{(1)}_{p,q,r(\alpha, \beta)}|^2 + |W^{(1)}_{p,q,r(\alpha', \beta')}|^2 = 0$, so that $|W^{(1)}_{p,q,r(\alpha, \beta)}|^2 = |W^{(1)}_{p,q,r(\alpha', \beta')}|^2 = 0$. The reason for this is that the triangle $\triangle_{p,q,r}^{(\alpha, \beta)}$ (and similarly for the triangle $\triangle_{p,q,r}^{(\alpha', \beta')}$) in $\mathcal{E}_1^{(12)}$ comes from the paths $p_l \rightarrow q_l \rightarrow r_{l+1} \rightarrow p_{l+1}$ in $\mathcal{E}_2^{(12)}$, which do not form a closed triangle.

From the type I frames $\stackrel{r}{\bullet} \rightrightarrows \stackrel{p}{\bullet}$ and $\stackrel{p}{\bullet} \rightrightarrows \stackrel{q}{\bullet}$, we obtain the equations
\begin{eqnarray}
\lambda_{1(\alpha)} \overline{\lambda_{1(\alpha')}} + \lambda_{2(\alpha)} \overline{\lambda_{2(\alpha')}} + \lambda_{3(\alpha)} \overline{\lambda_{3(\alpha')}} & = & 0, \label{eqn:restriction1a_on_lambdas-E1(12)}\\
\lambda_{1(\beta)} \overline{\lambda_{1(\beta')}} + \lambda_{2(\beta)} \overline{\lambda_{2(\beta')}} + \lambda_{3(\beta)} \overline{\lambda_{3(\beta')}} & = & 0, \label{eqn:restriction1b_on_lambdas-E1(12)}
\end{eqnarray}
where $W_{p,j_l,r(\xi)} = \lambda_{l(\xi)} |W_{p,j_l,r(\xi)}|$, for $\xi \in \{ \alpha, \alpha', \beta, \beta' \}$, $l=1,2,3$. Another restriction on the choice of phase is found from the type II frames $\stackrel{j_l}{\bullet} \rightarrow \stackrel{r}{\bullet} \leftarrow \stackrel{j_m}{\bullet}$, for $l \neq m$, $\mathrm{Re}(\lambda_{l(\alpha)} \lambda_{m(\alpha')} \overline{\lambda_{l(\alpha')}} \overline{\lambda_{m(\alpha)}}) = - 1/2$, and similarly for the type II frames $\stackrel{k_l}{\bullet} \rightarrow \stackrel{p}{\bullet} \leftarrow \stackrel{k_m}{\bullet}$, $l \neq m$, giving
\begin{eqnarray}
\lambda_{l(\alpha)} \lambda_{m(\alpha')} \overline{\lambda_{l(\alpha')}} \overline{\lambda_{m(\alpha)}} & = & - \frac{1}{2} + \varepsilon_{l,m} \frac{\sqrt{3}}{2}i, \label{eqn:restriction2a_on_lambdas-E1(12)}\\
\lambda_{l(\beta)} \lambda_{m(\beta')} \overline{\lambda_{l(\beta')}} \overline{\lambda_{m(\beta)}} & = & - \frac{1}{2} + \varepsilon'_{l,m} \frac{\sqrt{3}}{2}i, \label{eqn:restriction2b_on_lambdas-E1(12)}
\end{eqnarray}
where $\varepsilon_{l,m}, \varepsilon'_{l,m} \in \{ \pm 1 \}$. Lastly, from the type II frame consisting of the vertices $j_l$, $k_l$, $q$ and $r$ ($l=1,2,3$) we have
\begin{equation}\label{eqn:restriction3_on_lambdas-E1(12)}
\lambda_{l(\alpha)} \lambda_{l(\beta')} \overline{\lambda_{l(\alpha')}} \overline{\lambda_{l(\beta)}} = - \lambda_{(\alpha \beta')} \overline{\lambda_{(\alpha' \beta)}},
\end{equation}
where $W_{p,q,r(\xi_1, \xi_2)} = \lambda_{(\xi_1, \xi_2)} |W_{p,q,r(\xi_1, \xi_2)}|$, for $\xi_1 \in \{ \alpha, \alpha' \}$, $\xi_2 \in \{ \beta, \beta' \}$, $l=1,2,3$. Then for $l \neq m$,
$$\lambda_{l(\alpha)} \lambda_{m(\alpha')} \overline{\lambda_{l(\alpha')}} \overline{\lambda_{m(\alpha)}} = \lambda_{l(\beta)} \lambda_{m(\beta')} \overline{\lambda_{l(\beta')}} \overline{\lambda_{m(\beta)}},$$
and, from (\ref{eqn:restriction2a_on_lambdas-E1(12)}) and (\ref{eqn:restriction2b_on_lambdas-E1(12)}) we find $\varepsilon_{l,m} = \varepsilon'_{l,m}$. Substituting in for $\lambda_{l(\alpha)} \overline{\lambda_{l(\alpha')}}$ from (\ref{eqn:restriction2a_on_lambdas-E1(12)}) into (\ref{eqn:restriction1a_on_lambdas-E1(12)}), we see that $\varepsilon_{l,l+1} = \varepsilon_{m,m+1}$ for all $l,m = 1,2,3$, and that $\varepsilon_{l,l-1} = - \varepsilon_{l,l+1}$. Then the restrictions for the choice of phase are (\ref{eqn:restriction3_on_lambdas-E1(12)}) and
\begin{equation}\label{eqn:restriction_on_lambdas-E1(12)}
\lambda_{l(\alpha)} \lambda_{l+1(\alpha')} \overline{\lambda_{l(\alpha')}} \overline{\lambda_{l+1(\alpha)}} = \lambda_{l(\beta)} \lambda_{l+1(\beta')} \overline{\lambda_{l(\beta')}} \overline{\lambda_{l+1(\beta)}} = - \frac{1}{2} + \varepsilon \frac{\sqrt{3}}{2}i = e^{\varepsilon \frac{2 \pi i}{3}},
\end{equation}
where $\varepsilon \in \{ \pm 1 \}$.

Then we have obtained two orbifold solutions for the cell system for $\mathcal{E}_1^{(12)}$: $W^+$, $W^-$.

\begin{Thm}
The following solutions $W^+$, $W^-$ for the cells of $\mathcal{E}_1^{(12)}$ are inequivalent:
$$W_{i_l,j_l,k_l}^{\pm} = \sqrt{[2][3]}, \quad \;\; W_{p,j_l,k_l}^{\pm} = \sqrt{[3][4]},$$
$$W_{p,j_l,r(\alpha)}^{\pm} = \epsilon_l \sqrt{[2]}\sqrt{[2][4] \pm \sqrt{[2][4]}}, \quad \;\; W_{p,j_l,r(\alpha')}^{\pm} = \overline{\epsilon_l} \sqrt{[2]}\sqrt{[2][4] \mp \sqrt{[2][4]}},$$
$$W_{p,q,k_l(\beta)}^{\pm} = \epsilon_l \sqrt{[2]}\sqrt{[2][4] \mp \sqrt{[2][4]}}, \quad \;\; W_{p,q,k_l(\beta')}^{\pm} = \overline{\epsilon_l} \sqrt{[2]}\sqrt{[2][4] \pm \sqrt{[2][4]}},$$
$$W_{p,q,r(\alpha \beta')}^{\pm} = \frac{[4]}{\sqrt{[2]}} \sqrt{[2]^2 \mp \sqrt{[2][4]}}, \quad \;\; W_{p,q,r(\alpha' \beta)}^{\pm} = - \frac{[4]}{\sqrt{[2]}} \sqrt{[2]^2 \pm \sqrt{[2][4]}},$$
$$W_{p,q,r(\alpha \beta)}^{\pm} = W_{p,q,r(\alpha' \beta')}^{\pm} = 0,$$
for $l=1,2,3$, where $\epsilon_1 = 1$ and $\epsilon_2 = e^{2 \pi i/3} = \overline{\epsilon_3}$.
\end{Thm}
\emph{Proof:}
The solutions $W^+$, $W^-$ are not equivalent, as can be seen by considering (\ref{eqn:def-equvialence_of_W1,W2}) for the triangle $\triangle_{p,j_l,r}$. We have the following two equations, for $l=1,2,3$:
\begin{eqnarray*}
W^+_{p,j_l,r(\alpha)} & = & u_{p,j_l} u_{j_l,r} \left( u_{\alpha}(\alpha,\alpha) W^-_{p,j_l,r(\alpha)} + u_{\alpha}(\alpha,\alpha') W^-_{p,j_l,r(\alpha')} \right), \\
W^+_{p,j_l,r(\alpha')} & = & u_{p,j_l} u_{j_l,r} \left( u_{\alpha}(\alpha',\alpha) W^-_{p,j_l,r(\alpha)} + u_{\alpha}(\alpha',\alpha') W^-_{p,j_l,r(\alpha')} \right).
\end{eqnarray*}
So we require $u_{p,j_l}, u_{j_l,r} \in \mathbb{T}$ and a $2\times2$ unitary matrix $u_{\alpha}$ such that, for $l=1,2,3$,
\begin{eqnarray}
\epsilon_l \sqrt{[2]} \, x_+ & = & u_{p,j_l} u_{j_l,r} \left( u_{\alpha}(\alpha,\alpha) \epsilon_l \sqrt{[2]} \, x_- + u_{\alpha}(\alpha,\alpha') \overline{\epsilon}_l \sqrt{[2]} \, x_+ \right), \label{eqn:inequivalent_weights_for_E1(12)-1}\\
\overline{\epsilon}_l \sqrt{[2]} \, x_- & = & u_{p,j_l} u_{j_l,r} \left( u_{\alpha}(\alpha',\alpha) \epsilon_l \sqrt{[2]} \, x_- + u_{\alpha}(\alpha',\alpha') \overline{\epsilon}_l \sqrt{[2]} \, x_+ \right). \label{eqn:inequivalent_weights_for_E1(12)-2}
\end{eqnarray}
where $x_{\pm} = \sqrt{[2][4] \pm \sqrt{[2][4]}}$.
Equation (\ref{eqn:inequivalent_weights_for_E1(12)-1}) must hold for each $l=1,2,3$. On the left hand side we have $\epsilon_l$, hence we require $u_{\alpha}(\alpha,\alpha') = 0$ because $u_{\alpha}$ does not depend on $l$, and the difference in phase between $\epsilon_l$ and $\overline{\epsilon}_l$ is 0, $e^{-2 \pi i/3}$, $e^{2 \pi i/3}$ respectively for $l=1,2,3$ respectively. This difference in phase for each $l$ cannot come from $u_{p,j_l} u_{j_l,r}$ (although $u_{p,j_l}$, $u_{j_l,r}$ do depend on $l$) since in (\ref{eqn:inequivalent_weights_for_E1(12)-2}) the difference in phase is now 0, $e^{2 \pi i/3}$, $e^{-2 \pi i/3}$ respectively for $l=1,2,3$ respectively, so we would need $\overline{u_{p,j_l} u_{j_l,r}}$ to take care of the phase difference here, not $u_{p,j_l} u_{j_l,r}$. Then we have $u_{\alpha}(\alpha,\alpha) = \overline{u_{p,j_l} u_{j_l,r}} \; x_+ / x_-$, and similarly $u_{\alpha}(\alpha',\alpha) = 0$ and $u_{\alpha}(\alpha',\alpha') = \overline{u_{p,j_l} u_{j_l,r}} \; x_- / x_+$. But now $u_{\alpha}$ is not unitary.
\hfill
$\Box$

For $\mathcal{E}_1^{(12)}$, the cells $W^+$ above give the following representation of the Hecke algebra, where $l=1,2,3$ (mod 3):
\begin{eqnarray*}
U^{(i_l,k_l)} & = &  U^{(j_l,i_l)} \;\; = \;\; [2],
\qquad \qquad
U^{(k_l,j_l)} \;\; = \;\; \begin{array}{c} \scriptstyle i_l \\ \scriptstyle p \end{array} \left( {\begin{array}{cc}
                 \frac{[2]}{[3]} & \frac{\sqrt{[2][4]}}{[3]} \\
                 \frac{\sqrt{[2][4]}}{[3]} & \frac{[4]}{[3]}
               \end{array} } \right),
\end{eqnarray*}
\begin{eqnarray*}
U^{(r,j_l)} & = & \begin{array}{c} \scriptstyle p(\alpha) \\ \scriptstyle p(\alpha') \end{array} \left( {\begin{array}{cc}
                 \frac{[2]^2([2][4] + \sqrt{[2][4]})}{[3]^2[4]} & \frac{\overline{\epsilon}_l \sqrt{[2]^3}}{\sqrt{[3][4]}} \\
                 \frac{\epsilon_l \sqrt{[2]^3}}{\sqrt{[3][4]}} & \frac{[2]^2([2][4] - \sqrt{[2][4]})}{[3]^2[4]}
               \end{array} } \right), \\
\lefteqn{= \;\; U^{(k_l,q)} \quad \textrm{ with rows labelled by } p(\beta'), p(\beta), }
\end{eqnarray*}
\begin{eqnarray*}
U^{(j_l,p)} & = & \begin{array}{c} \scriptstyle k_l \\ \scriptstyle r(\alpha) \\ \scriptstyle r(\alpha') \end{array} \left( {\begin{array}{ccc}
                 \frac{1}{[2]} & \frac{\overline{\epsilon}_l \sqrt{[2][4] + \sqrt{[2][4]}}}{\sqrt{[2][3][4]}} & \frac{\epsilon_l \sqrt{[2][4] - \sqrt{[2][4]}}}{\sqrt{[2][3][4]}} \\
                 \frac{\epsilon_l \sqrt{[2][4] + \sqrt{[2][4]}}}{\sqrt{[2][3][4]}} & \frac{[2][4] + \sqrt{[2][4]}}{[3][4]} & \frac{\overline{\epsilon}_l \sqrt{[6]}}{\sqrt{[3][4]}} \\
                 \frac{\overline{\epsilon}_l \sqrt{[2][4] - \sqrt{[2][4]}}}{\sqrt{[2][3][4]}} & \frac{\epsilon_l \sqrt{[6]}}{\sqrt{[3][4]}} & \frac{[2][4] - \sqrt{[2][4]}}{[3][4]}
               \end{array} } \right), \\
\lefteqn{ = \;\; U^{(p,k_l)} \quad \textrm{ with rows labelled by } j_l, q(\beta'), q(\beta), }
\end{eqnarray*}
\begin{eqnarray*}
U^{(r,q)} & = & \begin{array}{c} \scriptstyle p(\alpha \beta) \\ \scriptstyle p(\alpha \beta') \\ \scriptstyle p(\alpha' \beta) \\ \scriptstyle p(\alpha' \beta') \end{array} \left( {\begin{array}{cccc}
                 0 & 0 & 0 & 0 \\
                 0 & \frac{[2]([2]^2 - \sqrt{[2][4]})}{[3]^2} & -\frac{\sqrt{[2]}}{\sqrt{[6]}} & 0 \\
                 0 & -\frac{\sqrt{[2]}}{\sqrt{[6]}} & \frac{[2]([2]^2 + \sqrt{[2][4]})}{[3]^2} & 0 \\
                 0 & 0 & 0 & 0
               \end{array} } \right),
\end{eqnarray*}
\begin{eqnarray*}
\lefteqn{ U^{(p,r)} \quad \textrm{ with labels } j_1, j_2, j_3, q(\beta), q(\beta') } \\
& = & \left( {\begin{array}{ccccc}
                 \frac{[3]}{[4]} & a_- & a_+ & - \sqrt{b_+} & \sqrt{b_-} \\
                 a_+ &  \frac{[3]}{[4]} & a_- & -\overline{\epsilon_2} \sqrt{b_+} & \epsilon_2 \sqrt{b_-} \\
                 a_- & a_+ & \frac{[3]}{[4]} & -\epsilon_2 \sqrt{b_+} & \overline{\epsilon}_2 \sqrt{b_-} \\
                 - \sqrt{b_+} & -\epsilon_2 \sqrt{b_+} & -\overline{\epsilon}_2 \sqrt{b_+} & \frac{[2]^2 + \sqrt{[2][4]}}{[2][3]} & 0 \\
                 \sqrt{b_-} & \overline{\epsilon}_2 \sqrt{b_-} & \epsilon_2\ \sqrt{b_-} & 0 & \frac{[2]^2 - \sqrt{[2][4]}}{[2][3]}
               \end{array} } \right) \\
& = & U^{(q,p)} \quad \textrm{ with labels } k_1, k_3, k_2, r(\alpha'), r(\alpha),
\end{eqnarray*}
where $a_{\pm} = (-[2]^2 \pm i\sqrt{[2][4]} \, )/[3][4]$, $b_{\pm} = ([2][4] \pm \sqrt{[2][4]} \, )/[3][4]^2$.

Our representation of the Hecke algebra is not equivalent to that given by Sochen for $\mathcal{E}_1^{(12)}$ in \cite{sochen:1991}, however we believe that there is a typographical error in Sochen's presentation and that the weights he denotes by $U^{(4,2_r)} = (U^{(3_r,6)})^{\ast}$ should be the complex conjugate of the one given. In this case, the representation of the Hecke algebra we give above can be shown to be equivalent by choosing a family of unitaries $u_{i_l,j_l} = u_{j_l,k_l} = u_{k_l,i_l} = u_{p,j_l} = u_{k_l,p} = u_{q,r} = 1$, $u_{j_l,r} = - \epsilon_l = \overline{u_{q,k_l}}$ and set the $2 \times 2$ unitary matrices $u_{\alpha}$, $u_{\beta}$ to be the identity matrix.

\section{$\mathcal{E}_5^{(12)}$}

\begin{figure}[tb]
\begin{center}
\includegraphics[width=70mm]{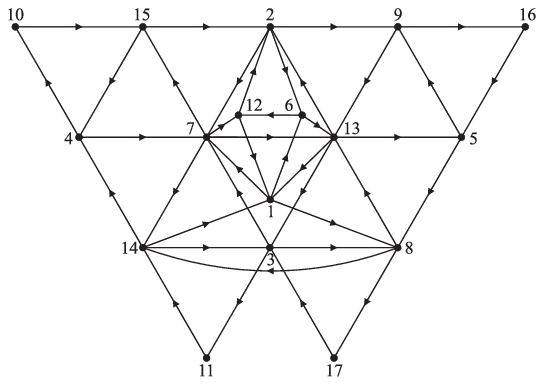}\\
 \caption{Labelled graph $\mathcal{E}_5^{(12)}$} \label{fig:Weights7}
\end{center}
\end{figure}

We label the vertices of $\mathcal{E}_5^{(12)}$ as in Figure \ref{fig:Weights7}. The Perron-Frobenius weights associated to the vertices are $\phi_1 = [3][6]/[2]$, $\phi_2 = \phi_3 = \phi_8 = \phi_{14} = [3][4]/[2]$, $\phi_4 = \phi_5 = \phi_9 = \phi_{15} = [3]$, $\phi_6 = \phi_{12} = [2][3]^2/[6] = [2]^2$, $\phi_7 = \phi_{13} = [3]^2[4]/[6] = [2][4]$, $\phi_{10} = \phi_{16} = 1$, $\phi_{11} = \phi_{17} = [4]/[2]$. The distinguished $\ast$-vertex is vertex 10.

With $W_{v_1,v_2,v_3} = \lambda_{v_1,v_2,v_3} |W_{v_1,v_2,v_3}|$, $\lambda_{v_1,v_2,v_3} \in \mathbb{T}$, we find two restrictions on the choice of phase
\begin{eqnarray}
\lambda_{1,6,12} \lambda_{2,7,12}\overline{ \lambda_{1,7,12}} \overline{\lambda_{2,6,12}} = - \lambda_{1,6,13} \lambda_{2,7,13}\overline{ \lambda_{1,7,13}} \overline{\lambda_{2,6,13}}, \label{eqn:restriction_on_lambdas_for_Weights_E5(12)-1}\\
\lambda_{1,7,14} \lambda_{1,8,13}\overline{ \lambda_{1,7,13}} \overline{\lambda_{1,8,14}} = - \lambda_{3,7,14} \lambda_{3,8,13}\overline{ \lambda_{3,7,13}} \overline{\lambda_{3,8,14}} . \label{eqn:restriction_on_lambdas_for_Weights_E5(12)-2}
\end{eqnarray}

\begin{Thm}
There is up to equivalence a unique set of cells for $\mathcal{E}_5^{(12)}$ given by
$$W_{1,6,12} = W_{4,10,15} = W_{5,9,16} = \sqrt{[2][3]},$$
$$W_{1,6,13} = W_{1,7,12} = [2]{\sqrt{[3][4]}},$$
$$W_{1,7,13} = W_{3,7,14} = W_{3,8,13} = W_{3,8,17} = W_{3,11,14} = W_{2,7,15} = W_{2,9,13}$$
$$ = W_{4,7,14} = W_{5,8,13} = \frac{[4]\sqrt{[3]}}{\sqrt{[2]}},$$
$$W_{1,8,14} = \frac{[4]\sqrt{[3][6]}}{\sqrt{[2]}}, \qquad \qquad \qquad W_{1,7,14} = W_{1,8,13} = \frac{\sqrt{[3][4][6]}}{\sqrt{[2]}},$$
$$W_{2,6,12} = [4]\sqrt{[2]}, \qquad \qquad \qquad W_{2,6,13} = W_{2,7,12} = [2]\sqrt{[4]},$$
$$W_{2,7,13} = - [4]\sqrt{[2]}, \qquad \qquad \qquad W_{3,7,13} = - [4]\sqrt{[6]},$$
$$W_{3,8,14} = \frac{[4]\sqrt{[6]}}{[2]}, \qquad \qquad \qquad W_{4,7,15} = W_{5,9,13} = \sqrt{[3][4]}.$$
$$\;$$
\end{Thm}
\emph{Proof:}
Let $W^{\sharp}$ be any other solution for the cells of $\mathcal{E}_5^{(12)}$. Then we have $W^{\sharp}_{v_1,v_2,v_3} = \lambda^{\sharp}_{v_1,v_2,v_3} |W_{v_1,v_2,v_3}|$, where the $\lambda^{\sharp}$'s satisfy the conditions (\ref{eqn:restriction_on_lambdas_for_Weights_E5(12)-1}) and (\ref{eqn:restriction_on_lambdas_for_Weights_E5(12)-2}). We need to find unitaries $u_{v_1,v_2} \in \mathbb{T}$ which satisfy $u_{7,13} u_{13,2} u_{2,7} \lambda^{\sharp}_{2,7,13} = -1$, $u_{7,13} u_{13,3} u_{3,7} \lambda^{\sharp}_{3,7,13} = -1$ and $u_{v_1,v_2} u_{v_2,v_3} u_{v_3,v_1} \lambda^{\sharp}_{v_1,v_2,v_3} = 1$ for all other triangles $\triangle_{v_1,v_2,v_3}$ on $\mathcal{E}_5^{(12)}$.
We choose $u_{2,7} = u_{2,9} = u_{3,8} = u_{3,11} = u_{6,13} = u_{7,13} = u_{7,14} = u_{8,13} = u_{8,17} = u_{9,16} = u_{10,15} = u_{12,1} = u_{12,2} = u_{13,5} = u_{14,7} = u_{15,2} = 1$,
$\; u_{5,8} = \overline{\lambda^{\sharp}_{5,8,13}}$,
$\; u_{7,12} = \overline{\lambda^{\sharp}_{2,7,12}}$,
$\; u_{7,15} = \overline{\lambda^{\sharp}_{2,7,15}}$,
$\; u_{11,14} = - \overline{\lambda^{\sharp}_{3,11,14}}$,
$\; u_{13,1} = \overline{\lambda^{\sharp}_{1,6,13}}$,
$\; u_{13,2} = - \overline{\lambda^{\sharp}_{2,7,13}}$,
$\; u_{13,3} = \overline{\lambda^{\sharp}_{3,8,13}}$,
$\; u_{14,4} = \overline{\lambda^{\sharp}_{4,7,14}}$,
$\; u_{17,3} = \overline{\lambda^{\sharp}_{3,8,17}}$,
$\; u_{1,7} = \lambda^{\sharp}_{2,7,12} \overline{\lambda^{\sharp}_{1,7,12}}$,
$\; u_{2,6} = - \lambda^{\sharp}_{2,7,13} \overline{\lambda^{\sharp}_{2,6,13}}$,
$\; u_{3,7} = - \lambda^{\sharp}_{3,8,13} \overline{\lambda^{\sharp}_{3,7,13}}$,
$\; u_{9,13} = - \lambda^{\sharp}_{2,7,13} \overline{\lambda^{\sharp}_{2,9,13}}$,
$\; u_{15,4} = \lambda^{\sharp}_{2,7,15} \overline{\lambda^{\sharp}_{4,7,15}}$,
$\; u_{4,10} = \lambda^{\sharp}_{4,7,15} \overline{\lambda^{\sharp}_{2,7,15}} \overline{\lambda^{\sharp}_{4,10,15}}$,
$\; u_{5,9} = - \lambda^{\sharp}_{2,9,13} \overline{\lambda^{\sharp}_{2,7,13}} \overline{\lambda^{\sharp}_{5,9,13}}$,
$\; u_{6,12} = - \lambda^{\sharp}_{2,6,13} \overline{\lambda^{\sharp}_{2,6,12}} \overline{\lambda^{\sharp}_{2,7,13}}$, \\
$\; u_{14,1} = \lambda^{\sharp}_{1,7,12} \overline{\lambda^{\sharp}_{1,7,14}} \overline{\lambda^{\sharp}_{2,7,12}}$,
$\;u_{14,3} = - \lambda^{\sharp}_{3,7,13} \overline{\lambda^{\sharp}_{3,7,14}} \overline{\lambda^{\sharp}_{3,8,13}}$, \\
$\; u_{1,6} = - \lambda^{\sharp}_{2,6,12} \lambda^{\sharp}_{2,7,13} \overline{\lambda^{\sharp}_{1,6,12}} \overline{\lambda^{\sharp}_{2,6,13}}$,
$\; u_{1,8} = \lambda^{\sharp}_{1,7,12} \overline{\lambda^{\sharp}_{1,7,13}} \overline{\lambda^{\sharp}_{1,8,13}} \overline{\lambda^{\sharp}_{2,7,12}}$, \\
$\; u_{8,14} = \lambda^{\sharp}_{1,7,14} \lambda^{\sharp}_{1,8,13} \overline{\lambda^{\sharp}_{1,7,13}} \overline{\lambda^{\sharp}_{1,8,14}}$
and $u_{16,5} = - \lambda^{\sharp}_{2,7,13} \lambda^{\sharp}_{5,9,13} \overline{\lambda^{\sharp}_{2,9,13}} \overline{\lambda^{\sharp}_{5,9,16}}$.
\hfill
$\Box$

For $\mathcal{E}_5^{(12)}$, we have the following representation of the Hecke algebra:
\begin{eqnarray*}
U^{(5,16)} & = & U^{(16,9)} \;\; = \;\; U^{(10,4)} \;\; = \;\; U^{(15,10)} \;\; = \;\; [2], \\
U^{(3,17)} & = & U^{(17,8)} \;\; = \;\; U^{(11,3)} \;\; = \;\; U^{(14,11)} \;\; = \;\; \frac{[2]}{[4]}, \\
U^{(2,15)} & = & U^{(4,14)} \;\; = \;\; U^{(8,5)} \;\; = \;\; U^{(9,2)} \;\; = \;\; \frac{[4]}{[3]},
\end{eqnarray*}
\begin{eqnarray*}
U^{(14,8)} & = & \begin{array}{c}  \scriptstyle 3 \\ \scriptstyle 1 \end{array} \left( {\begin{array}{cc}
                \frac{1}{[2]} & \frac{\sqrt{[3]}}{\sqrt{[2]}} \\
                \frac{\sqrt{[3]}}{\sqrt{[2]}} & [3]
               \end{array} } \right),
\end{eqnarray*}
\begin{eqnarray*}
U^{(12,7)} & = & \begin{array}{c} \scriptstyle 2 \\ \scriptstyle 1 \end{array} \left( {\begin{array}{cc}
                \frac{1}{[2]} & \frac{\sqrt{[3]}}{[2]} \\
                \frac{\sqrt{[3]}}{[2]} & \frac{[3]}{[2]}
               \end{array} } \right) \;\; = \;\; U^{(13,6)} \quad \textrm{ with rows labelled by } 2,1,
\end{eqnarray*}
\begin{eqnarray*}
U^{(3,13)} & = & \begin{array}{c} \scriptstyle 8 \\ \scriptstyle 7 \end{array} \left( {\begin{array}{cc}
                \frac{1}{[2]} & -\frac{\sqrt{[3]}}{[2]} \\
                -\frac{\sqrt{[3]}}{[2]} & \frac{[3]}{[2]}
               \end{array} } \right) \; = \;\; U^{(7,3)} \quad \textrm{ with rows labelled by } 14,13,
\end{eqnarray*}
\begin{eqnarray*}
U^{(5,13)} & = & \begin{array}{c} \scriptstyle 9 \\ \scriptstyle 8 \end{array} \left( {\begin{array}{cc}
                \frac{1}{[2]} & \frac{\sqrt{[4]}}{\sqrt{[2]^3}} \\
                \frac{\sqrt{[4]}}{\sqrt{[2]^3}} & \frac{[4]}{[2]^2}
               \end{array} } \right) \;\; = \;\; U^{(13,9)} \textrm{ with labels 5,2} \\
& = & U^{(7,4)} \textrm{ with labels 15,14} \;\; = \;\; U^{(15,7)} \quad \textrm{ with labels 4,2,}
\end{eqnarray*}
\begin{eqnarray*}
U^{(2,12)} & = & \begin{array}{c} \scriptstyle 7 \\ \scriptstyle 6 \end{array} \left( {\begin{array}{cc}
                 \frac{[2]}{[3]} & \frac{\sqrt{[2][4]}}{[3]} \\
                 \frac{\sqrt{[2][4]}}{[3]} & \frac{[4]}{[3]}
               \end{array} } \right) \;\; = \;\; U^{(6,2)} \textrm{ with labels 13,12} \\
& = & U^{(4,15)} \textrm{ with labels 10,7} \;\; = \;\; U^{(9,5)} \quad \textrm{ with labels 16,13,}
\end{eqnarray*}
\begin{eqnarray*}
U^{(1,14)} & = & \begin{array}{c} \scriptstyle 7 \\ \scriptstyle 8 \end{array} \left( {\begin{array}{cc}
                 \frac{[2]}{[3]} & \frac{[2]\sqrt{[4]}}{[3]} \\
                 \frac{[2]\sqrt{[4]}}{[3]} & \frac{[2][4]}{[3]}
               \end{array} } \right) \;\; = \;\; U^{(8,1)} \quad \textrm{ with labels 13,14,}
\end{eqnarray*}
\begin{eqnarray*}
U^{(12,6)} & = & \begin{array}{c} \scriptstyle 1 \\ \scriptstyle 2 \end{array} \left( {\begin{array}{cc}
                 \frac{[3]}{[2]^3} & \frac{[4]\sqrt{[3]}}{[2]^3} \\
                 \frac{[4]\sqrt{[3]}}{[2]^3} & \frac{[4]^2}{[2]^3}
               \end{array} } \right),
\end{eqnarray*}
\begin{eqnarray*}
U^{(1,12)} & = & \begin{array}{c} \scriptstyle 6 \\ \scriptstyle 7 \end{array} \left( {\begin{array}{cc}
                 \frac{1}{[6]} & \frac{\sqrt{[2][4]}}{[6]} \\
                 \frac{\sqrt{[2][4]}}{[6]} & \frac{[2][4]}{[6]}
               \end{array} } \right) \;\; = \;\; U^{(6,1)} \quad \textrm{ with labels 12,13,}
\end{eqnarray*}
\begin{eqnarray*}
U^{(13,8)} & = & \begin{array}{c} \scriptstyle 5 \\ \scriptstyle 3 \\ \scriptstyle 1 \end{array} \left( {\begin{array}{ccc}
                 \frac{1}{[2]} & \frac{1}{[2]} & \frac{\sqrt{[6]}}{[2]\sqrt{[4]}} \\
                 \frac{1}{[2]} & \frac{1}{[2]} & \frac{\sqrt{[6]}}{[2]\sqrt{[4]}} \\
                 \frac{\sqrt{[6]}}{[2]\sqrt{[4]}} & \frac{\sqrt{[6]}}{[2]\sqrt{[4]}} & \frac{[6]}{[2][4]}
               \end{array} } \right) \;\; = \;\; U^{(14,7)} \quad \textrm{ with labels 4,3,1,}
\end{eqnarray*}
\begin{eqnarray*}
U^{(3,14)} & = & \begin{array}{c} \scriptstyle 8 \\ \scriptstyle 7 \\ \scriptstyle 11 \end{array} \left( {\begin{array}{ccc}
                 \frac{1}{[2]} & \frac{1}{\sqrt{[3]}} & \frac{1}{\sqrt{[3]}} \\
                 \frac{1}{\sqrt{[3]}} & \frac{[2]}{[3]} & \frac{[2]}{[3]} \\
                 \frac{1}{\sqrt{[3]}} & \frac{[2]}{[3]} & \frac{[2]}{[3]}
               \end{array} } \right) \;\; = \;\; U^{(8,3)} \quad \textrm{ with labels 14,13,17,}
\end{eqnarray*}
\begin{eqnarray*}
U^{(2,13)} & = & \begin{array}{c} \scriptstyle 9 \\ \scriptstyle 7 \\ \scriptstyle 6 \end{array} \left( {\begin{array}{ccc}
                 \frac{1}{[2]} & -\frac{1}{\sqrt{[3]}} & \frac{\sqrt{[2]}}{\sqrt{[3][4]}} \\
                 -\frac{1}{\sqrt{[3]}} & \frac{[2]}{[3]} & -\frac{\sqrt{[2]^3}}{[3]\sqrt{[4]}} \\
                 \frac{\sqrt{[2]}}{\sqrt{[3][4]}} & -\frac{\sqrt{[2]^3}}{[3]\sqrt{[4]}} & \frac{[2]^2}{[3][4]}
               \end{array} } \right) \\
& = & U^{(7,2)} \quad \textrm{ with labels 15,13,12,}
\end{eqnarray*}
\begin{eqnarray*}
U^{(1,13)} & = & \begin{array}{c} \scriptstyle 8 \\ \scriptstyle 7 \\ \scriptstyle 6 \end{array} \left( {\begin{array}{ccc}
                 \frac{1}{[2]} & \frac{\sqrt{[4]}}{[2]\sqrt{[6]}} & \frac{\sqrt{[2]^3}}{\sqrt{[6]}} \\
                 \frac{\sqrt{[4]}}{[2]\sqrt{[6]}} & \frac{[4]}{[2][6]} & \frac{\sqrt{[2]^3[4]}}{[6]} \\
                 \frac{\sqrt{[2]^3}}{\sqrt{[6]}} & \frac{\sqrt{[2]^3[4]}}{[6]} & \frac{[2]^2}{[6]}
               \end{array} } \right) \\
& = & U^{(7,1)} \quad \textrm{ with labels 14,13,12,}
\end{eqnarray*}
\begin{eqnarray*}
U^{(13,7)} & = & \begin{array}{c} \scriptstyle 2 \\ \scriptstyle 3 \\ \scriptstyle 1 \end{array} \left( {\begin{array}{ccc}
                 \frac{1}{[2]} & \frac{\sqrt{[6]}}{\sqrt{[2]^3}} & -\frac{\sqrt{[3]}}{[2]^2} \\
                 \frac{\sqrt{[6]}}{\sqrt{[2]^3}} & \frac{[6]}{[2]^2} & -\frac{\sqrt{[3][6]}}{\sqrt{[2]^5}} \\
                 -\frac{\sqrt{[3]}}{[2]^2} & -\frac{\sqrt{[3][6]}}{\sqrt{[2]^5}} & \frac{[3]}{[2]^3}
               \end{array} } \right).
\end{eqnarray*}

\section{$\mathcal{E}^{(24)}$}\label{Sect:computation_of_weights-E(24)}

\begin{figure}[tb]
\begin{center}
\includegraphics[width=80mm]{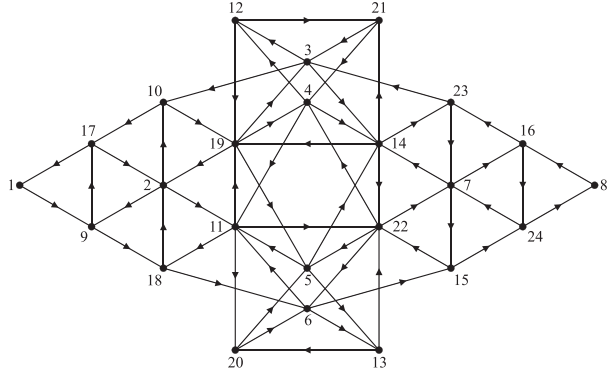}\\
 \caption{Labelled graph $\mathcal{E}^{(24)}$}\label{fig:Weights8}
\end{center}
\end{figure}

We label the vertices of the graph $\mathcal{E}^{(24)}$ as in Figure \ref{fig:Weights8}. The Perron-Frobenius weights are: $\phi_1 = \phi_8 = 1$, $\phi_2 = \phi_7 = [2][4]$, $\phi_3 = \phi_6 = [4][5]/[2]$, $\phi_4 = \phi_5 = [4][7]/[2]$, $\phi_9 = \phi_{16} = \phi_{17} = \phi_{24} = [3]$, $\phi_{10} = \phi_{15} = \phi_{18} = \phi_{23} = [3][4]/[2]$, $\phi_{11} = \phi_{14} = \phi_{19} = \phi_{22} = [3][5]$ and $\phi_{12} = \phi_{13} = \phi_{20} = \phi_{21} = [9]$. With $[a]=[a]_q$, $q=e^{i \pi/24}$, we have the relation $[4]^2 = [2][10]$.

The following cells follow from the $\mathcal{A}$ case: $|W_{1,9,17}|^2 = |W_{8,16,24}|^2 = [2][3]$, $|W_{2,9,17}|^2 = |W_{7,16,24}|^2 = [3][4]$, $|W_{2,9,18}|^2 = |W_{2,10,17}|^2 = |W_{7,15,24}|^2 = |W_{7,16,23}|^2 = [3]^2[4]$, $|W_{2,10,19}|^2 = |W_{2,11,18}|^2 = |W_{7,14,23}|^2 = |W_{7,15,22}|^2 = [3][4][5]$, $|W_{2,11,19}|^2 = |W_{7,14,22}|^2 = [3]^2[4][5]$ and $|W_{3,10,19}|^2 = |W_{3,14,23}|^2 = |W_{6,11,18}|^2 = |W_{6,15,22}|^2 = [3][4]^2[5]/[2]$.

The type II frame $\stackrel{2}{\bullet} \rightarrow \stackrel{19}{\bullet} \leftarrow \stackrel{4}{\bullet}$ gives $\phi_{11}^{-1} |W_{2,11,19}|^2 |W_{4,11,19}|^2 = [3][4]^2[5][7]$, and so we obtain $|W_{4,11,19}|^2 = [4][5][7]$. From the type I frame $\stackrel{11}{\bullet} \rightarrow \stackrel{19}{\bullet}$ we have the equation $|W_{2,11,19}|^2 + |W_{4,11,19}|^2 + |W_{5,11,19}|^2 = [2][3]^2[5]^2$, giving $|W_{5,11,19}|^2 = [4][5][7] = |W_{4,11,19}|^2$. Then by considering the type I frames $\stackrel{4}{\bullet} \rightarrow \stackrel{11}{\bullet}$ and $\stackrel{22}{\bullet} \rightarrow \stackrel{4}{\bullet}$, we find $|W_{4,14,22}|^2 = |W_{5,14,22}|^2 = |W_{4,11,19}|^2 = |W_{5,11,19}|^2$. Similarly $|W_{4,12,19}|^2 = |W_{4,14,21}|^2 = |W_{5,11,20}|^2 = |W_{5,13,22}|^2$ and $|W_{3,12,19}|^2 = |W_{3,14,21}|^2 = |W_{6,11,20}|^2 = |W_{6,13,22}|^2$, and the cells have a $\mathbb{Z}_2$ symmetry.

From type I frames we have the equations:
\begin{eqnarray}
|W_{4,11,19}|^2 + |W_{4,12,19}|^2 + |W_{4,14,19}|^2 & = & [3][4][5][7], \label{eqn:E(24)_Weights-II}\\
|W_{3,12,19}|^2 + |W_{4,12,19}|^2 & = & [2][3][5][9], \label{eqn:E(24)_Weights-III}\\
|W_{3,10,19}|^2 + |W_{3,12,19}|^2 + |W_{3,14,19}|^2 & = & [3][4][5]^2, \label{eqn:E(24)_Weights-IV}\\
|W_{3,14,19}|^2 + |W_{4,14,19}|^2 + |W_{5,14,19}|^2 & = & [2][3]^2[5]^2, \label{eqn:E(24)_Weights-V}\\
|W_{3,12,19}|^2 + |W_{3,12,21}|^2 & = & [4][5][9], \label{eqn:E(24)_Weights-VI}\\
|W_{3,12,21}|^2 + |W_{4,12,21}|^2 & = & [2][9]^2. \label{eqn:E(24)_Weights-VII}
\end{eqnarray}
The type II frame $\stackrel{11}{\bullet} \rightarrow \stackrel{19}{\bullet} \leftarrow \stackrel{12}{\bullet}$, gives $\phi_4^{-1} |W_{4,11,19}|^2 |W_{4,12,19}|^2 = [3]^2[5]^2[9]$, so $|W_{4,12,19}|^2 = [3]^2[5][9]/[2]$. Then using the equations (\ref{eqn:E(24)_Weights-II})-(\ref{eqn:E(24)_Weights-VII}) we obtain $|W_{4,14,19}|^2 = [5]^2[7]/[2]$, $\; |W_{3,12,19}|^2 = [3][5][9]/[2]$, $\; |W_{3,14,19}|^2 = [3]^2[5]^2/[2]$, $\; |W_{5,14,19}|^2 = [5][7][10]$, $\; |W_{3,12,21}|^2 = [5]^2[9]/[2]$ and $|W_{4,12,21}|^2 = [7][9]/[2]$.

With $W_{v_1,v_2,v_3} = \lambda_{v_1,v_2,v_3} |W_{v_1,v_2,v_3}|$, $\lambda_{v_1,v_2,v_3} \in \mathbb{T}$, we have the following restrictions on the $\lambda$'s:
\begin{eqnarray}
\lambda_{3,12,19} \lambda_{3,14,21} \overline{\lambda_{3,12,21}} \overline{\lambda_{3,14,19}} & = & - \lambda_{4,12,19} \lambda_{4,14,21} \overline{\lambda_{4,12,21}} \overline{\lambda_{4,14,19}}, \label{eqn:restriction_on_lambdas_for_Weights_E(24)-1}\\
\lambda_{4,11,22} \lambda_{4,14,19} \overline{\lambda_{4,11,19}} \overline{\lambda_{4,14,22}} & = & - \lambda_{5,11,22} \lambda_{5,14,19} \overline{\lambda_{5,11,19}} \overline{\lambda_{5,14,22}}, \label{eqn:restriction_on_lambdas_for_Weights_E(24)-2}\\
\lambda_{5,11,20} \lambda_{5,13,22} \overline{\lambda_{5,11,22}} \overline{\lambda_{5,13,20}} & = & - \lambda_{6,11,20} \lambda_{6,13,22} \overline{\lambda_{6,11,22}} \overline{\lambda_{6,13,20}}. \label{eqn:restriction_on_lambdas_for_Weights_E(24)-3}
\end{eqnarray}

\begin{Thm}
There is up to equivalence a unique set of cells for $\mathcal{E}^{(24)}$ given by
$$W_{1,9,17} = W_{8,16,24} = \sqrt{[2][3]}, \qquad \qquad \qquad W_{2,9,17} = W_{7,16,24} = \sqrt{[3][4]},$$
$$W_{2,9,18} = W_{2,10,17} = W_{7,15,24} = W_{7,16,23} = [3]\sqrt{[4]},$$
$$W_{2,10,19} = W_{2,11,18} = W_{7,14,23} = W_{7,15,22} = \sqrt{[3][4][5]},$$
$$W_{2,11,19} = W_{7,14,22} = [3]\sqrt{[4][5]},$$
$$W_{3,10,19} = W_{3,14,23} = W_{6,11,18} = W_{6,15,22} = \frac{[4]\sqrt{[3][5]}}{\sqrt{[2]}},$$
$$W_{4,11,19} = W_{4,14,22} = W_{5,11,19} = W_{5,14,22} = \sqrt{[4][5][7]},$$
$$W_{4,12,19} = W_{4,14,21} = W_{5,11,20} = W_{5,13,22} = \frac{[3]\sqrt{[5][9]}}{\sqrt{[2]}},$$
$$W_{3,12,19} = W_{3,14,21} = W_{6,11,20} = W_{6,13,22} = \frac{\sqrt{[3][5][9]}}{\sqrt{[2]}},$$
$$W_{3,14,19} = W_{6,11,22} = \frac{[3][5]}{\sqrt{[2]}}, \qquad \qquad \qquad W_{4,14,19} = W_{5,11,22} = \frac{[5]\sqrt{[7]}}{\sqrt{[2]}},$$
$$W_{5,14,19} = \sqrt{[5][7][10]}, \qquad \qquad \qquad W_{4,11,22} = - \sqrt{[5][7][10]},$$
$$W_{3,12,21} = W_{6,13,20} = - \frac{[5]\sqrt{[9]}}{\sqrt{[2]}}, \qquad \qquad \qquad W_{4,12,21} = W_{5,13,20} = \frac{\sqrt{[7][9]}}{\sqrt{[2]}}.$$
\end{Thm}
\emph{Proof:}
Let $W^{\sharp}$ be any solution for the cells of $\mathcal{E}^{(24)}$. Then $W^{\sharp}_{v_1,v_2,v_3} = \lambda^{\sharp}_{v_1,v_2,v_3} |W_{v_1,v_2,v_3}|$, where the $\lambda^{\sharp}$'s satisfy the conditions (\ref{eqn:restriction_on_lambdas_for_Weights_E(24)-1}), (\ref{eqn:restriction_on_lambdas_for_Weights_E(24)-2}) and (\ref{eqn:restriction_on_lambdas_for_Weights_E(24)-3}). We need to find unitaries $u_{v_1,v_2} \in \mathbb{T}$, for vertices $v_1$, $v_2$ of $\mathcal{E}^{(24)}$, such that $u_{12,21} u_{21,3} u_{3,12} \lambda^{\sharp}_{3,12,21} = -1$, $u_{13,20} u_{20,6} u_{6,13} \lambda^{\sharp}_{6,13,20} = -1$, $u_{11,22} u_{22,4} u_{4,11} \lambda^{\sharp}_{4,11,22} = -1$, and for all other triangles $\triangle_{v_1,v_2,v_3}$ on $\mathcal{E}^{(24)}$ we require $u_{v_1,v_2} u_{v_2,v_3} u_{v_3,v_1} \lambda^{\sharp}_{v_1,v_2,v_3} = 1$.
We make the following choices for the $u_{v_1,v_2}$:
$$u_{3,12} = u_{3,14} = u_{4,11} = u_{5,13} = u_{5,14} = u_{11,20}$$
$$ = u_{14,19} = u_{20,6} = u_{21,3} = u_{21,4} = u_{22,6} = 1,$$
$$u_{12,21} = - \overline{\lambda^{\sharp}_{3,12,21}}, \quad u_{14,21} = \overline{\lambda^{\sharp}_{3,14,21}}, \quad u_{19,3} = \overline{\lambda^{\sharp}_{3,14,19}}, \quad u_{19,5} = - \overline{\lambda^{\sharp}_{5,14,19}},$$
$$u_{4,12} = - \lambda^{\sharp}_{3,12,21} \overline{\lambda^{\sharp}_{4,12,21}}, \quad u_{4,14} = \lambda^{\sharp}_{3,14,21} \overline{\lambda^{\sharp}_{4,14,21}}, \quad u_{6,11} = \lambda^{\sharp}_{5,14,22} \overline{\lambda^{\sharp}_{6,11,20}},$$
$$u_{12,19} = \lambda^{\sharp}_{3,14,19} \overline{\lambda^{\sharp}_{3,12,19}}, \qquad u_{11,22} = \lambda^{\sharp}_{6,11,20} \overline{\lambda^{\sharp}_{5,14,22}} \overline{\lambda^{\sharp}_{6,11,22}},$$
$$u_{19,4} = \lambda^{\sharp}_{4,14,21} \overline{\lambda^{\sharp}_{3,14,21}} \overline{\lambda^{\sharp}_{4,14,19}}, \qquad u_{22,4} = - \lambda^{\sharp}_{5,14,22} \lambda^{\sharp}_{6,11,22} \overline{\lambda^{\sharp}_{4,11,22}} \overline{\lambda^{\sharp}_{6,11,20}},$$
$$u_{5,11} = - \lambda^{\sharp}_{4,11,22} \lambda^{\sharp}_{4,14,21} \lambda^{\sharp}_{5,14,22} \overline{\lambda^{\sharp}_{3,14,21}} \overline{\lambda^{\sharp}_{4,14,22}} \overline{\lambda^{\sharp}_{5,11,22}},$$
$$u_{11,19} = - \lambda^{\sharp}_{3,12,21} \lambda^{\sharp}_{3,14,19} \lambda^{\sharp}_{4,12,19} \overline{\lambda^{\sharp}_{3,12,19}} \overline{\lambda^{\sharp}_{4,11,19}} \overline{\lambda^{\sharp}_{4,12,21}},$$
$$u_{20,5} = - \lambda^{\sharp}_{3,14,21} \lambda^{\sharp}_{4,14,22} \lambda^{\sharp}_{5,11,22} \overline{\lambda^{\sharp}_{4,11,22}} \overline{\lambda^{\sharp}_{4,14,21}} \overline{\lambda^{\sharp}_{5,11,20}},$$
$$u_{22,5} = - \lambda^{\sharp}_{3,14,21} \lambda^{\sharp}_{4,14,22} \lambda^{\sharp}_{6,11,22} \overline{\lambda^{\sharp}_{4,11,22}} \overline{\lambda^{\sharp}_{4,14,21}} \overline{\lambda^{\sharp}_{6,11,20}},$$
$$u_{13,22} = - \lambda^{\sharp}_{4,11,22} \lambda^{\sharp}_{4,14,21} \lambda^{\sharp}_{6,11,20} \overline{\lambda^{\sharp}_{3,14,21}} \overline{\lambda^{\sharp}_{4,14,22}} \overline{\lambda^{\sharp}_{5,13,22}} \overline{\lambda^{\sharp}_{6,11,22}},$$
$$u_{14,22} = - \lambda^{\sharp}_{4,11,22} \lambda^{\sharp}_{4,14,21} \lambda^{\sharp}_{6,11,20} \overline{\lambda^{\sharp}_{3,14,21}} \overline{\lambda^{\sharp}_{4,14,22}} \overline{\lambda^{\sharp}_{5,14,22}} \overline{\lambda^{\sharp}_{6,11,22}},$$
$$u_{6,13} = - \lambda^{\sharp}_{3,14,21} \lambda^{\sharp}_{4,14,22}  \lambda^{\sharp}_{5,13,22} \lambda^{\sharp}_{6,11,22} \overline{\lambda^{\sharp}_{4,11,22}} \overline{\lambda^{\sharp}_{4,14,21}} \overline{\lambda^{\sharp}_{6,11,20}} \overline{\lambda^{\sharp}_{6,13,22}},$$
$$u_{13,20} = \lambda^{\sharp}_{4,11,22} \lambda^{\sharp}_{4,14,21}  \lambda^{\sharp}_{6,11,20} \lambda^{\sharp}_{6,13,22} \overline{\lambda^{\sharp}_{3,14,21}} \overline{\lambda^{\sharp}_{4,14,22}} \overline{\lambda^{\sharp}_{5,13,22}} \overline{\lambda^{\sharp}_{6,11,22}} \overline{\lambda^{\sharp}_{6,13,20}}.$$
The $u_{v_1,v_2}$ involving the vertices 1, 2, 7, 8, 9, 10, 15, 16, 17, 18, 23 and 24 are chosen in the same way as in the proof of uniqueness of the cells for the $\mathcal{A}$ graphs.
\hfill
$\Box$

For $\mathcal{E}^{(24)}$, we have the following representation of the Hecke algebra (we omit those weights which come from the $\mathcal{A}^{(24)}$ graph):
\begin{eqnarray*}
U^{(3,21)} & = & \begin{array}{c} \scriptstyle 12 \\ \scriptstyle 14 \end{array} \left( {\begin{array}{cc}
                \frac{[5]}{[4]} & -\frac{\sqrt{[3][5]}}{[4]} \\
                -\frac{\sqrt{[3][5]}}{[4]} & \frac{[3]}{[4]}
               \end{array} } \right) \;\; = \;\; U^{(12,3)} \quad \textrm{ with labels 21,19} \\
& = & U^{(6,20)} \quad \textrm{ with labels 13,11} \;\; = \;\; U^{(13,6)} \quad \textrm{ with labels 20,22,}
\end{eqnarray*}
\begin{eqnarray*}
U^{(19,12)} & = & \begin{array}{c} \scriptstyle 3 \\ \scriptstyle 4 \end{array} \left( {\begin{array}{cc}
                \frac{1}{[2]} & \frac{\sqrt{[3]}}{[2]} \\
                \frac{\sqrt{[3]}}{[2]} & \frac{[3]}{[2]}
               \end{array} } \right), \;\; = \;\; U^{(21,14)} \quad \textrm{ with labels 3,4} \\
& = & U^{(20,11)} \quad \textrm{ with labels 6,5} \;\; = \;\; U^{(22,13)} \quad \textrm{ with labels 6,5,}
\end{eqnarray*}
\begin{eqnarray*}
U^{(5,19)} & = & \begin{array}{c} \scriptstyle 11 \\ \scriptstyle 14 \end{array} \left( {\begin{array}{cc}
                \frac{[2]}{[3]} & \frac{\sqrt{[2][4]}}{[3]} \\
                \frac{\sqrt{[2][4]}}{[3]} & \frac{[4]}{[3]}
               \end{array} } \right),
\;\; = \;\; U^{(14,5)} \quad \textrm{ with labels 22,19,}
\end{eqnarray*}
\begin{eqnarray*}
U^{(4,22)} & = & \begin{array}{c} \scriptstyle 14 \\ \scriptstyle 11 \end{array} \left( {\begin{array}{cc}
                \frac{[2]}{[3]} & -\frac{\sqrt{[2][4]}}{[3]} \\
                -\frac{\sqrt{[2][4]}}{[3]} & \frac{[4]}{[3]}
               \end{array} } \right),
\;\; = \;\; U^{(11,4)} \quad \textrm{ with labels 19,22,}
\end{eqnarray*}
\begin{eqnarray*}
U^{(20,13)} & = & \begin{array}{c} \scriptstyle 6 \\ \scriptstyle 5 \end{array} \left( {\begin{array}{cc}
                \frac{[5]^2}{[2][9]} & -\frac{[5]\sqrt{[7]}}{[2][9]} \\
                -\frac{[5]\sqrt{[7]}}{[2][9]} & \frac{[7]}{[2][9]}
               \end{array} } \right),
\;\; = \;\; U^{(21,12)} \quad \textrm{ with labels 3,4,}
\end{eqnarray*}
\begin{eqnarray*}
U^{(4,21)} & = & \begin{array}{c} \scriptstyle 12 \\ \scriptstyle 14 \end{array} \left( {\begin{array}{cc}
                \frac{1}{[4]} & \frac{[3]\sqrt{[5]}}{[4]\sqrt{[7]}} \\
                \frac{[3]\sqrt{[5]}}{[4]\sqrt{[7]}} & \frac{[3]^2[5]}{[4][7]}
               \end{array} } \right), \;\; = \;\; U^{(12,4)} \quad \textrm{ with labels 21,19} \\
& = & U^{(5,20)} \quad \textrm{ with labels 13,11} \;\; = \;\; U^{(13,5)} \quad \textrm{ with labels 20,22,}
\end{eqnarray*}
\begin{eqnarray*}
U^{(19,14)} & = & \begin{array}{c} \scriptstyle 3 \\ \scriptstyle 4 \\ \scriptstyle 5 \end{array} \left( {\begin{array}{ccc}
                 \frac{1}{[2]} & \frac{\sqrt{[7]}}{[2][3]} & \frac{\sqrt{[7][10]}}{[3]\sqrt{[2][5]}} \\
                 \frac{\sqrt{[7]}}{[2][3]} & \frac{[7]}{[2][3]^2} & \frac{[7]\sqrt{[10]}}{[3]^2\sqrt{[2][5]}} \\
                 \frac{\sqrt{[7][10]}}{[3]\sqrt{[2][5]}} & \frac{[7]\sqrt{[10]}}{[3]^2\sqrt{[2][5]}} & \frac{[7][10]}{[3]^2[5]}
               \end{array} } \right),
\end{eqnarray*}
\begin{eqnarray*}
U^{(22,11)} & = & \begin{array}{c} \scriptstyle 6 \\ \scriptstyle 5 \\ \scriptstyle 4 \end{array} \left( {\begin{array}{ccc}
                 \frac{1}{[2]} & \frac{\sqrt{[7]}}{[2][3]} & -\frac{\sqrt{[7][10]}}{[3]\sqrt{[2][5]}} \\
                 \frac{\sqrt{[7]}}{[2][3]} & \frac{[7]}{[2][3]^2} & -\frac{[7]\sqrt{[10]}}{[3]^2\sqrt{[2][5]}} \\
                 -\frac{\sqrt{[7][10]}}{[3]\sqrt{[2][5]}} & -\frac{[7]\sqrt{[10]}}{[3]^2\sqrt{[2][5]}} & \frac{[7][10]}{[3]^2[5]}
               \end{array} } \right),
\end{eqnarray*}
\begin{eqnarray*}
U^{(19,11)} & = & \begin{array}{c} \scriptstyle 2 \\ \scriptstyle 4 \\ \scriptstyle 5 \end{array} \left( {\begin{array}{ccc}
                 \frac{[4]}{[5]} & \frac{[4]\sqrt{[7]}}{[3][5]} & \frac{[4]\sqrt{[7]}}{[3][5]} \\
                 \frac{[4]\sqrt{[7]}}{[3][5]} & \frac{[4][7]}{[3]^2[5]} & \frac{[4][7]}{[3]^2[5]} \\
                 \frac{[4]\sqrt{[7]}}{[3][5]} & \frac{[4][7]}{[3]^2[5]} & \frac{[4][7]}{[3]^2[5]}
               \end{array} } \right)
\; = \; U^{(22,14)} \quad \textrm{ with labels 7,4,5,}
\end{eqnarray*}
\begin{eqnarray*}
U^{(3,19)} & = & \begin{array}{c} \scriptstyle 10 \\ \scriptstyle 14 \\ \scriptstyle 12 \end{array} \left( {\begin{array}{ccc}
                \frac{[4]}{[5]} & \frac{\sqrt{[3]}}{\sqrt{[5]}} & \frac{\sqrt{[9]}}{[5]} \\
                \frac{\sqrt{[3]}}{\sqrt{[5]}} & \frac{[3]}{[4]} & \frac{\sqrt{[3][9]}}{[4]\sqrt{[5]}} \\
                \frac{\sqrt{[9]}}{[5]} & \frac{\sqrt{[3][9]}}{[4]\sqrt{[5]}} & \frac{[9]}{[4][5]}
               \end{array} } \right) \; = \; U^{(14,3)} \quad \textrm{ with labels 23,19,21} \\
& = & U^{(6,22)} \quad \textrm{ with labels 15,11,13} \;\; = \;\; U^{(11,6)} \quad \textrm{ with labels 18,22,20,}
\end{eqnarray*}
\begin{eqnarray*}
U^{(4,19)} & = & \begin{array}{c} \scriptstyle 11 \\ \scriptstyle 14 \\ \scriptstyle 12 \end{array} \left( {\begin{array}{ccc}
                 \frac{[2]}{[3]} & \frac{\sqrt{[2][5]}}{[3]\sqrt{[4]}} & \frac{\sqrt{[2][9]}}{\sqrt{[4][7]}} \\
                 \frac{\sqrt{[2][5]}}{[3]\sqrt{[4]}} & \frac{[5]}{[3][4]} & \frac{\sqrt{[5][9]}}{[4]\sqrt{[7]}} \\
                 \frac{\sqrt{[2][9]}}{\sqrt{[4][7]}} & \frac{\sqrt{[5][9]}}{[4]\sqrt{[7]}} & \frac{[3][9]}{[4][7]}
               \end{array} } \right) \\
& = & U^{(14,4)} \quad \textrm{ with labels 22,19,21} \;\; = \;\; U^{(5,22)} \quad \textrm{ with labels 14,11,13}\\
& = & U^{(11,5)} \quad \textrm{ with labels 19,22,20.}
\end{eqnarray*}

The Hecke representation given above cannot be equivalent to that given by Sochen in \cite{sochen:1991} for $\mathcal{E}^{(24)}$ as our weights $[U^{(14,4)}]_{19,19}$, $[U^{(14,4)}]_{21,21}$, $[U^{(11,5)}]_{20,20}$, $[U^{(11,5)}]_{22,22}$ and $[U^{(19,11)}]_{2,2}$ (as well as the corresponding weights under the reflection of the graph which sends vertices $1 \leftrightarrow 8$) have different absolute values to those given by Sochen (and there are no double edges on the graph). We do not believe that there exists two inequivalent solutions for the Hecke representation for $\mathcal{E}^{(24)}$, and that the differences must be due to typographical errors in \cite{sochen:1991}. \\ \\

\paragraph{Acknowledgements}

This paper is based on work in \cite{pugh:2008}. The first author was partially supported by the EU-NCG network in Non-Commutative Geometry MRTN-CT-2006-031962, and the second author was supported by a scholarship from the School of Mathematics, Cardiff University.

\end{document}